\documentclass[10pt,leqno]{article}

\usepackage{latexsym} 
\usepackage[all]{xy} 
\usepackage{amsmath}
\usepackage{amssymb} 
\usepackage{amsfonts} 
\usepackage{color}
\usepackage{theorem}

\newcommand {\OO}{\mathcal{O}} 

\newcommand {\DR}{\mathbf{DR}} 
\newcommand{\Pol}{\mathbf{Pol}} 
\newcommand{\pol}{\mathsf{Pol}}

\newcommand{\B}{\mathcal{P}} 
 
\newcommand{\C}{\mathcal{M}} 
\newcommand{\T}{\mathbf{T}} 
\newcommand{\ti}{\mathsf{T}} 
\newcommand{\Spec}{\mathbf{Spec}} 
\newcommand{\dg}{\mathbf{dg}} 
\newcommand{\dglie}{\mathbf{dgLie}}

\newcommand{\edg}{\epsilon-\mathbf{dg}}
\newcommand{\medg}{\epsilon-\mathbf{dg}^{gr}}
\newcommand{\eMgr}{\epsilon-\mathbf{M}^{gr}}

\newcommand{\mecdga}{\epsilon-\mathbf{cdga}^{gr}}

\newcommand{\cdga}{\mathbf{cdga}}

\newcommand{\dAff}{\mathbf{dAff}}

\newcommand{\dSt}{\mathbf{dSt}}

\newcommand{\s}{\infty}
\newcommand{\D}{\mathbb{D}}
\newcommand{\ealgrM}{\epsilon-\mathbf{CAlg}_{\C}^{gr}}
\newcommand{\ealgrm}{\epsilon-\mathbf{CAlg}_{\mathsf{M}}^{gr}}
\newcommand{\calg}{\mathbf{CAlg}}
\newcommand{\dr}{\mathsf{DR}}
\newcommand{\ind}{\mathbf{Ind}(\C)}
\newcommand{\m}{\mathsf{M}}
\newcommand{\qcoh}{\mathbf{QCoh}}

\newcommand{\red}{\textrm{red}}

\newtheorem{thm}{Theorem}[subsection]
\newtheorem{prop}[thm]{Proposition}

\newtheorem{df}[thm]{Definition}
\newtheorem{cor}[thm]{Corollary}

{\theorembodyfont{\rmfamily}\newtheorem{rmk}[thm]{Remark}}
\newtheorem{ex}[thm]{Example}
\newtheorem{exs}[thm]{Examples}

\begin{document}


\title{\textbf{ Symplectic and Poisson derived geometry\\ and deformation quantization}}

\author{T. Pantev, G. Vezzosi}

\date{}

\maketitle

\begin{abstract}
  We review recent results and ongoing investigation of the
  symplectic and Poisson geometry of derived moduli spaces, and
  describe applications to deformation quantization of such spaces.
\end{abstract}

\tableofcontents

\section*{Introduction} 
\addcontentsline{toc}{section}{Introduction}

From the vantage point of the timeline of the AMS Summer Institutes,
this contribution is a continuation and an update of B. To\"en's 2005
overview \cite{seattle}. Our goal here is to highlight some of the remarkable
developments in derived geometry that we witnessed in the past ten
years.

 Our main topics - symplectic and Poisson geometry on derived moduli
 spaces - are among the latest topics in the area. Even though the
 study of these topics is still in an early stage, it has already led
 to some exciting applications. Among the moduli theoretic
 applications, we might mention the fact that $(-1)$-shifted
 symplectic structures induce symmetric perfect obstruction theories
 in the sense of \cite{bf} (and actually all the geometrically
 interesting examples of symmetric obstruction theories arise this
 way), and the related important result that the Donaldson-Thomas
 moduli space is $(-1)$-shifted symplectic and Zariski locally
 isomorphic to the critical locus of a potential \cite{bbj}.

Section 1 explains and summarizes the main results from \cite{ptvv},
while sections 2 and 3 delve into the substance of \cite{cptvv}. The
aim of our review has been twofold: on one hand to convey the
intuition behind definitions, constructions, and proofs of the main
results, and on the other hand, to explain and motivate the slight
change of point of view going from \cite{ptvv} to \cite{cptvv}. The
upshot is that shifted Poisson geometry and its applications to
deformation quantization of derived categories require a new broader
perspective and new technical tools, i.e.  differential calculus in an
extremely general setting (Section~2.1), and formal localization
(Section~2.2).  Even though these powerful tools were created in order
to solve our specific problems related to shifted Poisson structures,
they constitute also a conceptual advance, very likely to become
relevant in other contexts and to different problems in derived
algebraic geometry.

Note that there is a parallel theory of \emph{shifted quadratic forms}
on derived moduli spaces, but we will not review it here (see
\cite{cliff} for a first investigation).

We now describe, section by section, the mathematical contents of this
paper in more details.

\

\noindent \textbf{Shifted symplectic structures.} A \emph{shifted
  symplectic structure} on a derived stack $X$ with a perfect
cotangent complex $\mathbb{L}_X$ is a structured self-duality of
$\mathbb{L}_X$ up to a shift, i.e. a quasi-isomorphism $\mathbb{T}_X
\simeq \mathbb{L}_X[n]$ induced by a closed $n$-shifted $2$-from on
$X$. The idea is an obvious generalization of the classical definition
of symplectic form, but with an important additional, purely derived
algebro-geometrical feature: for a shifted form on $X$, being closed
is \emph{not} a property but rather a datum.  In other words,
there is a canonical map, called the underlying-form map, from the
space $\mathcal{A}^{p, \textrm{cl}}(X,n)$ of $n$-shifted closed
$p$-forms on $X$, to the space $\mathcal{A}^{p}(X,n)$ of $n$-shifted
$p$-forms on $X$, but this map is not, in general, ``injective'' in
any reasonable sense of the word (e.g. not injective on the connected
components of these spaces). The space $\mathsf{Sympl}(X,n)$ is
exactly the subspace of $\mathcal{A}^{2, \textrm{cl}}(X,n)$ of closed
$2$-forms whose underlying $2$-from is non-degenerate, i.e. such that
the induced map $\mathbb{T}_X \to \mathbb{L}_X[n]$ is a
quasi-isomorphism. Shifted symplectic structures abound, in the sense
that many moduli spaces of interest to algebraic geometers and
topologists, such as the moduli spaces of principal bundles or perfect
complexes on algebraic Calabi-Yau manifolds or compact orientable
topological manifolds, have derived enhancements carrying natural
shifted symplectic structures. In Section 1 we give three general
existence results for shifted symplectic structures on derived moduli
stacks, leading to a long list of examples.

\

\noindent \textbf{Shifted Poisson structures.} Having at our disposal a
theory of shifted symplectic structures, it is natural to look for a more
general theory of \emph{shifted Poisson structures} on derived moduli
stacks. Actually, our original motivation for such a general theory 
came from the expected link between a shifted Poisson structure on a
derived stack and an induced deformation quantization of its dg-derived
category of perfect complexes. We will say more about this
motivation-application below. While classically,  setting up a theory of
Poisson varieties does not present more difficulties than setting up a
theory of symplectic varieties, in derived algebraic geometry the situation
is radically different. The usual bad functoriality properties of shifted
polyvectors (as opposed to the good functoriality of closed shifted forms)
together with the very delicate and intricate strictification problems related to
establishing a meaningful shifted Poisson algebra structure on them,
immediately made us realize that, outside the derived Deligne-Mumford case,
a full-fledged theory of shifted Poisson structures on derived Artin stacks
required new ideas and tools. If $X$ is a derived Artin stack, locally of
finite presentation over the base $k$ (always assumed to be a Noetherian
commutative $\mathbb{Q}$-algebra), then its cotangent complex
$\mathbb{L}_X$ is perfect, and we may consider the graded commutative
differential graded algebra
$$
Pol(X, m) := \bigoplus_{p\geq 0} \Gamma(X, Sym^{p}(\mathbb{T}_{X}[-m])
$$
of $m$-shifted polyvectors on $X$. Here, $m \in \mathbb{Z}$,
$\mathbb{T}_{X}$ is the tangent complex of $X$, i.e. the
$\mathcal{O}_X$-dual of $\mathbb{L}_X$, $\Gamma$ denotes the derived
functor of global sections (i.e. hypercohomology), and the external
grading, called \emph{weight grading}, is given by $p$. In order to define
shifted Poisson structures on $X$, we have to endow $Pol(X, m)$ with
a degree $m$ and weight $-1$ Lie bracket, making it into a graded 
 $\mathbb{P}_{m+1}$-commutative differential graded algebra over $k$. In
particular,  $Pol(X, m)[m]$ will be a graded dg-Lie algebra over
$k$, with a weight $-1$ Lie bracket. Then we could adopt the following
derived variation of a classical definition of Poisson structure, and put
$$
\mathsf{Poiss}(X,n):=\mathsf{Map}_{\mathsf{dglie}^{gr}_k}(k(2)[-1],
Pol(X,n+1)[n+1])
$$ for the space $\mathsf{Poiss}(X,n)$ of $n$-shifted Poisson
structures on $X$, for $n\in \mathbb{Z}$. Here $k(2)[-1]$ is the
graded dg-Lie algebra consisting of $k$ in pure cohomological degree
1, pure weight 2, and trivial bracket, obviously. If $X$ is a smooth
underived scheme, $n=0$, and we replace the mapping space
$\mathsf{Map}_{\mathsf{dglie}^{gr}_k}$ in the model category
$\mathsf{dglie}^{gr}_k$, with its Hom set
$Hom_{\mathsf{dglie}^{gr}_k}$, then we obtain the set of bivectors
$\pi$ on $X$, whose Schouten-Nijenhuis self bracket $[\pi, \pi]$ is
zero, i.e.  exactly the set of Poisson bivectors on $X$. The
functoriality problems mentioned above prevent any elementary, easy
extension of (a shifted version of) the usual Schouten-Nijenhuis
bracket when $X$ is a general derived Artin stack, locally of finite
presentation over $k$. Hence, there is no elementary, easy way of
making sense of the above definition of $\mathsf{Poiss}(X,n)$.

Our solution to this problem consists of two steps. First of all,  in Section
2.1 we build a very general theory of \emph{differential calculus},
including de Rham algebras and polyvectors, in an arbitrary symmetric
monoidal model $\s$-category $\C$ enriched over $k$-dg modules (and
satisfying suitable, mild hypotheses). In particular, for any
 commutative algebra $A$ in $\C$, we are able to make sense of the space
$\mathsf{Sympl}(A, n)$ of $n$-shifted symplectic structures on $A$, to
define a  $\mathbb{P}_{n+1}$-commutative differential graded algebra
$\mathsf{Pol}(A,n)$ over $k$ of $n$-shifted polyvectors on $A$, and hence
to make sense  of the space $\mathsf{Poiss}(A,n)$ of $n$-shifted Poisson
structures on $A$, as explained above. Moreover, we produce a general
comparison map from the the space $\mathsf{Poiss}^{\textrm{nd}}(A,n)$, of
suitably defined non-degenerate $n$-shifted Poisson structures on $A$ to
$\mathsf{Sympl}(A, n)$. The second step is what we call \emph{formal
localization} (Section~2.2), and it concerns essentially the study of the map
$X \to X_{DR}$, for $X$ derived Artin stack, locally of finite presentation
over $k$. Here $X_{DR}$ is the de Rham stack of $X$ (Definition \ref{red&dR}),
and the fiber of $X \to X_{DR}$ at a closed point $\overline{x}: \Spec \,
\mathbb{K} \to X_{DR}$ is the formal completion $\widehat{X_{x}}$ of $X$ at
the corresponding point $x \in X$ \footnote{Note that $X$ and $X_{DR}$ have
the same reduced points.}; hence $X \to X_{DR}$ can be viewed as the
\emph{family of formal completions} of $X$.

The remarkable properties of the map $X \to X_{DR}$, allow us to
associate to any derived Artin stack $X$, locally of finite
presentation over $k$, a commutative algebra $\B_{X}(\infty)$ in a
suitable symmetric monoidal model $\s$-category $\C_{X}$ constructed
out of $X$, such that:
\begin{itemize}
\item There is an equivalence of spaces $\mathsf{Sympl}(\B_X(\s), n) \simeq
\mathsf{Symp}(X, n)$.
\item After forgetting the bracket $\mathsf{Pol}(\B_{X}(\s),n)$ is
equivalent to $Pol(X,n)$, in other words we finally have a way of endowing
$Pol(X,n)$ with the structure of a graded  $\mathbb{P}_{n+1}$-commutative
differential graded algebra over $k$. In particular, it now makes sense to
define $\mathsf{Poiss}(X,n):=\mathsf{Map}_{\mathsf{dglie}^{gr}_k}(k(2)[-1],
\mathsf{Pol}(\B_{X}(\s),n+1)[n+1])$.
\item The canonical map
$\mathsf{Poiss}^{\textrm{nd}}(X,n)=\mathsf{Poiss}^{\textrm{nd}}(\B_{X}(\s),n)
\to \mathsf{Sympl}(\B_{X}(\s), n) \simeq \mathsf{Symp}(X, n)$ is an
equivalence of spaces.
\item The $\s$-category $\mathbf{Perf}(X)$ of perfect complexes on $X$ is
equivalent to a suitably defined $\s$-category of perfect
$\B_{X}(\s)$-dg-modules.
\end{itemize}

\

Let us briefly describe the commutative algebra object
$\B_{X}(\infty)$ in $\C_{X}$. As already observed, the canonical map $
X \longrightarrow X_{DR}$ realizes $X$ as a family of formal derived
stacks over $X_{DR}$, namely as the family of formal completions at
closed points of $X$. By \cite{lu} each of these formal completions is
determined by a dg-Lie algebra $\ell_x$. The collection of dg-Lie
algebras $\ell_x$ \emph{does not} fit together globally in a sheaf of
dg-Lie algebras over $X_{DR}$, simply because its underlying complex
is the shifted tangent complex $\mathbb{T}_{X}[-1]$ of $X$ (see \cite{hen}), which in
general does not admit a flat connection and thus does not descend to
$X_{DR}$. However, a remarkable consequence of derived formal
localization is that the Chevalley-Eilenberg complexes of $\ell_x$,
viewed as graded mixed commutative dg-algebras, \emph{do fit
  together} into a global object over $X_{DR}$. Up to a twist (by
$k(\s)$, see Section 2.1), this is exactly $\B_{X}(\s)$. Thus, formal
localization tells us how to express global geometric objects on $X$
as correspondingly sheafified objects on $X_{DR}$ related to
$\B_X(\s)$.

\

\noindent \textbf{Deformation quantization of derived categories.} One
of our main original motivations for developing a theory of shifted
symplectic and Poisson structures on derived moduli spaces was in fact
a prospective application to \emph{deformation quantization of derived
  categories} of perfect complexes. We are now able to obtain such
applications, and we will briefly describe them here (for more
details, see Section~3). We start by defining the deformation
quantization problem for $n$-shifted Poisson structures, whenever
$n\geq 0$. For every such $n$, we consider a
$\mathbb{G}_m$-equivariant $\mathbb{A}^1_k$-family of $k$-dg-operads
$\mathbb{BD}_{n+1}$ such that its $0$-fiber is the Poisson operad
$\mathbb{P}_{n+1}$ and its generic fiber is the $k$-dg-operad
$\mathbb{E}_{n+1}$ of chains of the topological operad of little
$(n+1)$-disks.  The general deformation quantization problem can then
be loosely stated as follows:

\

\noindent
\textbf{Deformation Quantization Problem. } \emph{Given a
$\mathbb{P}_{n+1}$-algebra stucture on an object in a $k$-linear symmetric monoidal $\s$-categorie,  does it exist a family
of
$\mathbb{BD}_{n+1}$-algebra structures such that its $0$-fiber is the
original $\mathbb{P}_{n+1}$-algebra structure ?}

\

To be more precise, let now $X$ be a derived Artin stack locally of
finite presentation over $k$, and equipped with an $n$-shifted Poisson
structure.  Using the formality of the $\mathbb{E}_{n}$-operad, and
the fact that for $n \geq 1$ the homology operad of $\mathbb{E}_{n+1}$
is $\mathbb{P}_{n+1}$, we can solve the deformation quantization
problem above for the $\mathbb{P}_{n+1}$-algebra structure on
$\B_X(\s)$.  This gives us, in particular, a
$\mathbb{G}_m$-equivariant $1$-parameter family of
$\mathbb{E}_{n+1}$-algebra structures on $\B_X(\s)$.

One of the main results of formal localization (Section 2.2) tells us
that the $\s$-category $\mathbf{Perf}(X)$ of perfect complexes on $X$ is equivalent to
the $\s$-category of (suitably defined) perfect $\B_X(\s)$-modules (in
$\C_{X}$).  We thus get a $1$-parameter
deformation of $\mathbf{Perf}(X)$ as an $\mathbb{E}_n$-monoidal
$\infty$-category, which we call the $n$\emph{-quantization of}
$X$. We also give a version of this result for $n<0$ (where of course
$\mathbb{E}_n$ will be replaced by $\mathbb{E}_{-n}$). In
contrast, the unshifted $n=0$ case for derived Artin stacks, which was
previously addressed for smooth varieties by \cite{ko, yekutieli}, is
not currently covered by our analysis and seems to require new
ideas.

Finally, in Section 3.2, we describe some examples of these
$n$-shifted quantizations, especially the quantization on a formal
neighborhood of a point, and of various derived moduli stacks of
$G$-local systems, for $G$ a complex reductive group. Many more
examples are awaiting a careful investigation.

\

\noindent \textbf{Acknowledgements.} First of all, we would like to
thank our co-authors D. Calaque, B. To\"en, and M. Vaqui\'e for the interesting mathematics we did together. 
We thank V. Melani and M. Porta for their
interesting questions that have hopefully led to a clearer text. We
are grateful to P. Safronov and N.  Rozenblyum for useful exchanges on
various topics treated in this review.  We are also grateful to the
organizers of the 2015 Summer Research Institute on Algebraic Geometry
in Salt Lake City, for their invitation to give our talks and to write
this review.

Tony Pantev was partially supported by NSF research grant DMS-1302242
and by grant \# 347070 from the Simons Foundation.  Gabriele Vezzosi is
a member of the GNSAGA-INDAM group (Italy) and of PRIN-Geometria delle
variet\'{a} algebriche (Italy). He would like to point out the serious
problems of scientific research in public universities in Italy, due
to the substantial lack of acknowledgement from our government
of the important cultural and social role of public research in a modern country, and to
the subsequent largely insufficient investment of government funds
into public research and universities.

\

\bigskip

\noindent \textbf{Background.} We will assume the reader has some
familiarity with derived algebraic geometry, for which useful reviews
are \cite{seattle}, and the more recent \cite{toenems}, while the
foundational works are  To\"en-Vezzosi's \cite{hagII}, J. Lurie's DAG series \cite{lu2}, and also
the recent \cite{luspec}, the last two being available at
http://www.math.harvard.edu/$\sim$lurie/ .  We will use both the
``old'' but sometimes still useful language and theory of model
categories (see e.g. \cite{hov, hir}), and the modern language and
theory of $\s$-categories (\cite{lutop, lualg}).

\

\smallskip

\noindent \textbf{Notations.}

\noindent
$\bullet$ Throughout this paper $k$ will denote a noetherian
commutative $\mathbb{Q}$-algebra.

\noindent
$\bullet$ We will use $(\infty, 1)$-categories \cite{lutop} as our
model for $\infty$-categories. They will be simply called
$\infty$-categories.

\noindent
$\bullet$ As a general rule, a model category is written in sans-serif
fonts $\mathsf{N}$, and we denote in bold-face fonts
$\mathbf{N}:=L(\mathsf{N})$ the $\infty$-category defined as the
homotopy coherent nerve of the Dwyer-Kan localization of
fibrant-cofibrant objects in $\mathsf{N}$ along its weak equivalences,
with the notable exceptions of the $\infty$-category of spaces,
denoted as $\T:= L(\mathsf{sSets})$, and of our base $\s$-category
$\C:= L(\mathsf{M})$ (Section 2). The passage from a model category to
the associated $\s$-category is a localization, and thus very similar
to the passage from the category of complexes in an abelian category
to the associated derived category. This is a good example to keep in
mind.

\noindent
$\bullet$ All symmetric monoidal categories we use will be
symmetric monoidal (bi)closed categories.

\noindent
$\bullet$ $\mathsf{dg}_k$
will denote the symmetric monoidal model category of (unbounded)
complexes of $k$-modules, with fibrations being degreewise surjective
amps, and weak equivalences being quasi-isomorphisms. The associated
$\s$-category will be denoted by $\dg_k$. Note that $\dg_k$ is then a
stable symmetric monoidal $\infty$-category (\cite[Definition
  2.0.0.7]{lualg}).

\noindent
$\bullet$ $\cdga_{k}$ will denote the $\s$-category
of
non-positively graded differential graded $k$-algebras (with
differential increasing the degree). Its objects will be frequently
called simply cdga's. For $A\in \cdga_{k}$, we will write $\pi_{i}\, A
:= H^{-i}(A)$ for any $i \geq 0$.

\noindent
$\bullet$ For $A \in \cdga_{k}$, we will denote either by
$\mathbf{L}(A)$ or $\mathbf{QCoh}(A)$ the $\s$-category of
$A$-dg-modules.

\noindent
$\bullet$ For $A \in \cdga_{k}$, we will denote by $\mathbf{Perf}(A)$
the full sub-$\s$-category of $\mathbf{QCoh}(A)$ consisting of perfect
$A$-dg-modules.

\noindent
$\bullet$ If $X$ is a derived geometric stack, we will denote by
$\mathbf{QCoh}(X)$ the $k$-linear symmetric monoidal dg-category of
quasi-coherent complexes on $X$.

\noindent
$\bullet$ If $X$ is a derived geometric stack, we will
denote by $\mathbf{Perf}(X)$
the symmetric monoidal sub-dg-category of $\mathbf{QCoh}(X)$
consisting of dualizable objects, i.e. perfect complexes over $X$.

\noindent
$\bullet$ If $X$ is a derived geometric
stack, we will denote by $\mathbf{Coh}(X)$ or
 the full sub-dg category of
$\mathbf{QCoh}(X)$ consisting of complexes whose cohomology sheaves
 are coherent over the truncation $\mathsf{t}_{0} X$.

\noindent
 $\bullet$ For a
morphism $A \to B$ of cdga's, the relative cotangent complex will be
denoted $\mathbb{L}_{B/A} \in \mathbf{L}(B)$. When $A=k$, we will
simply write $\mathbb{L}_B$ instead of $\mathbb{L}_{B/k}$.

\noindent
$\bullet$
For derived stacks, we follow the vocabulary of \cite{hagII}. In particular
derived Artin stacks $X$ will be higher derived stacks, unless stated
otherwise, and always have a cotangent complex, denoted as $\mathbb{L}_X \in \mathbf{QCoh}(X)$.
The acronym lfp means, as usual, locally finitely presented.

\noindent
$\bullet$
For a derived
stack $X$, $\Gamma(X, -)$ will always denote the \emph{derived}
functor of global sections on $X$ (i.e. hypercohomology).


\section{Shifted symplectic structures}\label{symplectic}

\subsection{Definitions}\label{DEF}
Let $\epsilon-\mathsf{dg}_{k}^{gr}$ be the category of \emph{graded
  mixed complexes of $k$-dg-modules}. Its objects are
$\mathbb{Z}$-families of $k$-dg-modules $\{E(p)\}_{p\in \mathbb{Z}}$,
equipped with dg-module maps $ \epsilon : E(p) \longrightarrow
E(p+1)[1], $ such that $\epsilon^2=0$, and the morphisms are
$\mathbb{Z}$-families of morphisms in $\mathsf{dg}_{k}$ commuting with
$\epsilon$.  This is a symmetric monoidal model category: weak
equivalences and cofibrations are defined weight-wise (i.e. with
respect to the external $\mathbb{Z}$-grading, that will be called the
\emph{weight} grading), the monoidal structure is defined by
$(E\otimes E')(p):=\bigoplus_{i+j=p}E(i)\otimes E'(j) $, and the
symmetry constraint \emph{does not} involve signs, but just swaps the
two factors in $E(i)\otimes E'(j)$.  Since our base ring $k$ has
characteristic zero, the category
$\mathsf{Comm}(\epsilon-\mathsf{dg}_{k}^{gr})=:\epsilon-\mathsf{cdga}_{k}^{gr}$
of commutative monoid objects in $\epsilon-\mathsf{dg}_{k}^{gr}$ is
again model category, with weak equivalences and fibrations inherited
via the forgetful functor to $\epsilon-\mathsf{dg}_{k}^{gr}$ (which is
then a right Quillen adjoint). According to our general conventions,
we will denote by $\epsilon-\mathbf{dg}_{k}^{gr}$ (respectively,
$\epsilon-\mathbf{cdga}_{k}^{gr}$), the $\infty$-category associated
to $\epsilon-\mathsf{dg}_{k}^{gr}$ (respectively to
$\epsilon-\mathsf{cdga}_{k}^{gr}$). Informally speaking,
$\epsilon-\mathbf{cdga}_{k}^{gr}$ is therefore the $\infty$-category
of $\{B(p) \in \mathbf{dg}_{k} \}_{p\in \mathbb{Z}}$ together with
mixed differential $\epsilon: B(p) \to B(p+1)[1],$ $\epsilon^{2}=0$,
and maps $B(p)\otimes B(q) \to B(p+q)$ which are unital, associative,
commutative, and suitably compatible with $\epsilon$.

The $\infty$-functor $\epsilon-\mathbf{cdga}^{gr}_{k} \to
\mathbf{cdga}_{ k } : \{B(p)\} \mapsto B(0)$ is accessible and
preserves limits, thus (\cite[Corollary~5.5.2.9]{lutop}) has a left adjoint
$\mathbf{DR}: \mathbf{cdga}_{ k }\to \epsilon-\mathbf{cdga}^{gr}_{k}$.

\begin{df}\label{d1.1}
The functor $\mathbf{DR}: \mathbf{cdga}_{ k }\to
\epsilon-\mathbf{cdga}^{gr}_{k}$ is called the \emph{de Rham algebra}
$\infty$-functor.
\end{df}

\begin{rmk}\label{relDR} If $A\in \cdga_k$, we can replace in the previous argument
  $\s$-category $\cdga_k$ with $A/\cdga_k$, and the $\s$-category
  $\epsilon-\mathbf{cdga}^{gr}_{k}$ with $\DR(A)/
  \epsilon-\mathbf{cdga}^{gr}_{k}$, and get a \emph{relative de Rham
    algebra} $\s$-functor $\DR(-/A)$.
\end{rmk}

One can prove that $\mathbf{DR}(A)\simeq Sym_{A}(\mathbb{L}_{A}[-1])$
in $\mathbf{cdga}^{gr}_{k}$ (i.e. as graded cdga's, by forgetting the
mixed differential defined on the rhs). In other words, the
construction $\mathbf{DR}$ yields the full derived de Rham complex of
$A$, including the de Rham differential. We are now able to define
\emph{shifted closed forms} on cdga's. For $m, n \in \mathbb{Z}$,
$k(m)[n]$ will denote the graded $k$-dg-module sitting in weight
degree $m$ and in cohomological degree $-n$.

\begin{df}\label{closed} Let $A \in \mathbf{cdga}_{ k }$
\begin{itemize}
\item The space of \emph{closed $n$-shifted $p$-forms} on $A$ is
  $\mathcal{A}^{p,cl}(A,n)
  :=\mathsf{Map}_{\epsilon-\mathbf{dg}_{k}^{gr}}(k(p)[-p-n],\mathbf{DR}(A))
  \in \mathbb{T}$. An element in $\pi_0 (\mathcal{A}^{p,cl}(A,n))$ is
  called a \emph{closed $n$-shifted $p$-form} on $A$.
\item The space of \emph{$n$-shifted $p$-forms} on $A$ is
  $\mathcal{A}^{p}(A,n):=\mathsf{Map}_{\mathbf{dg}_{k}}(k[-n],
  \wedge_{A}^{p}\mathbb{L}_{A}) \in \mathbb{T}$. An element in $\pi_0
  (\mathcal{A}^{p}(A,n))$ is called a \emph{$n$-shifted $p$-form} on
  $A$.
\item The induced map $u: \mathcal{A}^{p,cl}(A,n) \to
  \mathcal{A}^{p}(A,n)$ is called the \emph{underlying $p$-form map}.
\end{itemize}
\end{df}

\begin{rmk} Here is a more concrete description of the space of shifted (closed) forms.
  If $A \in \mathsf{cdga}_{k}$, and $A'\to A$ is a cofibrant
  replacement in $\mathsf{cdga}_{k}$, then $\oplus_{p\geq 0}
  \mathbb{L}_{A/k}^p =\oplus_{p\geq 0} \Omega_{A'/k}^p $ is a fourth
  quadrant bicomplex with vertical differential $d^v$ induced by
  $d_{A'}$, and horizontal differential $d^h$ given by the de Rham
  differential
  $$
  d^v:\Omega_{A'/k}^{p,i} \to \Omega_{A'/k}^{p,i+1} \,\,, \,\,
  d^h=d_{DR}:\Omega_{A'/k}^{p,i} \to \Omega_{A'/k}^{p+1,i}.
  $$
The Hodge filtration $F^{\bullet}$ defined  by $F^q(A):=
\oplus_{p\geq q}\Omega_{A'/k}^p$  is still a fourth quadrant bicomplex, and if we put $ \underline{\mathcal{A}}^{p,cl}(A, n) :=\mathrm{Tot}^{\prod}(F^p(A)[n+p] $, 
we have
$$
\mathcal{A}^{p,cl}(A, n)=|\underline{\mathcal{A}}^{p,cl}(A, n) | \,\,\, n \in \mathbb{Z}
$$
where $|E|$ denotes $\mathsf{Map}_{\mathbf{dg}_{k}}(k,E)$ i.e. the
Dold-Kan construction applied to the $\leq 0$-truncation of the
dg-module $E$, and $\mathrm{Tot}^{\prod}$ is the totalization by
products. In particular, we have a corresponding
\emph{Hodge tower} of dg-modules $$... \to \underline{\mathcal{A}}^{p,cl}(A, 0)[-p] \to
\underline{\mathcal{A}}^{p-1,cl}(A, 0)[1-p] \to ... \to \underline{\mathcal{A}}^{0,cl}(A, 0), $$
where, for any $p$, the cofiber of $ \underline{\mathcal{A}}^{p,cl}(A, 0)[-p] \to
\underline{\mathcal{A}}^{p-1,cl}(A, 0)[1-p]$ is equivalent to the dg-module $\underline{\mathcal{A}}^{p-1}(A,0)[1-p]:=  (\wedge_{A}^{p-1}\mathbb{L}_{A})[1-p]$ of $(1-p)$-shifted $(p-1)$-forms on $A$ (so that we have an equivalence $|\underline{\mathcal{A}}^{p-1}(A,0)[1-p] | \simeq \mathcal{A}^{p-1}(A,1-p)$ in $\T$).
Finally, let us observe that the rightmost dg-module $\underline{\mathcal{A}}^{0,cl}(A, 0)$ in the above Hodge tower, is
exactly Illusie's \emph{derived de Rham complex} of $A$
(\cite[ch. VIII]{ill}).
\end{rmk}

\begin{rmk}\label{generalM} Note  that
  the de Rham algebra functor, and hence the notion of (closed)
  shifted forms, makes sense when $\mathbf{dg}_{k}$ is replaced by
  (essentially) any symmetric monoidal stable $k$-linear
  $\infty$-category $\C$. The intermediate categories of interest will
  then be $\epsilon-\C^{gr}$ (generalizing
  $\epsilon-\mathbf{dg}_{k}^{gr}$), and $\ealgrM$ (generalizing
  $\epsilon-\mathsf{cdga}_{k}^{gr}$). For any $A \in
  \mathbf{CAlg}_{\C}$, this will yield a cotangent complex
  $\mathbb{L}^{\C}_{A} \in A-\mathbf{Mod}_{\C}$, a de Rham algebra
  functor $\mathbf{DR}^{\C}: \mathbf{CAlg}_{\C} \to \ealgrM$, and a
  space of $n$-shifted (closed) $p$-form $\mathcal{A}_{\C}^{p}(A,n)$
  ($\mathcal{A}_{\C}^{p,cl}(A,n)$), where the sub/superscript ${\C}$
  indicates that all the constructions are performed
  internally to $\C$. This level of generality and flexibility in the
  choice of the context for our differential calculus, will prove
  extremely useful in the rest of the paper. As relevant cases, the
  reader should keep in mind the case where
  $\C=\epsilon-\mathsf{dg}_{k}^{gr} $ or, more generally, diagrams in
  $\epsilon-\mathsf{dg}_{k}^{gr}$. We will come back to this generalization more systematically in Section \ref{diff}, and use it as an
  essential tool starting from Section \ref{sps}.
\end{rmk}

We are now ready to globalize the above construction to derived
stacks. We start by globalizing the de Rham algebra construction
(Definition \ref{d1.1}).  The functor $A \to \DR(A)$, and its relative
version (over a fixed base $B$, see Remark \ref{relDR}), are both
derived stacks (for the \'etale topology) with values in mixed graded
dga's, so we give the following

\begin{df}\label{DRglobal}
  \emph{(1)} Let $F \to \Spec \, B$ be a map in $\dSt_k$. The
  \emph{relative de Rham algebra of $F$ over $B$}
  $$
  \DR(F/B):=\lim_{\Spec\, C \to F} \DR(C/B) \in \mecdga_B
  $$
  (see Remark \ref{relDR}) where the limit is taken in the
  $\s$-category $\mecdga_B= B/\mecdga_k$ of graded mixed
$B$-linear cdgas, and over all morphisms $\Spec\, C \rightarrow F$ of
  derived stacks over $\Spec \, B$.

  \

  \noindent
  \emph{(2)} For an arbitrary map $F \to G$ in $\dSt_k$, we define the
  \emph{relative de Rham algebra of $F$ over $G$}
  as $$\DR(F/G):=\lim_{\Spec\, A \to G} \DR(F_A/ A) \in \mecdga_k \,
  ,$$ where $F_A$ denotes the base change of $F\to G$ along $\Spec\, A
  \to G$, and the limit is taken in the $\s$-category $\mecdga_k$.
\end{df}

\

\smallskip

We now globalize to derived stacks the notion of (closed) shifted
forms. From Definition \ref{closed} one deduces $\infty$-functors
$\mathcal{A}^{p,cl}(-,n): A \mapsto \mathcal{A}^{p,cl}(A,n)$, and
$\mathcal{A}^{p}(-,n)): A \mapsto \mathcal{A}^{p}(A,n)$ from $\cdga_k$
to $\mathbf{T}$. By \cite[Proposition 1.11]{ptvv}, these functors are
derived stacks (for the \'etale topology). This allows us to globalize
Definition \ref{closed} on an arbitrary derived stack.

\begin{df} Let $F$ be a derived stack. 
\begin{itemize}
\item The space of \emph{closed $n$-shifted $p$-forms} on $F$ is
  $\mathcal{A}^{p,cl}(F,n):=\mathsf{Map}_{\mathbf{dSt}_{k}}(F,
  \mathcal{A}^{p,cl}(-,n))$.
\item The space of \emph{$n$-shifted $p$-forms} on $F$ is
  $\mathcal{A}^{p}(F,n):=\mathsf{Map}_{\mathbf{dSt}_{k}}(F,
  \mathcal{A}^{p}(-,n))$.
\item The resulting induced map $u: \mathcal{A}^{p,cl}(F,n) \to
  \mathcal{A}^{p}(F,n)$ is called the \emph{underlying $p$-form map}.
\end{itemize}
\end{df}

\

\noindent
Note that, in general, the homotopy fiber of the underlying $p$-form
map $u$ can be non-trivial (i.e. not empty nor contractible). Hence
being closed is a datum rather than a property, for a general derived
stack.

\

\begin{rmk} (1) Equivalently, we have
  $\mathcal{A}^{p,cl}(F,n) \simeq\lim_{\Spec \, A \in (\dAff/F)^{op}}\mathcal{A}^{p,cl}(A,n)$, and 
$\mathcal{A}^{p}(F,n) \linebreak \simeq \lim_{\Spec \, A \in
    (\dAff/F)^{op}}\mathcal{A}^{p,cl}(A,n)$, where the limits are taken in
  the $\infty$-category of $\infty$-functors from $(\dAff/F)^{op}$ to
  $\mathbf{T}$.

  \

  \noindent
  (2) Also note that we have an equivalence
  $\mathcal{A}^{p,cl}(F,n)\simeq
  \mathsf{Map}_{\epsilon-\mathbf{dg}_{k}^{gr}}(k(p)[-p-n],\DR(F/k))$
  in $\T$.

  \

  \noindent
  (3) For $F= \Spec \, A$, the complex $\wedge^p\mathbb{L}_A$ has
  non-positive tor-amplitide, hence there are no non-trivial
  $n$-shifted $p$-forms on $F$, for $n >0 $. For $F=X$ an underived
  smooth scheme, a similar argument shows that $F$ admits no
  non-trivial $n$-shifted $p$-forms on $F$, for $n <0 $. If moreover
  $X$ is proper over $k$, then the degeneration of the Hodge-to-de
  Rham spectral sequence implies that the underlying $p$-form map is
  injective on $\pi_0$.
\end{rmk}

\

\noindent
For derived (higher) Artin stacks $F$, the space
$\mathcal{A}^{p}(F,n)$ has the following equivalent description
(smooth descent).

\begin{prop}\emph{\cite[Proposition 1.14]{ptvv}}\,
If $F \in \dSt$ is Artin, then we have an
equivalence
$$
\mathcal{A}^{p}(F,n) \simeq \mathsf{Map}_{\qcoh
  (F)}(\mathcal{O}_{F}, \wedge^{p}\mathbb{L}_{F}[n])\, ,
$$
functorial
in $F$.
\end{prop}

In particular $\pi_0 (\mathcal{A}^{p}(F,n)) \simeq \mathbb{H}^{n}(X,
\wedge^p \mathbb{L}_F)$, for $F$ Artin. Thus an $n$-shifted $2$-form
$\omega$ can be identified with a map $\omega : \mathcal{O}_F \to
\wedge^2 \mathbb{L}_F [n]$. If $F$ is moreover locally of finite
presentation over $k$ (so that its cotangent complex $\mathbb{L}_F$ is
perfect, i.e. dualizable in $\qcoh (F)$), we may associate to such an
$\omega$ an adjoint map $\omega^{\flat}: \mathbb{T}_{F} \to
\mathbb{L}_F [n]$, where $\mathbb{T}_F$ denotes the dual of
$\mathbb{L}_F$ in $\qcoh (F)$, and it is called the \emph{tangent
  complex} of $F$. An $n$-shifted $2$ form $\omega$ on such an $F$ is
said to be \emph{non-degenerate} if the map $\omega^{\flat}$ is an
equivalence, i.e. an isomorphism in the derived category of
quasi-coherent complexes on $F$, and we will denote by
$\mathcal{A}_{nd}^{2}(F,n)$ the subspace of $\mathcal{A}^{2}(F,n)$
consisting of connected components of non-degenerate forms.

\begin{df}\label{sympl}
Let $F$ be an derived Artin stack locally of finite presentation over
$k$. The \emph{space of $n$-shifted symplectic structures} on $F$ is
defined by the following pullback diagram in $\mathbb{T}$
$$
\xymatrix{\mathsf{Sympl}(F,n) \ar[r] \ar[d] &
  \mathcal{A}^{2,cl}(F,n) \ar[d] \\
  \mathcal{A}_{nd}^{2}(F,n) \ar[r] & \mathcal{A}^{2}(F,n)},
$$
and an element in $\pi_0 (\mathsf{Sympl}(F;n))$ is called a
\emph{$n$-shifted symplectic structure} on $F$.
\end{df}

\

In other words, an $n$-shifted symplectic structure $\omega$ on $F$ is
a closed $n$-shifted 2-form on $F$ whose underlying 2-form is
non-degenerate; in particular, $\mathbf{L}_F$ is self dual, up to a
shift. We use the word symplectic structure instead of symplectic form
because a shifted symplectic structure is a closed- 2-form, and with
respect to the underlying form, this consists of additional structure
rather than just being a property. 

The non-degeneracy condition entails a mixture of the (higher) stacky
(i.e. positive degrees in the cotangent complex) and derived
(i.e. negative degrees of the cotangent complex) nature of the stack
$F$, and in particular it poses severe restrictions on the existence
of shifted symplectic structures on a given stack. E.g. it is clear
that if $\mathbb{L}_F$ has perfect amplitude in $[a,b]$, then $F$ may
only support $(a+b)$-shifted symplectic structures. More precisely, it
is easy to check that for a smooth underived scheme $X$, not \'etale
over $k$, the space $\mathsf{Sympl}(X,n)$ is empty for $n\neq 0$, and
either empty or contractible for $n=0$, and moreover, the set of
connected components of $\mathsf{Sympl}(X,0)$ is in canonical
bijection with the set of usual algebraic symplectic forms on $X$
(\cite[p. 298]{ptvv}). So, we get nothing new for smooth underived
schemes, or more generally, smooth Deligne-Mumford stacks. However, we
will see in the following subsection that there are plenty of derived
schemes or stacks carrying interesting shifted symplectic forms.

\subsection{Existence theorems and examples}
\label{symplex}
We will now review the three basic existence theorems (Theorems
\ref{MAP}, \ref{Lagr}, \ref{Perf}, below) for shifted symplectic
structures established in \cite{ptvv}. In combination they give a long
list of non-trivial examples of shifted symplectic stacks.

The first interesting example of a shifted symplectic stack is the
classifying stack $\mathsf{B}G$ of a smooth affine reductive group
scheme over $k$. In this case ,we have (\cite[p. 299]{ptvv})
$$
\pi_0
(\mathsf{Sympl}(\mathsf{B}G,n))= \begin{cases} 0 & n\neq 2
  \\ Sym_{k}(\mathfrak{g}^{\vee})_{nd}^{G} & n=2
        \end{cases}
$$
where $\mathfrak{g}$ is the Lie algebra of $G$, and
$Sym_{k}(\mathfrak{g}^{\vee})_{nd}^{G}$ is the set of non-degenerate
$G$-invariant symmetric bilinear forms on $\mathfrak{g}$. At the level
of tangent complex $\mathbb{T}_{\mathsf{B}G,e} \simeq
\mathfrak{g}[1]$, the underlying 2-form corresponding to a
non-degenerate $G$-invariant symmetric bilinear form $\sigma:
Sym_k^2(\mathfrak{g}) \to k$ is given by the
composite
$$
\xymatrix{\mathbb{T}_{\mathsf{B}G,e} \wedge
  \mathbb{T}_{\mathsf{B}G,e} \ar[r]^-{\sim} & \mathfrak{g}[1] \wedge
  \mathfrak{g}[1] \ar[r]^-{\sim} & Sym_k^{2}(\mathfrak{g})[2]
  \ar[r]^-{\sigma[2]} & k[2], }
$$
where the central quasi-isomorphism
is given by d\'ecalage. For example, if $G=\mathrm{GL}_{n}$, the usual
map $(A,B) \mapsto \mathrm{tr}(AB),$ for $A, B$ $(n\times n)$ matrices
over $k$, yields a $2$-shifted symplectic form on
$\mathsf{B}\mathrm{GL}_n$.  This example will be vastly generalized in
Theorem \ref{Perf} below.

As a second example, for any $n\in \mathbb{Z}$, one has that the
$n$-shifted cotangent stack $\mathrm{T}^* F [n]:=$ \linebreak $\Spec
_{F}(Sym_{\mathcal{O}_F}(\mathbb{T}_{F}[-n]))$ of a derived
Deligne-Mumford stack $F$ lfp over $k$, is canonically $n$-shifted
symplectic via the de Rham differential of the canonical shifted
Liouville form (\cite[Proposition 1.21]{ptvv}). Recently, D. Calaque
has extended this result to derived Artin stacks lfp over $k$ \cite{cal2}.

The first general existence result for shifted symplectic form is an
enhanced derived version of the main result in \cite{aksz}.

\begin{thm}\label{MAP} Let $F$ be a derived Artin stack lfp over $k$, equipped with
  a $n$-shifted symplectic form, and let $X$ be an
  $\mathcal{O}$-compact derived stack equipped with a $d$-orientation
  $[X]: \mathbb{H}(X,\mathcal{O}_X) \to k[-d]$. If
  $\mathsf{MAP}_{\dSt}(X,F)$ is a derived Artin stack lfp over $k$,
  then it carries a canonical $(n-d)$-shifted symplectic structure.
\end{thm}

We direct the reader to \cite[2.1]{ptvv} for the definition of
$\mathcal{O}$-compact derived stack, and for the notion of
$d$-orientation on a $\mathcal{O}$-compact derived stack (i.e. for the
special properties of the map $[X]$ in the theorem), and to
\cite[Theorem 2.5]{ptvv} for a detailed proof. Here we will content ourselves with
a few comments.

First of all the class $\mathcal{O}$-compact derived stacks equipped
with a $d$-orientation includes compact smooth oriented topological
manifolds $M$ of dimension $d$ (identified with constant stacks with
value $M$, and where capping with the fundamental class gives the
$d$-orientation), Calabi-Yau varieties of complex dimension $d$ (where
the orientation is given by a trivialization of the canonical sheaf
followed by the trace map), and De Rham stacks $X=Y_{DR}$ for $Y$ a
smooth and proper Deligne-Mumford stack with connected geometric
fibers, and relative dimension $d/2$ over $\Spec \, k$ \footnote{The
  notion of de Rham stack will be defined and discussed in Section
  \ref{FL}} (where the $d$-orientation is induced by the choice of a
fundamental class in de Rham cohomology $\mathbb{H}_{dR}(Y/k,
\mathcal{O})$).

As a second comment we outline  the
proof of Thm \ref{MAP}. One first uses the evaluation map
$\mathsf{ev}: X \times \mathsf{MAP}_{\dSt}(X,F) \to F$ in order to
pull back the $n$-shifted symplectic structure $\omega$ on $F$, to a
closed form on $X \times \mathsf{MAP}_{\dSt}(X,F)$; this pullback is
then ``integrated along the fiber'' of the projection $X \times
\mathsf{MAP}_{\dSt}(X,F) \to \mathsf{MAP}_{\dSt}(X,F)$, and this
integrated form is shown to be $(n-d)$-shifted symplectic. The
possibility of defining an integration along the fiber $X$ is a
consequence of the definition of $d$-orientation on an
$\mathcal{O}$-compact stack (\cite[Definition 2.3]{ptvv}). Finally, we
observe that the general question of finding optimal conditions on $X$
and $F$ ensuring that $\mathsf{MAP}_{\dSt}(X,F)$ is a derived Artin
stack lfp over $k$ is delicate, the Artin-Lurie representability
criterion (even in the simplified version of \cite[Appendix]{hagII})
will give an answer in all the applications we will discuss below.

We come to the second existence theorem for shifted symplectic
structures. Before stating it, we need to define the notion of
\emph{lagrangian structure} on a map whose target is a shifted
symplectic stack. We start by defining what is an \emph{isotropic
  structure} on such a map.

\begin{df}\label{isotr}
Let $X$ and $F$ be derived Artin stacks lfp over $k$, $\omega$ a
$n$-shifted symplectic structure on $F$, and $f:X\to F$ a map.  The
\emph{space} $\mathsf{Isotr}(f; \omega)$ \emph{of isotropic structures
  on the map} $f$ relative to $\omega$ is the space of paths
$\Omega_{f^*\omega, 0}\mathcal{A}^{2, cl}(X, n)$ between $f^*\omega$
and $0$ in the space of $n$-shifted closed $p$-forms on $X$. An
element in $\gamma \in \pi_0 (\mathsf{Isotr}(f; \omega))$, i.e. a path
between $f^*\omega$ and $0$, is called an \emph{isotropic structure
  on} $f$ relative to $\omega$.
\end{df}
 
The idea, ubiquitous in all of derived geometry and higher category
theory, and that we already saw in action in the definition of closed
forms (Def \ref{closed}), is that it is not enough to say that there
exists a path between $f^*\omega$ and $0$ in $\mathcal{A}^{2, cl}(X,
n)$ (this would be a property), but one rather has to specify one  such
path (i.e. a datum).

Once an isotropic structure $\gamma$ is given, it makes sense to say
that it has the property of being non-degenerate, as follows. By
composition with the underlying 2-form map $u: \mathcal{A}^{2, cl}(X,
n) \to \mathcal{A}^{2}(X, n) $, the path $\gamma$ induces a path
$\gamma_{u}$ between $u(f^*\omega)=f^*(u(\omega))$ and $0$ in
$\mathcal{A}^{2}(X, n)$. Since $\mathbb{L}_X$ is perfect, by
adjunction, this yields in turn a path $\gamma_{u}^{\flat}$ between
$(f^*(u(\omega)))^{\flat}$ and $0$ in
$\mathsf{Map}_{\mathbf{Perf}(X)}(\mathbb{T}_X, \mathbb{L}_X [n])$,
where $(f^*(u(\omega)))^{\flat}$ is the
composite $$\xymatrix{\mathbb{T}_X \ar[r]^-{a^{\vee}} &
  f^*\mathbb{T}_F \ar[rr]^-{f^*(u(\omega)^{\flat})} & &
  f^*\mathbb{L}_F [n] \ar[r]^{a[n]} & \mathbb{L}_X[n]},$$ $a: f^*
\mathbb{L}_F \to \mathbb{L}_X$ being the canonical map induced by
$f$. If we denote by $t_{f,\omega}$ the composite $a[n] \circ
f^*(u(\omega)^{\flat})$, we thus obtain that $\gamma_{u}^{\flat}$ is a
homotopy commutativity datum for the square
$$
\xymatrix{\mathbb{T}_X
  \ar[r] \ar[d]_-{a^{\vee}} & 0 \ar[d] \\ f^* \mathbb{T}_F
  \ar[r]_-{t_{f,\omega}} & \mathbb{L}_X [n].}$$ In particular, if we
denote by $\mathbb{T}^{\perp}_{f, \omega}$ the \emph{pullback} in the
diagram $$\xymatrix{\mathbb{T}^{\perp}_{f, \omega} \ar[r] \ar[d] & 0
  \ar[d] \\ f^* \mathbb{T}_F \ar[r]_-{t_{f,\omega}} & \mathbb{L}_X
     [n]}
$$
(i.e. the kernel of $t_{f,\omega}$), we get a canonical
induced map $\theta_{\gamma}: \mathbb{T}_X \to \mathbb{T}^{\perp}_{f,
  \omega}$ in $\mathbf{Perf}(X)$.

\begin{df}\label{lagr} In the setting of Definition \ref{isotr}, an isotropic structure
  $\gamma \in \pi_0 (\mathsf{Isotr}(f; \omega))$ is called
  \emph{non-degenerate} or \emph{lagrangian} if the induced map
  $\theta_{\gamma}: \mathbb{T}_X \to \mathbb{T}^{\perp}_{f, \omega}$
  is an equivalence in $\mathbf{Perf}(X)$ (i.e. an isomorphism in the
  underlying derived/homotopy category). The \emph{space}
  $\mathsf{Lagr}(f; \omega)$ of \emph{lagrangian structures on} $f:
  X\to F$ \emph{relative to} $\omega$ is the subspace of
  $\mathsf{Isotr}(f; \omega)$ consisting of connected components of
  lagrangian structures.
\end{df}

\begin{rmk}\label{rmks} (1) It is easy to check that if $X$ and $F$ are
  underived smooth schemes, $\omega$ is a usual (i.e. $0$-shifted)
  symplectic structure on $F$, and $f$ is a closed immersion, then
  $\mathsf{Lagr}(f; \omega)$ is either empty or contractible, and it
  is in fact contractible iff $X$ is a usual smooth lagrangian
  subscheme of $F$ via $f$. The nondegeneracy condition ensures that
  the dimension of $X$ is then half of the dimension of $F$, and in
  fact $\mathbb{T}^{\perp}_{f, \omega}$ is then quasi-isomorphic to
  the usual symplectic orthogonal of $\mathbb{T}_X= \mathrm{T}_X$ in
  $\mathbb{T}_F= \mathrm{T}_X$. Thus, the notion of lagrangian
  structure reduces to the usual notion of lagrangian subscheme in
  this case.
  
  \

  \noindent
  (2) By rephrasing Definition \ref{lagr}, an isotropic
  structure $\gamma \in \pi_0 (\mathsf{Isotr}(f; \omega))$ is
  lagrangian iff the (homotopy) commutative
  square $$\xymatrix{\mathbb{T}_X \ar[r] \ar[d]_-{a^{\vee}} & 0 \ar[d]
    \\ f^* \mathbb{T}_F \ar[r]_-{t_{f,\omega}} & \mathbb{L}_X [n]}$$
  is actually a pullback square. But the
  square $$\xymatrix{f^*\mathbb{T}_F \ar[r]^-{t_{f,\omega}}
    \ar[d]_-{f^*(u(\omega)^{\flat})} & \mathbb{L}_{X}[n]
    \ar[d]^-{\mathrm{id}} \\ f^*\mathbb{L}_F [n] \ar[r]_-{a[n]} &
    \mathbb{L}_X [n]}$$ is a pullback because $f^*(u(\omega)^{\flat})$
  is an equivalence ($u(\omega)$ being non-degenerate), hence $\gamma$
  is lagrangian iff the outer square in $$\xymatrix{\mathbb{T}_X
    \ar[r] \ar[d]_-{a^{\vee}} & 0 \ar[d] \\ f^* \mathbb{T}_F \ar[d]
    _-{f^*(u(\omega)^{\flat})} \ar[r]^-{t_{f,\omega}} & \mathbb{L}_X
       [n] \ar[d]^-{\mathrm{id}} \\ f^*\mathbb{L}_F [n] \ar[r]_-{a[n]}
       & \mathbb{L}_X [n]}$$ is cartesian, i.e. iff the induced
  canonical map $\rho_{\gamma}: \mathbb{T}_X \to \ker (a[n]) \simeq
  \mathbb{L}_{f}[n-1]$ is an equivalence i.e. iff the shifted dual map
  $\Theta_{\gamma}:= \rho_{\gamma}^{\vee} [n-1]: \mathbb{T}_f \to
  \mathbb{L}_X [n-1]$ is an equivalence. This shows the equivalence
  between Definition \ref{lagr} and \cite[Definition
    2.8]{ptvv}.

  \

  \noindent
  (3) As first noticed by D. Calaque, \emph{shifted
    symplectic structures are particular instances of lagrangian
    structures} (a fact that is obviously false inside usual,
  $0$-shifted and underived algebraic geometry). In fact, let $n \in
  \mathbb{Z}$, $X$ be a derived stack lfp over $k$, $f:X \to \Spec \,
  k$ be the structure map, and let us endow $\Spec \, k$ with the
  unique $(n+1)$-shifted symplectic structure $\omega_{n+1}^{0}$ (note
  that $\mathbb{L}_{k}\simeq \Omega_{k/k}^1[0] \simeq 0[0]$, therefore
  for any $m\in \mathbb{Z}$, there is a unique $m$-shifted symplectic
  form given by the the shift of the zero form). Now, if $\gamma \in
  \pi_0 (\mathsf{Lagr}(f; \omega))$ as in Definition \ref{lagr}, then
  $\gamma$ is actually a loop at $0$ inside $\mathcal{A}^{2, cl}(X,
  n+1)$, and its class $[\gamma] \in \pi_1 (\Omega_{0,0}
  (\mathcal{A}^{2, cl}(X, n+1)))\simeq \pi_0 (\mathcal{A}^{2, cl}(X,
  n))$\footnote{This isomorphism is perhaps more familar to the reader
    in the following form: a self-homotopy $h$ of the zero map of
    complexes $0: E \to F$ is the same thing as map \emph{of
      complexes} $E \to F[1]$.}  gives an $n$-shifted closed 2-form
  $\omega_{\gamma}$ on $X$. The non-degeneracy condition on $\gamma$
  is equivalent to the fact that the induced map $\Theta_{\gamma}:
  \mathbb{T}_f \simeq \mathbb{T}_X \to \mathbb{L}_X [n-1]$ of point
  (2) in this Remark, is an equivalence. But it is easy to check that
  this map is exactly $\omega^{\flat}_{\gamma}$, hence
  $\omega_{\gamma}$ is indeed an $n$-shifted symplectic form on
  $X$. By using again that $\pi_1 (\Omega_{0,0} (\mathcal{A}^{2,
    cl}(X, n+1)))\simeq \pi_0 (\mathcal{A}^{2, cl}(X, n))$, and
  running the previous argument backwards, we get an equivalence of
  spaces $\mathsf{Sympl}(X,n) \simeq \mathsf{Lagr}(f: X\to \Spec \, k
  , \omega_{n+1}^{0})$.
\end{rmk}

The link between shifted symplectic structures and lagrangian
structures expressed in Remark \ref{rmks} (3) extends to the case of
lagrangian intersections as follows.
 
 \begin{thm}\label{Lagr}
 Let $n \in \mathbb{Z}$, $(F, \omega)$ be a $n$-shifted symplectic
 stack, $f_i: X_i \to F$, $i=1,2$ maps between derived Artin stacks
 lfp over $k$, and $\gamma_i$ lagrangian structures on $f_i$ relative
 to $\omega$, $i=1,2$. Then, there is a canonical induced
 $(n-1)$-shifted symplectic structure on the fiber product $X_1
 \times_{F} X_2$.
 \end{thm}

Recall that, according to our conventions, all fiber products of
derived stacks are taken in the $\infty$-category of derived stacks,
and are therefore implicitly derived fiber products.  We refer the
reader to \cite[Theorem 2.9]{ptvv} for a proof of Theorem
\ref{Lagr}. We will only give a sketch of the argument in the
classical case i.e. for $n=0$, $F, X_1,X_2$ underived smooth schemes,
and $f_i$ closed immersions, $i=1,2$, in order to convey the main idea
of why Theorem \ref{Lagr} is true. Under our hypotheses, $X_1$ and
$X_2$ are usual lagrangian smooth subschemes of the smooth symplectic
scheme $F$. If $X_{12}:=X_1 \times_{F} X_2$ denotes the (derived)
intersection of $X_1$ and $X_2$ inside $F$, we may pull back the
closed form $\omega$ to $X_{12}$ in two different ways, i.e. using
$f_1$ or $f_2$, and get two closed forms $\omega_1$ and $\omega_2$ on
$X_{12}$. Now, $X_{12}$ is a derived fiber product, hence these two
pullbacks come equipped with a canonical path $\gamma$ between them
inside $\mathcal{A}^{2,cl}(X_{12},0)$. On the other hand, as $X_1$ and
$X_2$ are lagrangian subschemes we have $\omega_1 =\omega_2 =0$, so
that $\gamma$ is in fact a loop at $0$ in
$\mathcal{A}^{2,cl}(X_{12}, 0)$. Since, $\pi_1 (
\mathcal{A}^{2, cl}(X_{12},0) ; 0)\simeq \pi_0 (\mathcal{A}^{2,
  cl}(X_{12}, -1))$, the class $[\gamma] \in \pi_0 (\mathcal{A}^{2,
  cl}(X_{12}, -1))$ defines a $(-1)$-shifted closed 2-form on
$X_{12}$, whose non-degeneracy follows easily from the same property
for $\omega$. Thus $[\gamma]$ is a $(-1)$-shifted symplectic structure
on the derived intersection $X_{12}$. The appearence of a $(-1)$-shift here
perhaps explains why this phenomenon, even though arising from a
completely classical situation in usual (algebraic) symplectic
geometry, was not observed before.

\begin{rmk}\label{crit}
  The following special case of \ref{Lagr} is particularly relevant
  for applications to Donaldson-Thomas invariants (see Theorem
  \ref{joyce} below). Let $X$ be a smooth underived scheme and $g: X
  \to \mathbb{A}_{k}^1$ a smooth function. We may embed $X$ inside its
  cotangent bundle either via the zero section or via the differential
  $dg: X \to \mathrm{T}^*X$, and both these immersions are lagrangian
  with respect to the standard symplectic structure on
  $\mathrm{T}^*X$. The derived intersection of these two lagrangians
  is called the derived critical locus of $g$, and is denoted by
  $\mathsf{dCrit}(g)$. Note that its truncation is the usual
  scheme-theoretic critical locus of $g$. Now, Theorem \ref{Lagr}
  endow $\mathsf{dCrit}(g)$ with a canonical $(-1)$-shifted symplectic
  structure $\omega_{g}$. One can rather easily give an explicit
  description of $\omega_{g}$ using Koszul resolutions (see,
  e.g. \cite{vez1}), and observe that a similar result holds by
  replacing $dg$ with an arbitrary closed 1-form on $X$. However,
  since derived critical loci are particularly important as local
  models of $(-1)$-shifted symplectic structures, the generalization
  to derived zero schemes of closed 1-forms has not yet received much
  attention. Also notice that the previous construction can be easily
  generalized to build derived zero loci of sections of arbitrary
  locally free sheaves on $X$.
\end{rmk}

\begin{rmk}\label{damienboundary}
  One may combine Theorem \ref{MAP} and Definition \ref{lagr} together
  with a notion of \emph{relative orientation}, in order to get a
  relative version of Theorem \ref{MAP}. This is due to D. Calaque
  (\cite[2.3]{cal}). Here is an outline of the construction. Given a
  map of derived stacks $b: B \to X$, and a perfect complex $E$ on
  $X$, we let $\mathbb{H}(X\textrm{rel}B,E)$ be the fiber of the
  induced map $b_E: \Gamma(X,E)\to \Gamma (B, b^* E)$. We define a
  \emph{relative $d$-orientation on} $b$ is a map $\eta_{b}:
  \mathbb{H}(X\textrm{rel}B,\OO_X) \to k[-d]$ satisfying the following
  two non-degeneracy properties. First of all, we assume that $B$ is
  $\OO$-compact, and we require that the composite map $\Gamma(B,
  \OO_B) \to \mathbb{H}(X\textrm{rel}B,\OO_X)[1] \to k[-d+1]$ defines
  a $(d-1)$-orientation on $B$. Then, for any $E\in \mathbf{Perf}(X)$,
  we ask that the induced map $$\xymatrix{\Gamma(X, E) \otimes
    \mathbb{H}(X\textrm{rel}B, E^{\vee}) \ar[r]^-{\alpha} &
    \mathbb{H}(X\textrm{rel}B, \OO_X) \ar[r]^-{\eta_b} & k[-d] }$$
  yields, by adjunction, an equivalence $\Gamma(X, E) \simeq
  \mathbb{H}(X\textrm{rel}B, E^{\vee})[-d]$. Here $\alpha$ is the map
  induced on the vertical fibers of $$\xymatrix{\Gamma(X, E)\otimes
    \Gamma(X, E^{\vee}) \ar[rr]^-{\textrm{tr}} \ar[d]_-{\textrm{id}
      \otimes b_{E^{\vee}} } & & \Gamma(X, \OO_X) \ar[d]^-{b_{\OO_X}}
    \\ \Gamma(X, E)\otimes \Gamma(B, b^*E^{\vee})
    \ar[rr]^-{\textrm{tr}\, (b_E \otimes \textrm{id})} & &\Gamma(B,
    \OO_B).}$$ Given a $d$-orientation $\eta_b$ on the map $b: B \to
  X$, a $(n+1)$-shifted symplectic stack $Y'$, a map $\ell: Y\to Y'$,
  and a lagrangian structure on $f$, we may consider the derived
  stack $$\mathsf{MAP}_{\dSt}(b,f):= \mathsf{MAP}_{\dSt}(B, Y)
  \times_{\mathsf{MAP}_{\dSt}(B,Y')} \mathsf{MAP}_{\dSt}(X, Y')$$ of
  arrows from $b$ to $f$. The generalization of Theorem \ref{MAP} to
  this relative situation says that if $\mathsf{MAP}_{\dSt}(b, \ell)$
  and $\mathsf{MAP}_{\dSt}(X, Y)$ are derived Artin stacks lfp over
  $k$, then $\mathsf{MAP}_{\dSt}(b, \ell)$ has a canonical
  $(n-d+1)$-shifted symplectic form, and there is a canonical
  lagrangian structure on the natural map $\mathsf{MAP}_{\dSt}(X, Y)
  \to \mathsf{MAP}_{\dSt}(b, \ell)$. If we take $B$ to be empty (so
  that $\eta_b$ is just a $d$-orientation on $X$), and $Y=(\Spec \, k,
  \omega_{n+1}^{0})$, we have $\mathsf{MAP}_{\dSt}(b, \ell)\simeq
  (\Spec \, k, \omega_{n-d+1}^{0})$, and by Remark \ref{rmks} (3), we
  get back Theorem \ref{MAP}. But we may also take $Y=(\Spec \, k,
  \omega_{n+1}^{0})$, and $B$ arbitrary (non-empty), and we get a
  lagrangian structure on the restriction map $\mathsf{MAP}_{\dSt}(X,
  Y) \to \mathsf{MAP}_{\dSt}(B, Y)\simeq \mathsf{MAP}_{\dSt}(b,
  \ell)$, where $\mathsf{MAP}_{\dSt}(B, Y)$ is $(n-d+1)$-symplectic
  (consistently with Theorem \ref{MAP}, since $B$ is $(d-1)$-oriented
  by hypothesis, and $Y$ is $n$-shifted symplectic by Remark
  \ref{rmks} (3)). In other words, what we gain in this case, is that
  restriction to the boundary (for maps to a fixed shifted symplectic
  target) is endowed with a lagrangian structure. Examples of relative
  orientations includes topological examples (where $b$ is the
  inclusion of the boundary in a compact oriented topological
  $d$-manifold with boundary), and algebro-geometric ones where $b$ is
  the inclusion of the derived zero locus (Remark \ref{crit}) of a
  section of the anti-canonical bundle of a smooth projective variety
  of dimension $d$. When $B$ is a $K3$ surface inside a Fano $3$-fold,
  this might be compared with \cite[Proposition 2.2]{ty}. For more
  details, we address the reader to \cite[3.2.2]{cal}, and
  \cite[p. 227]{toenems}.

\end{rmk}

\

The last general existence theorem for shifted symplectic structures
is a generalization of the $2$-shifted symplectic structure on
$\mathsf{B}\mathrm{GL}_n$ described at the beginning of this
subsection. \\ Let $\mathsf{Perf}$ be the derived stack classifying
perfect complexes. It can be defined as the functor sending a cdga $A$
to the nerve of the category of cofibrant perfect $A$-dg-modules with
morphisms only the quasi-isomorphisms (as an $\infty$-functor it sends
$A$ to the coherent nerve of the Dwyer-Kan localization of the
previous category). The truncation of $\mathsf{Perf}$ is the (higher)
stack first introduced and studied in \cite{hs}. Though
$\mathsf{Perf}$ is not strictly speaking a derived Artin stack lfp
over $k$, it is quite close to it: it is \emph{locally geometric},
i.e. it is a union of open derived Artin substacks lfp over $k$,
e.g. $\mathsf{Perf}\simeq \cup_{n\geq 0} \mathsf{Perf}^{[-n,n]}$,
where $\mathsf{Perf}^{[-n,n]}$ is the derived stack classifying
perfect complexes of Tor-amplitude contained in $[-n,n]$
(\cite[Proposition 3.7]{tova}). In particular, it makes sense to ask
wether $\mathsf{Perf}$ carries a shifted symplectic structure, as we
know its substack $\mathsf{Perf}^{[-0,0]}=\mathsf{B}\mathrm{GL} \simeq
\coprod_{n\geq 0} \mathsf{B}\mathrm{GL}_n$ does. As shown in
\cite[Theorem 2.12]{ptvv}, the answer is affirmative.

 \begin{thm}\label{Perf}
 The derived stack $\mathsf{Perf}$ has a canonical $2$-shifted symplectic structure.
 \end{thm}

We will briefly give an idea of the proof, and address the readers to
\cite[Theorem 2.12]{ptvv} for all details. First of all, by
definition, $\mathsf{Perf}$ carries a universal perfect complex
$\mathcal{E}$, and we consider its perfect
$\mathcal{O}_{\mathsf{Perf}}$-Algebra of endomorphisms
$\mathcal{B}:=\mathbb{R}\underline{End}_{\mathcal{O}_{\mathsf{Perf}}}(\mathcal{E})\simeq
\mathcal{E}^{\vee} \otimes_{\mathcal{O}_{\mathsf{Perf}}}
\mathcal{E}$. One checks that $\mathbb{T}_{\mathsf{Perf}} \simeq
\mathcal{B}[1]$, and thus gets a well defined
map
$$
\xymatrix{\omega^{0}: \mathbb{T}_{\mathsf{Perf}}
  \wedge_{\mathcal{O}_{\mathsf{Perf}}} \mathbb{T}_{\mathsf{Perf}}
  \ar[r]^-{\sim} & Sym_{\mathcal{O}_{\mathsf{Perf}}} ( \mathcal{B})[2]
  \ar[r]^-{\textrm{mult}} & \mathcal{B}[2] \simeq \mathcal{E}^{\vee}
  \otimes_{\mathcal{O}_{\mathsf{Perf}}} \mathcal{E}[2]
  \ar[r]^-{\mathrm{tr[2]}} & \mathcal{O}_{\mathsf{Perf}}[2] }
$$
where
$tr$ denotes the trace (or evaluation) map for the perfect complex
$\mathcal{E}$. So, $\omega_{0}$ is a well defined, and non-degenerate
$2$-shifted 2-form on $\mathsf{Perf}$.  So, it only remains to  show
that there exists a $2$-shifted \emph{closed} 2-form $\omega$ on
$\mathsf{Perf}$, whose underlying 2-form is $\omega^0$. By
\cite{tove}, or \cite[Theorem 2.1]{hoy}, one can prove that the
weight $2$ component of the refined Chern character
$\mathsf{Ch}(\mathcal{E})$ (with values in negative cyclic homology)
as defined in \cite{chern} is indeed a $2$-shifted closed 2-form on
$\mathsf{Perf}$ whose underlying 2-form is $\frac{1}{2} \omega^0$,
thus $\omega := 2 \mathsf{Ch}_{2}(\mathcal{E})$ is indeed a
$2$-shifted symplectic form on $\mathsf{Perf}$ whose underlying 2-form
is $\omega^0$.

\

By combining Theorem \ref{MAP}, \ref{Lagr}, and \ref{Perf}, we get the
following (non-exhaustive) list of geometrically interesting,
classes of examples of $n$-shifted symplectic derived stacks:\\

\noindent
$\bullet$ the derived stack $\mathsf{Perf}(X):=\mathsf{MAP}_{\dSt}(X,
\mathsf{Perf}) $ of perfect complexes on a Calabi-Yau variety of
dimension $d$ ($n=2-d$).

\noindent
$\bullet$ the derived stack $\mathsf{Perf}(M) :=
\mathsf{MAP}_{\dSt}(\mathrm{const}(M), \mathsf{Perf})$ of perfect
complexes on a compact oriented topological manifold $M$ of dimension
$d$ ($n=2-d$).

\noindent
$\bullet$ the derived stack $\mathbb{R}\mathsf{Vect}_{n}(X)$ of rank
$\ n$ vector bundles on a Calabi-Yau variety of dimension $d$
($n=2-d$).

\noindent
$\bullet$ the derived stack $\mathbb{R}\mathsf{Vect}_{n}(M)$ of rank $
n$ vector bundles on a compact oriented topological manifold $M$ of
dimension $d$ ($n=2-d$).

\noindent
$\bullet$ the derived stack $\mathbb{R}\mathsf{Bun}_{G}(X)$ of
$G$-torsors on a Calabi-Yau variety of dimension $d$, for $G$
reductive ($n=2-d$).

\noindent
$\bullet$ the derived stack $\mathbb{R}\mathsf{Loc}_{G}(M)$ of
$G$-local systems on a compact oriented topological manifold $M$ of
dimension $d$ ($n=2-d$)

\noindent
$\bullet$ $\mathsf{MAP}_{\dSt}(X,Y)$, for $X$ a smooth and proper
$d$-dimensional Calabi-Yau scheme, and $Y$ a smooth symplectic scheme
($n=-d$).

\noindent
$\bullet$ $\mathsf{MAP}_{\dSt}(X,\mathrm{T}^*Y[m])$, for $X$ a smooth
and proper $d$-dimensional Calabi-Yau scheme, and $Y$ a derived Artin
stack lfp over $k$ ($n=m-d$).

\noindent
$\bullet$ iterated derived free loop spaces
$\mathsf{MAP}_{\dSt}(\mathrm{const}((S^1)^d),Y)$ of a smooth
symplectic scheme $Y$, and more generally,
$\mathsf{MAP}_{\dSt}(\mathrm{const}(M),Y)$, for $M$ a compact oriented
topological $d$-manifold, and $Y$ a smooth symplectic scheme ($n=-d$).

\noindent
$\bullet$ the derived stack $\mathbb{R}Loc^{DR}_{G}(X)$ of flat
$G$-bundles on $X$, a $d$-dimensional smooth and proper
Deligne-Mumford stack ($n= 2-2d$).

\noindent
$\bullet$ the derived moduli stack
$\mathsf{MAP}_{\dSt}(X_{\textrm{Dol}},\mathsf{Perf})$ of Higgs fields
on a proper and smooth Deligne-Mumford stack of dimension $d$
($n=2-2d$).

\noindent
$\bullet$ the derived moduli stack $\mathcal{M}_{T}$ of compact
objects in a $d$-Calabi-Yau dg-category $T$, e.g. the so-called
\emph{non-commutative $K3$} sub-dg-category of the derived dg-category
of a cubic $4$-fold ($n=2-d$, and $d=2$ for a non-commutative
$K3$).\footnote{See e.g. \cite[3.3]{ks} for a definition of a
  $d$-Calabi--Yau dg-category, and \cite[5.3]{toenems} for a sketch of
  a proof of this result.}
 
\begin{rmk} If $X$ is a proper and smooth scheme over $k$, and $D$ is a
  smooth Calabi-Yau divisor in the anti-canonical class, then Remark
  \ref{damienboundary} together with Thms. \ref{MAP} and \ref{Perf}
  give us a lagrangian structure on the restriction map
  $\mathsf{Perf}(X)= \mathsf{MAP}_{\dSt}(X, \mathsf{Perf}) \to
  \mathsf{MAP}_{\dSt}(D, \mathsf{Perf})=\mathsf{Perf}(D)$, The same is
  true for the restriction map $\mathbb{R}\mathsf{Bun}_{G}(X) \to
  \mathbb{R}\mathsf{Bun}_{G}(D)$ between the derived stack of torsors
  under a smooth reductive group scheme $G$ over $k$.
\end{rmk}

As a sample consequence of Theorem \ref{MAP}, \ref{Lagr}, and
\ref{Perf}, we recall the following important result (\cite[Corollary
  5.19]{bbj}) by Brav-Bussi-Joyce, establishing the existence of a
local \emph{algebraic} potential for Donaldson-Thomas moduli spaces
attached to Calabi-Yau $3$-folds.

\begin{thm}\label{joyce}\emph{(Brav-Bussi-Joyce)}
The moduli space $\mathcal{M}_{DT}$ of simple coherent sheaves or of
complexes of coherent sheaves on a Calabi-Yau $3$-fold is
Zariski-locally isomorphic, as a $(-1)$-shifted symplectic derived
scheme, to the derived critical locus (as in Remark \ref{crit}) of a
regular function $f: U \to \mathbb{A}_{k}^{1}$ over a smooth
$k$-scheme $U$.
\end{thm}

The function in the statement is called the \emph{Donaldson-Thomas
  potential}. Our existence theorems combine to give
$\mathcal{M}_{DT}$ a $(-1)$-shifted symplectic structure, and the
authors achieve Theorem \ref{joyce} by proving a local structure
theorem (Darboux style) for derived schemes $X$ endowed with a
$(-1)$-shifted symplectic structure: any such $X$ is Zariski-locally a
derived critical locus of a regular function on a smooth scheme. A
similar statement in the $\mathbb{C}$-analytic category was proved
before by Joyce-Song \cite[Theorem 5.4]{js} using Gauge Theory, and a
version valid formally locally at any point of $\mathcal{M}_{DT}$, by
Kontsevich-Soibelman \cite[Section 3.3]{ks}. Obviously, Theorem
\ref{joyce} is a considerably stronger and more precise result.


\section{Shifted Poisson structures}\label{poiss}

Setting up a notion of shifted Poisson structure for sufficiently
general derived Artin stacks turns out to be much more complicated
than the case of shifted symplectic structures, described in the
previous section. On the other hand, a shifted Poisson structure on a
derived stack $F$ is exactly the right structure that controls the
deformation quantization of the $\infty$-category of perfect complexes
on $F$.  Therefore, in order to establish deformation quantizations for
all the shifted symplectic derived moduli spaces listed at the end of
Section \ref{symplectic}, one also needs a comparison between shifted
symplectic structures and (non-degenerate) shifted Poisson
structures. Unfortunately, this comparison, which classically takes no
more than two lines, is rather tricky in the derived setting, due to
the fact that all the structures involved in the comparison are weak
ones i.e. defined up-to-homotopy. The general theory of shifted
Poisson structures, a comparison with shifted symplectic structures,
and applications to deformation quantization of derived moduli spaces,
have all been developed recently in \cite{cptvv}, which is more than
100 pages long. In this Section we give a summary of the main
constructions and results from \cite{cptvv}, and a guide to its
reading.

\subsection{Differential calculus}\label{diff}
In order to be able to define and study shifted Poisson structures on
derived Artin stacks, we will need to have at our disposal a machinery
of derived differential calculus (de Rham complex, shifted polyvectors
etc.) working in sufficiently general setting. Thus, let $\C$ be a
 stable presentable symmetric monoidal $\infty$-category which is
$k$-linear i.e. enriched over the $\infty$-category $\dg_k$
(\cite[Definition 2.0.0.7]{lualg}). In this paper, $\C$ will always be
obtained as the coherent nerve of the Dwyer-Kan localization of a
$\mathsf{dg}_k$-enriched symmetric monoidal model category
$\mathsf{M}$ satisfying some additional technical properties for
which we address the reader to \cite[1.1 and 1.2]{cptvv} \footnote{By \cite{ns}, it is likely that \emph{any} $k$-linear stable presentable symmetric monoidal $\infty$-category can be obtained this way.}. And 
we will suppose that the enrichment is
induced by a symmetric monoidal functor $\mathsf{dg}_k \to \m$.

For our present purposes, it will be enough to keep in mind some of
the examples we will be most frequently interested in: $\C$ could be
$\dg_k$ itself, the $\infty$-category $\dg^{gr}_k$ of
$\mathbb{Z}$-graded dg-modules over $k$, the $\infty$-category
$\edg_k$ of mixed dg-modules over $k$, the $\infty$-category $\medg_k$
of $\mathbb{Z}$-graded mixed dg-modules over $k$, or more generally,
any category of diagrams in the previous examples.

We denote by $\mathsf{CAlg}_{\mathsf{M}}$ the model category of
commutative algebras in $\mathsf{M}$, and by $\mathbf{CAlg}_{\C}=
L(\mathsf{CAlg}_{\mathsf{M}})$, the corresponding
$\infty$-category.

The symmetric monoidal model category $\epsilon-\m^{gr}$ of mixed
graded objects in $\m$ is defined by replacing $\mathsf{dg}_k$ with
$\m$ in the definition at the beginning of Section \ref{DEF} (and with
the cohomological shift defined here as $P[1]:=P^{k[-1]}$, for $P\in
\m^{gr}$, using the co-tensor enrichment of $\m$ over
$\mathsf{dg}_k$). The model category of commutative monoids in
$\epsilon-\m^{gr}$ is denoted by
$\epsilon-\mathsf{CAlg}^{gr}_{\mathsf{M}}$, and called the model
category of graded mixed algebras in $\m$.  The corresponding
$\s$-categories will be denoted by $\epsilon-\C^{gr}$, and
$\epsilon-\mathbf{CAlg}^{gr}_{\C}$. Note that we have a canonical
equivalence of $\s$-categories $\mathbf{CAlg}(\epsilon-\C^{gr}) \simeq
\epsilon-\mathbf{CAlg}^{gr}_{\C}$ (\cite[Theorem 4.5.4.7]{lualg}).

\

\subsubsection{De Rham theory}
As explained in Remark \ref{generalM},
for any $A \in \calg_{\C}$, we have the associated de Rham algebra
$\DR^{\C}(A)$,
a mixed graded algebra in $\C$, where $\DR^{\C} : \calg_{\C} \to
\ealgrM$ is the left adjoint to the $\s$-functor sending $B \in
\ealgrM$ to its weight $0$ part $B(0)$. Note that if $A \in
\cdga_{k}$, then $\DR(A)$ of Definition \ref{d1.1} is exactly
$\DR^{\C}(A)$, with $\C= \dg_k$. There is an
$\s$-functor
$$
A-\mathbf{Mod}_{\C} \longrightarrow \T \, ,  \qquad P
\longmapsto \mathsf{Der}(A,P) := \mathsf{Map}_{\calg_{\C}/A}(A,
A\oplus P)
$$
which is co-represented by an $A$-module
$\mathbb{L}^{\C}_{A}$ in $\C$, called the \emph{cotangent complex} or
\emph{cotangent object} of $A$. As in Section \ref{DEF}, one can prove
(\cite[Proposition 1.3.12]{cptvv}) that the canonical
map $$Sym_{A}(\mathbb{L}^{\C}_{A}[-1]) \longrightarrow \DR^{\C}(A)$$
is functorial in $A$, and an equivalence in $\calg_{\C}^{gr}$
(i.e. forgetting the mixed structure in the target). In other words,
the construction $\DR^{\C}$ endows the graded algebra
$Sym_{A}(\mathbb{L}^{\C}_{A}[-1])$ with a canonical (weak) mixed
differential, the (weak) de Rham differential.

For the sake of
brevity, we will omit  the completely
analogous definitions and results in the \emph{relative setting}, i.e. for
morphisms in $\calg_{\C}$. We will instead say  a few words about
\emph{strict models} of the above constructions, i.e. inside the
\emph{model} category $\mathsf{M}$. First of all, if $A' \in
\calg_{\mathsf{M}} $, then the functor
$$
A'-\mathbf{Mod}_{\mathsf{M}}
\longrightarrow \T \, , \qquad
P \longmapsto Der(A',P) :=
Hom_{\calg_{\mathsf{M}}/A'}(A', A'\oplus P)
$$
is co-represented by an
$A'$-module $\Omega_{\mathsf{M}, A'}^{1}$ in $\mathsf{M}$ called the
module of \emph{K\"ahler differentials} of $A'$. If $A \in \calg_{\C}$
and $QA$ is a cofibrant model for $A$ in $\calg_{\mathsf{M}}$, then we
have a canonical equivalence $\Omega_{\mathsf{M}, QA}^{1} \simeq
\mathbb{L}^{\C}_{A}$ in $A-\mathbf{Mod}_{\C} \simeq
QA-\mathbf{Mod}_{\C}$. Furthermore, the functor $\ealgrm \to
\calg_{\mathsf{M}}$ selecting the weight $0$ component, has a left
adjoint $\DR^{\mathsf{M}}$, thus for any $A' \in \calg_{\mathsf{M}}$,
$\DR^{\mathsf{M}}(A')$ is a strict mixed graded algebra in
$\mathsf{M}$. Still for an arbitrary $A' \in \calg_{\mathsf{M}}$, we
also have a functorial
\emph{isomorphism} $$Sym_{A'}(\Omega_{\mathsf{M}, A'}^{1}[-1])
\longrightarrow \DR^{\mathsf{M}}(A')$$ in the category
$\calg^{gr}_{\mathsf{M}}$ of strict graded algebras in
$\mathsf{M}$. This precisely says that $Sym_{A'}(\Omega_{\mathsf{M},
  A'}^{1}[-1])$ has a strict mixed structure, the de Rham
differential. If $A \in \calg_{\C}$ and $QA$ is a cofibrant model for
$A$ in $\calg_{\mathsf{M}}$, then there is a an equivalence $
\DR^{\mathsf{M}}(QA) \simeq \DR^{\C}(A)$ in $\ealgrM$,
i.e. $\DR^{\mathsf{M}}(QA)$, or equivalently
$Sym_{A'}(\Omega_{\mathsf{M}, QA}^{1}[-1])$ with the induced mixed
structure, is a strict model for $\DR^{\C}(A)$ in $\ealgrM$.

\

\noindent \textbf{Differential forms.}  Again as in Section \ref{DEF},
as suggested in Remark \ref{generalM}, we may define (closed) shifted
differential forms for commutative algebras in $\C$.
\begin{df} Let $A \in \calg_{\C}$. \begin{itemize}
  \item The space of \emph{closed $n$-shifted $p$-forms} on $A$
    is
    $$
    \mathcal{A}_{\C}^{p,cl}(A,n):=
    \mathsf{Map}_{\epsilon-\C^{gr}}(\mathbf{1}_{\C}(p)[-p-n],\DR^{\C}(A))  \in \mathbb{T}.
    $$
An element in $\pi_0
    (\mathcal{A}_{\C}^{p,cl}(A,n))$ is called a \emph{closed
      $n$-shifted $p$-form} on $A$.
\item The space of \emph{$n$-shifted $p$-forms} on $A$ is
  $\mathcal{A}_{\C}^{p}(A,n):=\mathsf{Map}_{\C}(\mathbf{1}_{\C}[-n],
  \wedge_{A}^{p}\mathbb{L}^{\C}_{A}) \in \mathbb{T}$.
  An element in $\pi_0 (\mathcal{A}_{\C}^{p}(A,n))$
  is called a \emph{$n$-shifted $p$-form} on $A$.
\item Since $Sym_{A}(\mathbb{L}^{\C}_{A}[-1]) \simeq \DR^{\C}(A)$,
  there is an induced map $u: \mathcal{A}^{p,cl}(A,n) \to \mathcal{A}^{p}(A,n)$,
  called the \emph{underlying $p$-form map}.
\end{itemize}
\end{df}

If $\omega$ is a $n$-shifted $2$-form on $A$, and we assume that
$\mathbb{L}^{\C}_{A}$ is dualizable in $A-\mathbf{Mod}_{\C}$, then
$u(\omega)$ induces a map $u(\omega)^{\flat}: \mathbb{T}^{\C}_{A} \to
\mathbb{L}^{\C}_{A}[n]$ in $A-\mathbf{Mod}_{\C}$, where
$\mathbb{T}^{\C}_{A}$ denotes the dual of $\mathbb{L}^{\C}_{A}$. We
say that $\omega$ is \emph{non-degenerate} if $u(\omega)^{\flat}$ is
an equivalence.

\begin{df}\label{symplabstrdef} Let $A \in \calg_{\C}$.
  The \emph{space of $n$-shifted symplectic structures on} $A$ is the
  subspace $\mathsf{Sympl}(A,n)$ of $\mathcal{A}_{\C}^{p,cl}(A,n)$
  whose connected components consist of non degenerate forms.
\end{df}

\begin{rmk} Note that, even if the notation $\mathsf{Sympl}(A,n)$
  does not record $\C$, this space obviously depends on the
  category $\C$ inside which we are working. 
In Section \ref{FL}, we will explain how the abstract Definition
\ref{symplabstrdef} gives back the definition of a shifted symplectic
structure on a derived Artin stack (Definition \ref{sympl}).
\end{rmk}

\

\subsubsection{Polyvectors}
Let $\mathbb{L}ie^{gr}_{k}$ be the graded Lie operad in
$\mathsf{dg}_k$, where the bracket operation has degree $-1$. We let
$\mathsf{dgLie}_{k}^{gr}$ be the model category of
$\mathbb{L}ie^{gr}_{k}$-algebras, i.e. graded Lie-algebras in
$\mathsf{dg}_k$, where the (strict) Lie bracket has weight $-1$; we
denote by $\mathbf{dgLie}_{k}^{gr}$, the corresponding $\s$-category
$L(\mathsf{dgLie}_{k}^{gr})$. By using our enriching functor
$\mathsf{dg}_k \to \mathsf{M}$, we get an operad
$\mathbb{L}ie^{gr}_{\mathsf{M}}$ in $\mathsf{M}$. Taking algebras with
respect to this operad, we have $\mathsf{Lie}_{\m}^{gr}$, the model
category of graded Lie-algebras in $\m$, and $\mathbf{Lie}_{\C}^{gr}$,
the corresponding $\s$-category
$L(\mathsf{Lie}_{\m}^{gr})$.

Similarly, for $n \in \mathbb{Z}$, we
will denote by $\mathbb{P}_n$ the dg-operad (i.e. the operad in
$\mathsf{dg}_k$) whose algebras are Poisson cdga's with a bracket of
degree $(1-n)$. Recall that, for $n>1$, $\mathbb{P}_n$ can be
identified with the operad of chains of the topological $n$-little
disks operad $\mathbb{E}_n$ (see \cite{cohen}). We will also be
interested in a graded version of $\mathbb{P}_n$, denoted by
$\mathbb{P}^{gr}_n$: this is an operad in $\mathsf{dg}^{gr}_k$, whose
algebras have multiplication of weight (= external grading) $0$, and
bracket of weight $-1$. The corresponding model categories of algebras
will be denoted by $\mathbb{P}_n -\mathsf{Alg}_{\mathsf{dg}_k} =
\mathbb{P}_n -\mathsf{cdga}_k$, and $\mathbb{P}^{gr}_n
-\mathsf{Alg}_{\mathsf{dg}_k} = \mathbb{P}_n -\mathsf{cdga}^{gr}_k$;
the associated $\s$-categories by $\mathbb{P}_n
-\mathbf{Alg}_{\mathsf{dg}_k} = \mathbb{P}_n -\mathbf{cdga}_k$, and
$\mathbb{P}^{gr}_n -\mathbf{Alg}_{\mathsf{dg}_k} = \mathbb{P}_n
-\mathbf{cdga}^{gr}_k$.

By using our enriching functor $\mathsf{dg}_k \to \mathsf{M}$, we thus
get operads $\mathbb{P}_{\mathsf{M}, n}$ in $\mathsf{M}$, and
$\mathbb{P}^{gr}_{\mathsf{M}, n}$ in $\mathsf{M}^{gr}$. Note that
$\mathbb{P}^{gr}_{\mathsf{M}, n}$-algebras are commutative algebras
$B$ in $\m^{gr}$ endowed with a bracket operation $[-,-]_{p,q}:
B(p)\otimes_{\m}B(q) \to B(b+q-1)[1-n]$, which is a graded
bi-derivation, and endows $B[n-1]$ with a Lie algebra structure inside
$\m^{gr}$. The corresponding $\s$-categories
$L(\mathbb{P}_{\mathsf{M}, n}-\mathsf{Alg})$ and
$L(\mathbb{P}^{gr}_{\mathsf{M}, n}-\mathsf{Alg})$ of algebras over
these operads, will be denoted by $\mathbb{P}_{\C, n}-\mathbf{Alg}$
and $\mathbb{P}^{gr}_{\C, n}-\mathbf{Alg}$, respectively.

For $A' \in \mathsf{CAlg}_{\m}$, $p \geq 0$, and $m \in \mathbb{Z}$,
we define the object (in $\m$) of $m$-shifted, degree $p$
\emph{polyvectors} $\mathsf{Pol}^{\m}(A',n)(p)$ on $A'$, as
follows. We consider $\ti^{(p)}(A', m) \in \m$ the sub-object of
$\underline{Hom}_{\m}(A'^{\otimes^{p} }, A'[mp])$ (where
$\underline{Hom}_{\m}$ is the internal Hom-object in $\m$) consisting
of (shifted) multiderivations (\cite[Section 1.4.2.]{cptvv}). There is
a natural action of the symmetric group $\Sigma_p$ on
$\underline{Hom}_{\m}(A'^{\otimes^{p} }, A'[mp])$ induced by its
standard action on $A'^{\otimes^{p}}$ and by the $(-1)^{m} \times
\mathrm{Sign}$ action on $A'[mp]$. this action restricts to an action
on $\ti^{(p)}(A', m)$, and we denote by $\ti^{(p)}(A',
m)^{\Sigma_p}$ the $\m$-object of $\Sigma_p$-invariants. By
standard conventions, we put $\ti^{(0)}(A', m)= A'$, and notice that
$\ti^{(1)}(A', m)= \underline{Hom}_{\m}(\Omega_{\m, A'}^{1}, A'[m])$.

\begin{df}\label{polydef}
  For $A' \in \mathsf{CAlg}_{\m}$, and $n \in \mathbb{Z}$, we define
  the $\m^{gr}$\emph{-object of} $n$\emph{-shifted polyvectors on}
  $A'$ to be
  $$
  \mathsf{Pol}^{\m}(A',n):= \bigoplus_{p}\ti^{(p)}(A',
  -n)^{\Sigma_p}.
  $$
\end{df}

As in the classical case, there is a ``composition'' (by insertion) of
shifted multiderivations and there is a Schouten-Nijenhuis-like bracket of
shifted multiderivations, so that $\mathsf{Pol}^{\m}(A',n)$
\emph{becomes a graded} $\mathbb{P}_{ n+1}$\emph{-algebra} in $\m$,
i.e. an object in $\mathbb{P}_{\mathsf{M}, n+1}-\mathsf{Alg}$.

\begin{rmk}\label{yok0}
  If $A'$ is such that $\Omega_{\m, A'}^{1}$ is (strictly) dualizable
  in $A'-\mathsf{Mod}_{\m}$, then it is easy to prove that there is an
  isomorphism $\ti^{(p)}(A', m)^{\Sigma_p} \simeq
  Sym^{p}(\mathrm{T}_{\m, A'}[-m])$ in $\m$, and more generally
  $\mathsf{Pol}^{\m}(A',n) \simeq \oplus_{p\geq 0}
  Sym^{p}(\mathrm{T}_{\m, A'}[-m])$ in $\mathbb{P}_{\mathsf{M},
    n+1}-\mathsf{Alg},$ where $\mathrm{T}_{\m, A'}$ is the $A'$-dual
  of $\Omega_{\m, A'}^{1}$ (and, on the rhs the multiplication is the
  one induced by $Sym$, while the bracket is the one canonically
  induced by the Lie bracket in $\mathrm{T}_{\m, A'}$). Note that,
  however, we have \emph{not} used $\mathrm{T}_{\m, A'}$ directly, in
  order to \emph{define} $\mathsf{Pol}^{\m}(A',n)$. We used
  multiderivations instead, and the the possible definitions agrees
  only under the hypotheses that $\Omega_{\m, A'}^{1}$ is (strictly)
  dualizable.
\end{rmk}

The construction $A' \mapsto \mathsf{Pol}^{\m}(A',n)$ is not fully
functorial, since we can neither pullback nor pushforward multiderivations in
general. However it is possible to define a restricted functoriality
(\cite[Lemma 1.4.13]{cptvv}) at the level of $\s$-categories, enabling
us to give the following

\begin{df}\label{Polydef}
  Let $n \in \mathbb{Z}$, and $\mathbf{CAlg}^{\textrm{f\'et}}_{\C}$
  the sub $\s$-category of $\mathbf{CAlg}_{\C}$ whose maps are only
  the formally \'etale ones (i.e.  maps $A\to B$ such that
  $\mathbb{L}_{B/A} \simeq 0$): Then there is a well-defined
  $\s$-functor $$\mathbf{Pol}^{\C}(-,n):
  \mathbf{CAlg}^{\textrm{f\'et}}_{\C} \longrightarrow \mathbb{P}_{\C,
    n+1}-\mathbf{Alg},$$ such that if $A \in \mathbf{CAlg}_{\C}$ and
  $A'$ is a fibrant-cofibrant replacement of $A$ inside
  $\mathsf{CAlg}_{\m}$, then we have an equivalence
  $\mathbf{Pol}^{\C}(A,n) \simeq \mathsf{Pol}^{\m}(A',n)$ in
  $\mathbb{P}_{\C, n+1}-\mathbf{Alg}$.
\end{df}

\begin{rmk}\label{yok} If $A \in \mathbf{CAlg}_{\C} $ is such that
  $\mathbb{L}^{\C}_{A}$
  is dualizable in $A-\mathsf{Mod}_{\C}$, then one deduces from Remark
  \ref{yok0} an equivalence $\mathbf{Pol}^{\C}(A,m) \simeq
  \oplus_{p\geq 0} Sym^{p}(\mathbb{T}^{\C}_{A}[-m])$.
\end{rmk}

We are now in a position to give the definition of a shifted Poisson
structure. Recall that if $A \in \mathbf{CAlg}_{\C}$, then
$\mathbf{Pol}^{\C}(A,n) \in \mathbb{P}_{\C, n+1}-\mathbf{Alg}$, so
that, in particular, $\mathbf{Pol}^{\C}(A,n)[n] \in
\mathbf{Lie}_{\C}^{gr}$.

\begin{df}\label{spoissdef}
  Let $n \in \mathbb{Z}$, and $A \in \mathbf{CAlg}_{\C}$. The space of
  $n$-shifted Poisson structures on $A$ is $$\mathbf{Poiss}(A,n):=
  \mathsf{Map}_{\mathbf{Lie}_{\C}^{gr}}(\mathbf{1}_{\C}[-1](2),
  \mathbf{Pol}^{\C}(A,n+1)[n+1]),$$ where $\mathbf{1}_{\C}[-1](2)$ is
  the Lie algebra object in $\C$ given by $\mathbf{1}_{\C}$ sitting in
  pure cohomological degree $1$, and pure weight $2$, with (strictly)
  trivial bracket.
\end{df}

This definition mimicks the classical one: if $\m=
\mathsf{dg}_k$, and we replace $ \mathsf{Map}$ by strict $Hom$-set in
$\mathsf{dgLie}_{k}^{gr}$, in the above definition, we get exactly the
\emph{set} of Poisson bivectors on $A$ (the bracket being trivial on
$\mathbf{1}[-1](2) = k[-1](2)$, the image of $1 \in k$ gives a
biderivation on $\pi$ on $A$, such that $[\pi,\pi]=0$).

\

\subsubsection{Standard realizations over $k$}

For the theory of shifted Poisson structures on derived stacks that we
develop in Section \ref{sps} we will systematically need to pass to
global sections. This procedure can be implemented already at the
level of generality of the present section, by considering what we
call \emph{standard realizations} over $k$. A more technical notion of
realization over $k$ (the \emph{Tate realization}) will be discussed later on
in this Section.

One of our standing hypotheses on the base symmetric monoidal
$\mathsf{dg}_{k}$ model category $\m$ is  that the unit
$\mathbf{1}_{\m}$ is cofibrant. Hence we get a Quillen adjunction,
with $ \mathbf{1}_{\m}\otimes -$ left adjoint,
$$
\xymatrix{
  \mathbf{1}_{\m}\otimes - : \mathsf{dg}_k \ar@<-.5ex>[r] &
  \ar@<-.5ex>[l] \m : \underline{Hom}_{k}(\mathbf{1}_{\m}, -) }
$$
where $\otimes$ denotes the tensor $\mathsf{dg}_k$-enrichment of $\m$,
and $\underline{Hom}_{k}$ is the Hom-enrichment of $\m$ in
$\mathsf{dg}_k$. This induces a derived $\s$-adjunction on the
associated $\s$-categories
$$
\xymatrix{ \mathbf{1}_{\m}\otimes - :
  \mathbf{dg}_k \ar@<-.5ex>[r] & \ar@<-.5ex>[l] \C :
  \mathbb{R}\underline{Hom}_{k}(\mathbf{1}_{\m}, -)\, , }
$$
and we
define the \emph{standard realization functor} as the right adjoint
$|-|:=\mathbb{R}\underline{Hom}_{k}(\mathbf{1}_{\m}, -)$. 
Actually
this is the first of a series of realization functors induced on
various categories of algebras in $\m$. In fact, since
$\mathbf{1}_{\m}$ is a comonoid object in $\m$ (hence in $\C$), $|-|$
is actually a lax symmetric monoidal $\s$-functor, hence it is a
\emph{right adjoint} on various functor on categories of ``algebras''
in $\C$. Our notation will always be $|-|$ for each of these induced
realization functors. For example, we have
$$
|-|: \mathbf{CAlg}_{\C}
\to \mathbf{cdga}_k \, , \,\,\,\,\, |-|: \mathbf{CAlg}^{gr}_{\C} \to
\mathbf{cdga}^{gr}_k \, , \,\,\,\,\, |-|: \epsilon-\mathbf{CAlg}_{\C}
\to \epsilon-\mathbf{cdga}_k\, , \,\,\,\,\, |-|:
\epsilon-\mathbf{CAlg}^{gr}_{\C} \to \epsilon-\mathbf{cdga}^{gr}_k\,
, $$ $$|-|: \mathbf{Lie}_{\C} \to \mathbf{dgLie}_k \, , \,\,\,\, |-|:
\mathbf{Lie}^{gr}_{\C} \to \mathbf{dgLie}^{gr}_k \, , \,\,\,\, |-|:
\mathbb{P}_{\C, n}-\mathbf{Alg} \to \mathbb{P}_{n}-\mathbf{cdga}_k\, ,
\,\,\,\, |-|: \mathbb{P}^{gr}_{\C, n}-\mathbf{Alg} \to
\mathbb{P}_{n}-\mathbf{cdga}^{gr}_k .
$$
In particular we get realizations
$$ | \mathbb{L}^{\C}_{B/A} | =: \mathbb{L}_{B/A} \in \mathbf{dg}_k \,
, \,\,\,\, |\DR^{\C}(B/A) | =: \DR(B/A) \in
\epsilon-\mathbf{cdga}^{gr}_k\, , \,\,\,\,  |\mathbf{Pol}^{\C}(A,n)| =:\mathbf{Pol}(A,n)
 \in \mathbb{P}_{n+1}-\mathbf{cdga}_k^{gr}.
$$
Note that when $\C= \dg_k$, all this realization functor are
(equivalent to) the identity functors.\\ 

Later on, we will consider
realization functors on category of diagrams in $\C$, and it will be
useful to recall the following fact. The $\dg_k$-enriched Hom in a
category of diagrams $I^{op} \to \C$,
satisfies
$$
\underline{Hom}_{\C^{I^{op}}}(\mathbf{1}_{\C^{I^{op}}}, F)
\simeq \lim_{x \in I^{op}} \underline{Hom}_{\C}(\mathbf{1}_{\C},
F(x)),
$$
since the monoidal unit $\mathbf{1}_{\C^{I^{op}}}$ is given
by the constant $I^{op}$-diagram at $\mathbf{1}_{\C}$. This
observation is the reason for the appearance, in Thm \ref{yeah}, of
(derived) global sections of various algebras over the derived affine
site of a derived stack.
\begin{ex}
  \emph{As examples of the use of standard realization, by definition
    of $|-|$ as a right adjoint, we have, for any $A \in
    \mathbf{CAlg}_{\C}$, the following equivalences in
    $\T$ $$\mathbf{Poiss}(A,n)\simeq
    \mathsf{Map}_{\mathbf{dgLie}_{k}^{gr}}(k[-1](2),
    \mathbf{Pol}(A,n+1)[n+1]),$$
    $$
    \mathcal{A}_{\C}^{p,cl}(A,n) \simeq
    \mathsf{Map}_{\epsilon-\dg_{k}^{gr}}(k(p)[-p-n],\DR(A)) \, ,
    \,\,\,\,\, \mathcal{A}_{\C}^{p}(A,n) \simeq \mathsf{Map}_{\dg_k}(k[-n],
    \wedge_{|A|}^{p}\mathbb{L}_{A}).
    $$
 }
 \end{ex}
 
\

\subsubsection{Comparison between Poisson algebras
  and shifted Poisson pairs} We need a few notations.  For an
arbitrary $\s$-category $\mathbf{C}$, $\mathcal{I}(\mathbf{C})$ will
denote its \emph{moduli space} or \emph{maximal}
$\infty$\emph{-subgroupoid} of equivalences (so that
$\mathcal{I}(\mathbf{C}) \in \T$), $\mathbf{C}^{\Delta[1]}$ the
$\s$-category of morphisms in $\mathbf{C}$, and $\textrm{ev}_0,
\textrm{ev}_1 : \mathbf{C}^{\Delta[1]} \to \mathbf{C}$ the source and
target $\s$-functors. For a model category $\mathsf{C}$, we denote by
$\mathsf{C}^{W}_{cf}$ the category of cofibrant-fibrant objects in
$\mathsf{C}$ with morphisms given by weak equivalences.

The nerve of $\mathsf{C}^{W}_{cf}$ is called the \emph{moduli space}
or \emph{classifying space} (of objects and equivalences) of
$\mathsf{C}$.  It is a general fact (going back essentially to Dwyer
and Kan) that we have an equivalence of spaces between the moduli
space of $\mathsf{C}$ and the space $\mathcal{I}(\mathbf{C})$. In
other words, we may (and will) identify the moduli space of
$\mathsf{C}$ with the moduli space of the associated $\s$-category
$\mathbf{C}$. We will be interested, below, in the moduli spaces of the 
$\s$-categories $\mathbb{P}_{\C, n}-\mathbf{Alg}$,
$\mathbf{Lie}_{\C}^{gr}$, and $\mathbf{CAlg}_{\C}$, associated to the
model categories $\mathbb{P}_{\m, n+1}-\mathsf{Alg}$,
$\mathsf{Lie}_{\m}^{gr}$, and $\mathsf{CAlg}_{\m}$, respectively.

For $A \in \mathbb{P}_{\C, n}-\mathbf{Alg}$, we will write
$\mathbb{P}_{\C, n}(A)$ for the space of $\mathbb{P}_{n}$ algebra
structures on $A$ having the given underlying commutative algebra
structure. More precisely, there is a forgetful $\s$-functor $U:
\mathbb{P}_{\C, n}-\mathbf{Alg} \to \mathbf{CAlg}_{\C}$ (forgetting
the bracket structure), and the fiber of $U$ at a given $A \in
\mathbf{CAlg}_{\C}$ is an $\s$-groupoid, i.e. a space, that we denote
by $\mathbb{P}_{\C, n}(A)$.

We have a canonical $\s$-functor
$$ v: \mathcal{I}(\mathbf{CAlg}_{\C}) \ni A \longmapsto
(\mathbf{1}[-1](2), \mathbf{Pol}^{\C}(A, n)[n]) \in
\mathbf{LieAlg}_{\C}^{gr} \times \mathbf{LieAlg}_{\C}^{gr}
$$
(note that $A \mapsto \mathbf{Pol}^{\C}(A, n)$ is functorial with
respect to equivalences in $\mathbf{CAlg}_{\C}$, since obviously
equivalences are formally \'etale). Borrowing ideas from
\cite[3.1]{cliff}, we define the space
$\mathbf{Poiss}^{\mathrm{eq}}_{\C,n}$ as the following pullback of in
$\T$
$$
\xymatrix{\mathbf{Poiss}^{\mathrm{eq}}_{\C, n} \ar[r] \ar[d]_-{q} &
  \mathcal{I}((\mathbf{LieAlg}_{\C}^{gr})^{\Delta[1]})
  \ar[d]^-{(\mathcal{I}(\mathrm{ev}_0) ,\mathcal{I}(\mathrm{ev}_1))} \\
  \mathcal{I}(\mathbf{CAlg}_{\C})
  \ar[r]_-{\mathcal{I}(v)} & \mathcal{I}(\mathbf{LieAlg}_{\C}^{gr}
  \times \mathbf{LieAlg}_{\C}^{gr}  )\, .}
$$

Thus, informally speaking, $\mathbf{Poiss}^{\mathrm{eq}}_{\C, n}$ is
the moduli space of $n$\emph{-shifted Poisson pairs} $(A, \pi)$,
consisting of a commutative algebra $A$ in $\C$ together with a
$n$-shifted Poisson structure on $A$ (see Definition \ref{spoissdef}).
More precisely, let $\mathsf{Poiss}^{\textrm{eq}}_{\mathsf{M}, n}$ be
the category whose objects are pairs $(A, \pi)$ where $A$ is a
fibrant-cofibrant object in $\mathsf{CAlg}_{\mathsf{M}}$, and $\pi$ is
a $n$-shifted Poisson structure on $A$, i.e a map
$\mathbf{1}_\mathsf{M}[-1](2) \to
\mathsf{Pol}^{\mathsf{M}}(A,n+1)[n+1] $ in the homotopy category of
$\mathbf{LieAlg}_{\mathsf{M}}^{gr}$, and whose morphisms $(A,\pi) \to
(A', \pi')$ are weak equivalences $u: A \to A'$ in
$\mathsf{CAlg}_{\mathsf{M}}$ such that the diagram
$$
\xymatrix{ &
  \mathsf{Pol}^{\mathsf{M}}(A,n+1)[n+1]
  \ar[dd]^-{\mathsf{Pol}^{\mathsf{M}}(u, n+1)[n+1]}
  \\ \mathbf{1}_\mathsf{M}[-1](2) \ar[ru]^-{\pi} \ar[rd]_-{\pi'} &
  \\ & \mathsf{Pol}^{\mathsf{M}}(A',n+1)[n+1] }
$$
is commutative in
the homotopy category of $\mathsf{LieAlg}_{\mathsf{M}}^{gr}$. Then,
the nerve of $\mathsf{Poiss}^{\mathrm{eq}}_{\mathsf{M}, n}$ is
equivalent to $\mathbf{Poiss}^{\mathrm{eq}}_{\C, n}$.  


There is a (strict) functor $w: \mathbb{P}_{\mathsf{M},
  n+1}-\mathbf{Alg}_{cf}^{W} \to
\mathsf{Poiss}^{\mathrm{eq}}_{\mathsf{M}, n}$, sending a strict
$\mathbb{P}_{n+1}$-algebra $B$ in $\mathsf{M}$ to the pair $(B, \pi)$,
where $\pi$ is induced, in the standard way, by the (strict) Lie
bracket on $B$ (since the bracket is strict, it is a strict
biderivation on $B$). Restriction to weak equivalences (between
cofibrant-fibrant objects) in $\mathbb{P}_{\mathsf{M},
  n+1}-\mathbf{Alg}$ ensures this is a functor, and note that objects
in the image of $w$ are, by definition, \emph{strict} pairs $(B,
\pi)$, i.e. the shifted Poisson structure $\pi:
\mathbf{1}_\mathsf{M}[-1](2) \to \mathsf{Pol}^{\mathsf{M}}(A,n+1)[n+1]
$ is an actual morphism in $\mathbf{LieAlg}_{\mathsf{M}}^{gr}$ (rather
than a map in its homotopy category). We have a commutative diagram of
functors
$$
\xymatrix{ \mathbb{P}_{\mathsf{M},
    n+1}-\mathbf{Alg}_{cf}^{W} \ar[dr]_-{u} \ar[rr]^-{w} &
  &\mathsf{Poiss}^{\mathrm{eq}}_{\mathsf{M}, n} \ar[dl]^-{q} \\ &
  (\mathsf{CAlg}_{\mathsf{M}})_{cf}^{W} & }
$$
where $u$ forgets the
bracket structure, and $q$ is the functor $(A,\pi) \mapsto A$. Taking
the nerves of the previous diagram, we get a commutative diagram in
$\T$ (where we have kept the same name for the maps)
$$
\xymatrix{
  \mathcal{I}(\mathbb{P}_{\C, n+1}-\mathbf{Alg}) \ar[dr]_-{u}
  \ar[rr]^-{w} & &\mathbf{Poiss}^{\mathrm{eq}}_{\C, n} \ar[dl]^-{q}
  \\ & \mathcal{I}(\mathbf{CAlg}_{\C}) & }
$$
Note that $u$, and $q$
are both surjective, since they both have a section given by the
trivial bracket or strict shifted Poisson structure. Moreover, the
fiber of $q$ at $A \in \mathbf{CAlg}_{\C}$ is obviously equivalent to
the space $\mathsf{Poiss}_{\C}(A,n)$ of $n$-shifted Poisson structures
on $A$ (Definition \ref{spoissdef}). The following is a slight
enhancement of \cite[Theorem 3.2]{mel}.

\begin{thm}\label{melext}
  The map of spaces
  $w: \mathcal{I}(\mathbb{P}_{\C, n+1}-\mathbf{Alg})
  \to \mathbf{Poiss}^{\mathrm{eq}}_{\C, n} $ is an equivalence.
\end{thm}
\noindent \textbf{Proof}.  It is enough to prove that for any
fibrant-cofibrant $A \in \mathsf{CAlg}_{\mathsf{M}}$, the map induced
by $w$ between the $u$ and $q$ fibers over $A$ is an equivalence. But,
since these fibers are, by definition, $\mathbb{P}_{\C, n+1}(A)$ (the
moduli space of $\mathbb{P}_n$-algebra structures on $A$ having the
given underlying commutative algebra structure), and
$\mathsf{Poiss}_{\C}(A,n)$ (the space of $n$-shifted Poisson structure
on $A)$, this is exactly \cite[Theorem 3.2]{mel}. \hfill $\Box$

\

\smallskip

\noindent
For future reference, we also state the following immediate
consequence.
\begin{cor}\label{mel}  \emph{(Melani)}
  For any $A\in \mathbf{CAlg}_{\C}$, the map $w$ of Theorem
  \ref{melext} induces an equivalence \linebreak
  $w_A: \mathbb{P}_{\C, n+1}(A)
  \simeq \mathsf{Poiss}_{\C}(A,n)$ between the moduli space of
  $\mathbb{P}_n$-algebra structures on $A$ (having the given
  underlying commutative algebra structure), and the moduli space of
  $n$-shifted Poisson structures on $A$.
\end{cor}
\noindent \textbf{Proof}. As mentioned above, this is exactly
\cite[Theorem 3.2]{mel}. \hfill $\Box$

\

\smallskip

\begin{rmk}
  (1) Theorem 3.2 in \cite{mel} is stated for the model category
  $\mathsf{M}$ of bounded above cochain complexes of $k$-modules, but
  the proof is general and it extends immediately to our general
  $\mathsf{M}$. The original statement seems moreover to require a
  restriction to those cdga's having a dualizable cotangent
  complex. This is due to the fact that the author uses the tangent
  complex (i.e. the dual of the cotangent complex) in order to
  identify derivations. However, the actual proof produces an
  equivalence between (weak, shifted) Lie brackets and (weak)
  biderivations. Therefore if one identifies derivations using the
  linear dual of the symmetric algebra of the cotangent complex, the
  need to pass to the tangent complex disappears, and the result
  holds, with the same proof and without the assumption of the
  cotangent complex being dualizable (see Remark \ref{yok}). This is
  the main reason we adopted Definition \ref{polydef} and
  \ref{Polydef} as our definition of polyvectors.

  \

  \noindent
  (2) Melani's proof of \cite[Theorem 3.2]{mel} also shows that the
  natural map from the classifying space $\mathcal{I}(\mathbb{P}_{\C,
    n+1}-\mathbf{Alg})$ of strict $\mathbb{P}_{n+1}$-algebras in
  $\mathsf{M}$ to the classifying space $\mathcal{I}(\mathbb{P}_{\C,
    n+1, \infty}-\mathbf{Alg})$ of weak $\mathbb{P}_{n+1}$-algebras in
  $\mathsf{M}$ (where the operad $\mathbb{P}_{\C, n+1, \infty}$ is any
  cofibrant resolution of $\mathbb{P}_{\C, n+1}$ in the model category
  of operads in $\C$) is an equivalence.
\end{rmk}

\

\smallskip

\noindent
We also have the following, easier, strict analog of Corollary \ref{mel}.

\begin{prop} \label{melstrict} \emph{\cite[Proposition 1.4.8]{cptvv}}
  For any $A \in \mathsf{CAlg}_{\mathsf{M}}$, there is a natural bijection 
$$
\mathbb{P}_{\m ,n}(A) \simeq Hom_{Lie_{\mathsf{M}}^{gr}}(\mathbf{1}_{\mathsf{M}}(2)[-1],
\mathsf{Pol}^{\mathsf{M}}(A,n)[n])
$$
between the set of (strict) $\mathbb{P}_n$-algebra structures on $A$ in $\mathsf{M}$,
and the set of morphisms $\mathbf{1}(2)[-1] \to
\mathsf{Pol}^{\mathsf{M}}(A,n)[n]$  of Lie algebra
objects in $\mathsf{M}^{gr}$.
\end{prop}

As an immediate consequence of Theorem \ref{melext}, we get the
following, useful, strictification result

\begin{cor} Any weak shifted Poisson pair in $\mathsf{Poiss}^{\textrm{eq}}_{\mathsf{M}, n}$
  is equivalent to a strict pair.
\end{cor}
\noindent \textbf{Proof}. By Theorem \ref{melext}, an object $(A, \pi)
\in \mathsf{Poiss}^{\textrm{eq}}_{\mathsf{M}, n}$ (i.e. an a priori
weak pair), is equivalent to a pair of the form $w(B)$, where $B \in
\mathbb{P}_{\mathsf{M}, n+1}-\mathbf{Alg}$ (i.e. is a strict
$\mathbb{P}_{n+1}$-algebra in $\mathsf{M}$), whose underlying
commutative algebra is weakly equivalent to $A$ in
$\mathsf{CAlg}_{\mathsf{M}}$. We conclude by observing that objects in
the image of $w$ are always strict pairs. \hfill $\Box$

\

\smallskip

\noindent \textbf{The $\DR$-to-$\Pol$ construction.}  Let $A' \in
\mathbb{P}_{\mathsf{M}, n}-\mathbf{Alg}.$ By Proposition
\ref{melstrict}, the $ \mathbb{P}_{n}$ algebra structure on the
underlying commutative algebra of $A'$ is encoded by a strict map $\pi
: \mathbf{1}_\mathsf{M}[-1](2) \to \pol^{\mathsf{M}}(A',n)[n]$ of
graded Lie algebras in $\mathsf{M}$. Since the weight $q$ part of
$\pol^{\mathsf{M}}(A',n)$ is, by definition, $\ti^{(q)}(A',
-n)^{\Sigma_q}$, $\pi$ induces a map \linebreak $\mathbf{1}_M \to \ti^{(2)}(A',
-n)^{\Sigma_2}[n+1]$ that we still denote by $\pi$. Write
$$
[ -,-]_{p,q}: \pol^{\mathsf{M}}(A',n)(p) \otimes_{\mathsf{M}}
\pol^{\mathsf{M}}(A',n)(q) \to \pol^{\mathsf{M}}(A',n)(p+q-1)[-n]
$$
for
the Lie bracket part of the graded $\mathbb{P}_{n+1}$-algebra
structure on $\pol^{\mathsf{M}}(A',n)$, the family of composite maps
in $\mathsf{M}$
$$
\xymatrix{
\ti^{(q)}(A',-n) \ar[drrr]_-{\epsilon_{q}} & \hspace{-2.5pc} \simeq & \hspace{-2.5pc} 
\mathbf{1}_M \otimes_{\mathsf{M}} \ti^{(q)}(A',-n)
  \ar[r]^-{\pi \otimes \mathrm{id}} & \ti^{(2)}(A',
  -n)^{\Sigma_2}[n+1] \otimes_{\mathsf{M}}
  \ti^{(q)}(A',-n)^{\Sigma_{q}} \ar[d]^-{[-,-]_{2,q}[n+1]} \\
& & & 
  \ti^{(q+1)}(A', -n)^{\Sigma_{q+1}}[1], }
$$
for $q\in \mathbb{N}$, is
easily verified to be the mixed differential of a mixed graded algebra
structure on $ \pol^{\mathsf{M}}(A',n)$ inside
$\mathsf{M}$.\\ Moreover, since $\pol^{\mathsf{M}}(A',n)(0)= A'$, by
the universal property of $\dr^{\mathsf{M}}(A')$ the identity $A' \to
A'$ induces a map $$\phi_{A', \pi}: \dr^{\mathsf{M}}(A')
\longrightarrow \pol^{\mathsf{M}}(A',n)$$ of mixed graded algebras in
$\mathsf{M}$.

\begin{rmk}
  The above construction of the mixed differential $\epsilon=
  (\epsilon_q)_{q}$ is a generalization of the classical construction
  associating to a classical Poisson bivector $\pi$ on a smooth
  manifold, the differential $d_{\pi}:=[\pi, -]$ on polyvectors, where
  $[-,-]$ is the Schouten-Nijenhuis bracket, and $d^2_{\pi}=0$ is equivalent to
  the bivector $\pi$ being Poisson.
\end{rmk}

\

 A slight elaboration of this construction (by choosing strict models,
 see \cite[1.4.3]{cptvv} for details), yields the following derived
 version \begin{itemize}
 \item functors $\DR^{\C}, \Pol^{\C} (-, n) :
   (\mathbb{P}_{\C , n+1}-\mathbf{Alg})^{eq} \to \ealgrM$
\item a morphism $\Phi: \DR^{\C} \to \Pol^{\C} (-, n)$ betwen the above functors.
\end{itemize} 
 Here $(\mathbb{P}_{\mathcal{M}, n+1}-\mathbf{Alg})^{eq}$ is the
 $\infty$-category of $\mathbb{P}_{n+1}$-algebras in $\C$ with only
 equivalences as morphisms (this ensures that $\Pol^{\C} (-, n)$ is
 indeed a functor), and we slightly abuse notation by writing
 $\DR^{\C}$ for the composition of the usual $\DR^{\C}$ with the
 forgetful functor $(\mathbb{P}_{\C, n+1}-\mathbf{Alg})^{eq} \to
 \mathbf{CAlg}_{\C}$.

 \subsubsection{From non-degenerate
   Poisson algebra structures to symplectic structures} Classically,
 one way of stating that a Poisson structure (on a
 smooth scheme or manifold $X$) is non degenerate is by declaring that
 the analog of the above map $\Phi$ establishes an isomorphism of
 mixed graded algebras between the de Rham algebra and the algebra of
 polyvectors. Analogously, we say that $A \in \mathbb{P}_{\mathcal{M},
   n+1}-\mathbf{Alg}$ is \emph{non degenerate} if $\Phi_{A} :
 \DR^{\C}(A) \to \Pol^{\C} (A, n)$ is an equivalence of mixed graded
 algebras in $\C$. For such an $A$, by Corollary \ref{mel}, we also
 get a map $\alpha_{A}: \mathbf{1}_{\C}(2) \to \Pol^{\C} (A, n)[n+1]$
 in $\eMgr$. By putting these  together, we get a diagram
 $$
 \xymatrix{
   \DR^{\C}(A)[n+1] \ar[r]^-{\Phi_{A}} & \Pol^{\C} (A, n)[n+1] &
   \mathbf{1}_{\C}(2) \ar[l]_-{\alpha_{A}} }
 $$
 in $\eMgr$, exhibiting
 both $\DR^{\C}(A)[n+1]$ and $\mathbf{1}_{\C}(2)$ as objects in the
 overcategory $\eMgr/\Pol^{\C} (A, n)[n+1]$. We can then give the
 following
 \begin{df}
   Let $A \in \mathbb{P}_{\mathsf{M}, n+1}-\mathbf{Alg}$. The
   \emph{space of closed $2$-shifted forms compatible with} the given
   $\mathbb{P}_{n}$\emph{-structure on} $A$ is the space of lifts of
   $\alpha_A$ along $\Phi_A$, i.e. the space $$\mathsf{Comp}_{\C}^{2,
     cl}(A,n) := \mathsf{Map}_{\eMgr/\Pol^{\C} (A,
     n)[n+1]}(\mathbf{1}_{\C}(2), \DR^{\C}(A)[n+1]).$$
\end{df}
The reasons for using the words ``closed forms'' in the previous
definition are the following. First of all there is a ``forgetful''
map $$\mathsf{Comp}_{\C}^{2, cl}(A,n) \longrightarrow
\mathsf{Map}_{\eMgr}(\mathbf{1}_{\C}(2), \DR^{\C}(A)[n+1]) =
\mathcal{A}_{\C}^{2, cl}(A,n-1)$$ to the actual space of closed
$2$-shifted forms on $A$. Moreover, if $A$ happens to be \emph{non
  degenerate}, then $\mathsf{Comp}_{\C}^{2, cl}(A,n-1)$ is
contractible (since $\Phi_A$ is an equivalence), hence there is a
unique closed $(n-1)$-shifted $2$-form $\omega_A \in
\pi_0(\mathcal{A}_{\C}^{2, cl}(A,n-1))$ on $A$, via the above
``forgetful'' map. Moreover, by definition, such an $\omega_A$ is non
degenerate since $A$ is: $\omega$ is thus a symplectic structure on
$A$. If we denote by $\mathbb{P}^{nd}_{\C, n}(A)$ the subspace of
$\mathbb{P}_{\C, n}(A)$ whose connected components consist of non
degenerate elements, we get the following
\begin{prop}\label{W}
  For $A\in \mathbf{CAlg}_{\C}$ the above construction yields a
  well-defined map of spaces \ \qquad
  \linebreak 
  $W_A: \mathbb{P}^{nd}_{\C, n}(A) \to \mathsf{Sympl}_{\C}(A, n-1)$.
\end{prop}
There is  a parallel (and in fact  equivalent) notion of non-degenerate
shifted Poisson structure in $\C$. Let $B \in \mathbf{CAlg}_{\C}$ such
that $\mathbb{L}^{\C}_{B}$ is a dualizable in $B-\mathbf{Mod}_{\C}$. An
$n$-shifted Poisson structure $\pi \in \pi_0
(\mathsf{Poiss}_{\C}(B,n))$ defines a map $\mathbf{1}_{\C} \to
Sym_B^{2}(\mathbb{T}^{\C}_{B}[-n-1])[n+2]$ in $\C$, and thus induces,
by adjunction, a map $\mathbb{L}^{\C}_B \to \mathbb{T}^{\C}_B [-n]$ in
$B-\mathbf{Mod}_{\C}$: we say that $\pi$ is \emph{non degenerate} if
this map is an equivalence. The subspace of $\mathsf{Poiss}_{\C}(B,n)$
whose connected components consist of non degenerate elements will be
denoted by $\mathsf{Poiss}^{nd}_{\C}(B,n)$. The notions of non
degeneracy for Poisson algebras and for Poisson structures coincide in
the following sense

\begin{prop}\label{W'}
  Let $B \in \mathbf{CAlg}_{\C}$ such that $\mathbb{L}^{\C}_{B}$ is a
  dualizable in $B-\mathbf{Mod}_{\C}$. The equivalence $w_B$ of
  Corollary \ref{mel}, restricts to an equivalence
  $\mathbb{P}^{nd}_{\C, n}(B)\simeq
  \mathsf{Poiss}^{nd}_{\C}(B,n-1)$. As a consequence of Proposition
  \ref{W}, we thus get a map $$W_B': \mathsf{Poiss}^{nd}_{\C}(B,n-1)
  \longrightarrow \mathsf{Sympl}_{\C}(B, n-1)$$
\end{prop}

The map $W'_B$ in Proposition \ref{W'} is called the \emph{comparison
  map} between non degenerate shifted Poisson structures and shifted
symplectic structures. We will study it for derived Artin stacks in
Section \ref{poiss-sympl}.

\begin{rmk} It is possible that  the map $W'_B$ in Proposition \ref{W'}
  will turn out to be an equivalence, for very general $\C$. We have
  proven this in our geometric case of interest (see Theorem
  \ref{comparison}). The difficulty in general stems from the fact
  that while for the source of $W_{B}'$ we have Theorem \ref{mel} and
  Proposition \ref{melstrict}, and we thus are able to perform the
  $\DR$-to-$\Pol$ construction and produce the map $W'_{B}$, we don't
  have anything similar for the target of $W'_{B}$. One runs into
  serious difficulties even just trying to construct an inverse
  equivalence to $W'_{B}$ at the level of connected components. The
  first obstacle is that a shifted symplectic structure is, by
  definition, a weak map (i.e. a map in the relevant homotopy
  category). Even if we could strictify this map (getting a strictly
  closed shifted $2$-form), we are still left with the problem that
  non degeneracy is a weak property, i.e. the property of a map being
  a quasi-isomorphism, and therefore cannot be readily used to build a
  strict Lie bracket on the de Rham algebra. One might be able to
  overcome these difficulties in general but we do not know how to do
  this at the moment.
\end{rmk}

\

\subsubsection{Tate realizations over $k$}\label{tate}
In this section, we will concentrate on the special case $\m=
\epsilon-\mathsf{dg}^{gr}_k$, with its associated $\s$-category $\C=
\medg$, and in the next Section we will apply the definitions and
results obtained here to categories of diagrams in $\C$.

The unit of the  symmetric monoidal category $\m$ is $k(0)$,
i.e. the complex $k[0]$ sitting in pure weight $0$, with the
trivial mixed differential. The enriching symmetric monoidal functor
is given by $\mathsf{dg}_k \to \m \, : V \mapsto V(0)$, the tensor
enrichment is then given by $V\otimes E:= V(0)\otimes_{\m} E$, for
$E\in \m$. The enriched hom object is thus
$$
\underline{Hom}_{k}(E,E') :=
\mathrm{Z}_{\epsilon} (\underline{Hom}_{\C}(E,E')(0)) \in
\mathsf{dg}_k \, ,
$$
where $\underline{Hom}_{\m}$ denotes the internal
Hom in $\m$, and, for $F \in \m$, we write
$\mathrm{Z}_{\epsilon}(F(0)) \in \dg_{k}$ for the kernel of the map of
$k$-dg-modules $\epsilon : F(0) \to F(1)[1]$.

The standard realization functor $|-|: \C \to \dg_k$, for $\C= \medg$,
is in some sense unsatisfactory since an easy computation
(\cite[Proposition 1.5.1]{cptvv}) shows that $|E| \simeq \prod_{p\geq
  0}E(p)$, for $E \in \C$, so that all negative weights are lost under
standard realization. An obvious way to modify $|-|$ and
remedy this flaw is to consider $|E|^{t}:= \textrm{colim}_{i\geq
  0}\prod_{p\geq -i}E(p)$, instead. This new functor $|-|^t: \C \to
\dg_k$ will be called the \emph{Tate realization functor} for $\C$. By
definition, there is a canonical morphism $|-| \to |-|^t$ of
$\s$-functors. One can show, exactly as for $|-|$, that $|-|^t$ is lax
symmetric monoidal as well, so that it is inherited by  categories of
algebras. In particular, we also get \emph{Tate realization
  functors}
$$
|-|^t : \epsilon-\mathbf{CAlg}^{gr}_{\C} \to
\epsilon-\mathbf{cdga}^{gr}_k\, , \,\,\,\, |-|^t :
\mathbb{P}^{gr}_{\C, n}-\mathbf{Alg} \to
\mathbb{P}_{n}-\mathbf{cdga}^{gr}_k \,.
$$
As in the linear case, there are canonical morphisms $|-| \to
|-|^t$ of $\s$-functors between realizations on the level of these
algebra structures.

\  

Let us put this into a broader perspective, and
relate the Tate realization to a standard realization (on a different
category). Let us start by the fact that there is an equivalence in
$\dg_k$
$$
\underline{\mathcal{H}om}_{k}(k(i), k(i+1)) =
\mathbb{R}\underline{Hom}_{k}(k(i), k(i+1)) \simeq k[0]\,
$$
(where
$\underline{\mathcal{H}om}_{k}$ denotes the $\dg_k$-enriched Hom in
$\C$), giving rise to the following canonical ind-object in
$\C$
$$
k(\s) := \{ k(0) \to k(1) \to \dots \to k(i) \to k(i+1) \to
\cdots \} \, \in \mathbf{Ind}(\C).
$$
One can then show that the
\emph{Tate} realization functor for $\C$ is related to the
\emph{standard} realization $|-|_{\mathbf{Ind}}$ for
$\mathbf{Ind}(\C)$, by
$$
|-|^t \simeq |\, - \otimes k(\s) \,
|_{\mathbf{Ind}} \, : \C \longrightarrow \dg_k \,
$$
where we have
implicitly used the canonical functor $\C \to \mathbf{Ind}(\C)$,
sending an object $E$ to the constant ind-diagram in $\C$ with value
$E$.

Moreover, since $k(i) \otimes k(j) \simeq k(i+j)$, $k(\s)$ is
a commutative monoid object in $\mathbf{Ind}(\C)$, hence $$A\in
\mathbf{CAlg}_{\C} \,\, \Rightarrow \, \, A(\s) := A\otimes k(\s) \in
k(\s)/\mathbf{CAlg}_{\ind}.$$ Therefore it make sense to consider the
relative objects
$$
\DR^{\ind}(A(\s)/k(\s)) \in
\epsilon-\mathbf{Alg}^{gr}_{\ind} \, , \,\,\, \,
\mathbf{Pol}^{\ind}(A(\s)/k(\s), n) \in \mathbb{P}_{\ind,
  n+1}-\mathbf{Alg}^{gr} \, ,
$$
and their standard realizations (on
algebras in $\ind$)
$$
\DR(A(\s)/k(\s)) \in
\epsilon-\mathbf{cdga}^{gr}_{k} \, , \,\,\, \,
\mathbf{Pol}(A(\s)/k(\s), n) \in
\mathbb{P}_{n+1}-\mathbf{cdga}_{k}^{gr} \, .
$$
For de Rham algebras
and polyvectors, we have the following comparison result
\begin{prop}\label{dunno}
  If $A \in \mathbf{CAlg}_{\C}$, then we have a canonical
  equivalence $$\DR^t (A):= |\DR (A)|^t \simeq \DR(A(\s)/k(\s))$$ in
  $\epsilon-\mathbf{cdga}^{gr}_{k}$. If moreover $\mathbb{L}^{\C}_{A}$
  is dualizable in $A-\mathbf{Mod}_{\C}$, then we have a canonical
  equivalence $$\mathbf{Pol}^t (A, n) := | \mathbf{Pol}(A, n) |^t
  \simeq \mathbf{Pol}(A(\s)/k(\s)) $$in
  $\mathbb{P}_{n+1}-\mathbf{cdga}_{k}^{gr}$.
\end{prop}

\begin{rmk}\label{sometimesiso}
As already observed in the linear case, in general, none of the
morphisms $\DR(A) \to \DR^t(A)$, $\mathbf{Pol} (A, n) \to
\mathbf{Pol}^t (A, n)$ are equivalences. If $A \in
\epsilon-\mathbf{CAlg}^{gr}_{\C}$ happens to have only non-negative
weights (this will be the case in our application to derived stacks),
then also $\mathbb{L}_{A}$ will have the same property, and $\DR(A)
\to \DR^t(A)$ will indeed be an equivalence. However, the dual to
$\mathbb{L}_{A}$, if existing, will have postive weights, so that
$\mathbf{Pol} (A, n) \to \mathbf{Pol}^t (A, n)$ will \emph{not} be an
equivalence, even in this case. So, at least for applications to
derived stacks, while the introduction of the Tate realization will
not make any difference for $\DR$, it will definitely do for
$\mathbf{Pol}$, and indeed the interesting realization will be
$\mathbf{Pol}^t$ rather than the standard one.
\end{rmk}

\subsection{Formal derived stacks and formal localization}\label{FL}

A crucial ingredient in the theory of shifted Poisson structures on
general derived Artin stacks is the method of  \emph{formal
  localization}.  Formal localization is interesting
in its own right as a new, very powerful tool that will prove useful
in many other situations inside derived geometry, especially in order
to globalize tricky constructions and results, whose extension from
the local case presents obstructions that only vanish formally
locally. An example is given by obstructions living in de Rham
cohomology (even, say, on a smooth scheme $X$). Suppose we wish to
glue some construction that we can perform ``locally'' on $X$, and we
know that obstructions to globalize live in de Rham cohomology of $X$
(e.g. we would like to globalize an algebraic version of the Darboux
lemma). Since de Rham cohomology never vanishes Zariski or \'etale
locally, it is going to be hard or impossible, depending on the
specific problem, to suppress the obstructions and glue with respect
to these topologies. On the other hand, for any $x \in X$, the de Rham
cohomology of the formal completion $\widehat{X_{x}}$ does vanish,
so we may try to glue the construction performed on the family
$\{\widehat{X_{x}}\, |\, x\in X \}$ to a construction on $X$. But in
order to do this we need a result telling us when and how we are able
to glue objects along the family of formal completions. This is exactly the
content of formal localization and below we will sketch how it works and
why it is useful. 

\

\noindent The following conventions will be adopted throughout
this section. A cdga $A$ is called \emph{almost finitely presented} if
$H^0(A)$ is a $k$-algebra of finite type, and each $H^{i}(A)$ is a
finitely presented $H^{0}(A)$-module. We will write  $\dAff_{k}$ for 
the opposite $\s$-category of almost finitely presented cdga's, and we
will simply refer to its objects as \emph{derived affine schemes}
without mentioning the finite presentation condition. In particular,
when writing $\Spec\, A$, we implicitly assume that $\Spec\, A$ is an
object of $\dAff_{k}$. The $\s$-category $\dAff_{k}$ is equipped with
its usual \'etale topology of \cite[Definition 2.2.2.3]{hagII}, and the
corresponding $\s$-topos of derived stacks will be denoted by
$\dSt_{k}$. Its objects will simply be called \emph{derived stacks},
instead of the more precise but longer \emph{locally almost finitely
  presented derived stacks over $k$}. With these conventions, an
algebraic derived $n$-stack $X$ will have a smooth atlas by objects in
$\dAff_{k}$, and in particular, it will have a bounded above cotangent
complex in $\mathsf{Coh}(X)$.\\

\subsubsection{Formal derived stacks}

As the name suggests, formal localization deals with \emph{formal
  derived stacks}, which we now define.

\begin{df}\label{d10-}
A \emph{formal derived stack} is
a derived stack $F$
satisfying the following conditions.

\begin{enumerate}

\item $F$ is \emph{nilcomplete} i.e. for all $\Spec\, A \in \dAff_{k}$, the canonical map 
$$F(A)\longrightarrow \lim_{k}F(A_{\leq k})$$
(induced by the Postnikov tower of $A$) is
an equivalence in $\T$.

\item $F$ is \emph{infinitesimally cohesive}
  i.e. for all cartesian squares of almost finitely presented
  $k$-cdga's in non-positive degrees
$$\xymatrix{
  B \ar[r] \ar[d] & B_{1} \ar[d] \\
  B_{2} \ar[r] & B_{0},}$$ such that each $\pi_{0}(B_i)
\longrightarrow \pi_{0}(B_0)$ is surjective with nilpotent kernel,
then the induced square
$$\xymatrix{
F(B) \ar[r] \ar[d] & F(B_{1}) \ar[d] \\
F(B_{2}) \ar[r] & F(B_{0}),}$$
is cartesian in $\T$.
\end{enumerate}
\end{df}

\begin{rmk} 
  (1) The property of being infinitesimally cohesive is a derived
  variation of the Schlessinger condition in classical deformation
  theory (\cite{sch}). In particular, one can show that \emph{any
    derived Artin stack} $F$ \emph{is a formal derived stack}
  (\cite[Appendix]{hagII}), and it is actually \emph{cohesive}
  i.e. sends any diagram as in \ref{d10-} $(2)$, with the nilpotency
  condition possibly omitted, to pullbacks in $\T$ (\cite[DAG IX,
    Corollary~6.5]{lu2} and \cite[DAG XIV, Lemma 2.1.7]{lu2}).

  \

  \noindent
  (2) A small limit of formal derived stacks is a formal derived
  stack.
\end{rmk}

\

\smallskip

\noindent
Let the $\s$-functor $i: \mathbf{alg}^{\textrm{red}}_{k} \longrightarrow
\cdga_{k}$ be the inclusion of the full reflective sub $\s$-category
of \emph{reduced discrete} objects (i.e. $R\in \cdga_{k}$ such that
$R$ is discrete and $R \simeq H^0(R)$ is a usual reduced $k$-algebra).
The $\s$-functor $i$ has a left adjoint
$$(-)^{red}: \cdga_{k} \longrightarrow
\mathbf{alg}^{\textrm{red}}_{k} \, , \, \, A \longmapsto
A^{red}:=H^0(A)/\textrm{Nilp}(H^0 (A),$$ 
and it is easy to
verify that we get an induced $\s$-functor 
$i^*: \dSt_{k} \longrightarrow
\mathbf{St}_{\textrm{red}, k}$,  
where $\mathbf{St}_{\textrm{red}, k}$ is the $\s$-category of stacks on
$(\mathbf{alg}^{\textrm{red}}_{k})^{op}$ for the \'etale topology.
Now $i^*$ has both a right adjoint $i_*$,
and a left adjoint $i_!$, both fully faithful, and  $i_!i^*$ is left adjoint to $i_*i^*$.

\begin{df}\label{red&dR} \begin{enumerate}
\item The functor 
$(-)_{DR} := i_*i^* : \dSt_{k} \longrightarrow
  \dSt_{k}$
is called the \emph{de Rham stack functor}. By
  adjunction,  for any $F\in \dSt_{k}$, we have a canonical natural
  map $q_F: F \to F_{DR}$. 
\item The functor 
$(-)_{\textrm{red}} := i_! i^* : \dSt_{k} \longrightarrow \dSt_{k}$ 
is called the \emph{reduced stack functor}. By adjunction,  for any
$F\in \dSt_{k}$, we have a canonical natural map $\iota_{F} :
F_{\textrm{red}} \to F$.  
\item Let $f : F \longrightarrow G$ be a morphism
in $\dSt_{k}$. We define the \emph{formal completion} $\widehat{G}_{f}$ 
\emph{of} $G$ \emph{along the morphism} $f$ as the fibered product in
$\dSt_{k}$: 
$$\xymatrix{\widehat{G}_{f} \ar[r] \ar[d] &
  F_{DR} \ar[d]^-{f_{DR}} \\ G \ar[r]_-{q_{G}} & G_{DR}.}$$  
\end{enumerate}
\end{df}

\

\noindent Since the left adjoint to $i$ is $(-)^{red}$, then it is
easy to see that 
$F_{DR}(A) \simeq F(A^{\red})$, and  $(\Spec \, A)_{\textrm{red}} \simeq \Spec
\, (A^{\red}),$ 
for any $A\in \cdga_{k}$. Therefore
$\widehat{G}_{f}(A)=G(A) \times_{G(A^{red})}F(A^{red})$, for $f: F \to
G$ in $\dSt_{k}$. We already observed that $(-)_{DR}$ is right adjoint
to $(-)_{\textrm{red}}$, as functors $\dSt_{k}\to \dSt_{k}$. We list a few elementary
properties of these constructions:  


%

\begin {itemize}
\item $F_{DR}$ is a formal derived stack for any $F \in \dSt_{k}$.
\item If $G$ is a formal derived stack, the formal completion
  $\widehat{G}_{f}$, along any map $f:F\to G$ in $\dSt_{k}$, is again a
  formal derived stack. 
\item For any $F \in \dSt_{k}$, if $j: \mathsf{t}_0 F \to F$ denotes
  the canonical map in $\dSt_{k}$ from the truncation of $F$ to $F$,
  then the canonical map $\widehat{F}_{j} \to F$ is an equivalence.
\end{itemize}


Our main object of study in the next section, will be the map $q: F \to
F_{DR}$. If $\mathbb{K}$ is a field, and $x: \Spec \, \mathbb{K} \to
F_{DR}$ is a point (by adjunction, this is the same thing as a
$\mathbb{K}$-point of $F$, since $\mathbb{K}$ is reduced), then the
fiber of $q$ at $x$ is exactly the (classical) formal completion
$\widehat{X_x}$ of $X$ at the closed point $x$. More generally, we
have

\begin{prop}\label{fibersDR}
  Let $F \in \dSt_{k}$, $\Spec\, A \in \dAff_{k}$, and
  $\overline{u}:\Spec\, A \longrightarrow F_{DR}$, corresponding to a morphism $u : \Spec\,
  A^{red} \longrightarrow F$. Then the  base-change derived stack
  $F\times_{F_{DR}}\Spec\, A$ is equivalent to the formal completion
  $\widehat{(\Spec\, A \times F)}_{(i,u)}$ of the graph morphism
$$(i,u) : \Spec\, A^{red} \longrightarrow \Spec\, A \times F,$$
where $i : \Spec\, A^{red} \longrightarrow \Spec\, A$ is the natural map.
\end{prop} 

\begin{rmk}\label{easy}
  Suppose that $F \in \dSt_k$ has a cotangent complex (e.g. $F$ is a
  derived Artin stack). Then $\mathbb{L}_{F/F_{DR}}$ exists, and we
  have an equivalence $\mathbb{L}_{F} \simeq \mathbb{L}_{F/F_{DR}}$ in
  $\mathbf{QCoh}(F)$. In fact, by the transitivity sequence for the
  map $q: F \to F_{DR}$. it is enough to show that
  $\mathbb{L}_{F_{DR}} \simeq 0$. But this follows immediately from
  the equivalences: $F_{DR}(A \oplus E) \simeq F((A\oplus E)^{\red})
  \simeq F(A^{\red}) \simeq F_{DR}(A),$ for any $A \in \cdga_{k}^{\leq
    0}$, and any $E \in \dg^{\leq 0}_{A}$.
\end{rmk}

\subsubsection{Formal localization for $X \to X_{DR}$}\label{xtoxdr}
The general theory of formal localization is developed in detail in
Section 2 of \cite{cptvv}. Instead of giving a complete account, we
will content ourselves with the application of the general theory to
our main case of interest, i.e. to the map $q: X \to
X_{DR}$. Throughout this Section, $X$ will be derived Artin stack lfp
over $k$ (hence with a perfect cotangent complex).\\

\noindent \textbf{The map $X\to X_{DR}$ as a family of formal derived
  stacks.} First of all, let us observe that $q: X \to X_{DR}$ is an
\emph{algebraisable family of perfect formal derived stacks}, i.e. for
any derived point $\Spec \, A \to X_{DR}$ the corresponding fiber $X_A
\to \Spec \, A$ of $q$ has the following properties:
\begin{enumerate}
\item $X_A$ is a formal derived stack, and the canonical map
  $(X_A)_{\red} \to \Spec \, A_{\red}$ is an equivalence in $\dSt_k$.
\item The relative cotangent complex $\mathbb{L}_{X_A/A}$ is perfect
  (by Remark \ref{easy} and base change)
\item $X_A$ has a cohomologically bounded above coherent cotangent
  complex $\mathbb{L}_{X_A}$ (i.e. for any $x_B: \Spec \, B \to X_A$,
  the fiber $x_B^*\mathbb{L}_{X_A}$ is a cohomologically bounded above
  coherent $B$-dg-module).
\item $X_A$ is equivalent to the formal completion of $X \times
  \Spec\, A$ along the map $\Spec \, A_{\red} \to X \times \Spec\, A$,
  induced by the chosen derived point $\Spec \, A \to X_{DR}$, and the
  canonical map $\Spec\, A_{\red} \to \Spec\, A$ (see Proposition
  \ref{fibersDR}).
\end{enumerate}

Properties 1-3 caracterize a family of perfect formal derived stacks
over $A$, while property 4 says that the family is algebraisable
(\cite[Section 2.1, 2.2.]{cptvv}). \\

\

\smallskip

\noindent \textbf{Crystalline structure sheaf and shifted principal
  parts.}  Let us consider the $\s$-functor
$$
\D : \dAff^{op}_{k}
\longrightarrow \epsilon-\cdga^{gr}_{k}\, , \qquad  A \longmapsto
\DR(A_{\textrm{red}}/A)
$$
(where $\DR(A^{\textrm{red}}/A)$ is defined
in Definition \ref{d1.1}, via Remark \ref{relDR}). Note that $\D(A)
\simeq Sym_{A_{\red}}(\mathbb{L}_{A_{\red}/A}[-1])$ in $\cdga_k^{gr}$
(Section \ref{DEF}). The functor $\D$ satisfies \'etale descent, and
thus we have an induced $\s$-functor $$\D : \dSt^{op}_{k}
\longrightarrow \epsilon-\cdga^{gr}_{k}\, , \,\, F \longmapsto
\lim_{\Spec\, A \to F} \D(A).$$ We consider the following prestacks of
mixed graded cdga's on $\dAff/X_{DR}$ $$\begin{aligned}
  \D_{X_{DR}}:=\D(\OO_{X_{DR}}) : & (\dAff_{k}/X_{DR})^{op}
  \longrightarrow \mecdga_{k}, \qquad (\Spec\, A \rightarrow X_{DR})
  \longmapsto \D(A), \\ \mathcal{P}_{X} : & (\dAff_{k}/X_{DR})^{op}
  \longrightarrow \mecdga_k, \qquad (\Spec\, A \rightarrow X_{DR})
  \longmapsto \D(X_A).
\end{aligned}
$$

\noindent Note that there is a natural equivalence
$\mathcal{P}_{X}(\Spec\, A \rightarrow X_{DR}) \simeq Sym_{A_{\red}}
(\mathbb{L}_{\Spec \, A_{\red}/X_A} [-1])$ in $\cdga_k^{gr}$
(\cite[Proposition 2.2.6]{cptvv}).

\begin{df}\label{crys}
The prestack $\D_{X_{DR}}$ on $X_{DR}$ is called the \emph{crystalline
  structure sheaf} for $X$. The prestack $\mathcal{P}_X$ on $X_{DR}$
is called the prestack of \emph{principal parts} for $X$.
\end{df}

\

\begin{rmk}\label{Chev} An alternative interpretation  of
$\B_{X}$ can be given as follows. As already observed, the canonical
  map $ X \longrightarrow X_{DR}$ realizes $X$ as a family of formal
  derived stacks over $X_{DR}$, namely as the family of formal
  completions at closed points of $X$. By \cite{lu} these formal
  completions are determined by a dg-Lie algebra $\ell_x$. The dg-Lie
  algebra $\ell_x$ itself does not extend globally as a sheaf of
  dg-Lie algebras over $X_{DR}$, simply because its underlying complex
  is $\mathbb{T}_{X}[-1]$, the shifted tangent complex of $X$ (\cite{hen}), does
  not descend to $X_{DR}$. However, a remarkable consequence of
  derived formal localization (Theorem \ref{yeah}) is that the
  Chevalley-Eilenberg complexes of $\ell_x$, for $x \in X$, viewed as
  a graded mixed commutative dg-algebras, do glue to a global object
  over $X_{DR}$. This is exactly $\B_{X}$. Then, the Formal
  Localization Theorem \ref{yeah} tells us exactly how to express
  global geometric objects on $X$ as correspondingly sheafified
  objects on $X_{DR}$ related to $\B_X$.
\end{rmk}

Note that, by functoriality of $\D$, we have a natural morphism
$\D_{X_{DR}} \to \mathcal{P}_X$ of prestacks of mixed graded cdga's on
$X_{DR}$.  In particular, if we consider the $\s$-category $\C'$ of
prestacks on $\dAff_{k}/X_{DR}$ with values in
$\epsilon-\mathbf{dg}^{gr}_{k}$, then $\D_{X_{DR}} \in
\mathbf{CAlg}(\C')$, and $\mathcal{P}_X \in
\D_{X_{DR}}/\mathbf{CAlg}(\C')$. We let $\C:=
\D_{X_{DR}}-\mathbf{Mod}_{\C'}$.  Then $\mathcal{P}_X \in
\mathbf{CAlg}(\C) \simeq \D_{X_{DR}}/\mathbf{CAlg}(\C')$, and, for any
$n \in \mathbb{Z}$, we may consider (Section
\ref{diff})
$$
\mathbf{Pol}^{\C'}(\mathcal{P}_{X}/\D_{X_{DR}}, n)=
\mathbf{Pol}^{\C}(\mathcal{P}_{X}, n) \in \mathbb{P}^{gr}_{\C,
  n+1}-\mathbf{Alg} \,, \qquad 
\DR^{\C'}(\mathcal{P}_{X}/\D_{X_{DR}})= \DR^{\C}(\mathcal{P}_{X}) \in
\ealgrM \, .
$$
We will also consider the following prestacks on
$\dAff/X_{DR}$ obtained by Tate
realizations:
$$
\mathbf{\underline{Pol}}^{t}(\mathcal{P}_{X}/\D_{X_{DR}},
n): (\dAff_{k}/X_{DR})^{op} \longrightarrow
\mathbb{P}_{n+1}-\mathbf{cdga}_{k}^{gr} \, , (\Spec \, A \to X_{DR})
\mapsto | \mathbf{Pol}^{\medg}(\mathcal{P}_{X}(A)/\D_{X_{DR}}(A), n)
|^{t}$$ $$\underline{\DR}^{t}(\mathcal{P}_{X}/\D_{X_{DR}}) :
(\dAff_{k}/X_{DR})^{op} \longrightarrow
\epsilon-\mathbf{cdga}_{k}^{gr} \, , (\Spec \, A \to X_{DR}) \mapsto |
\DR^{\medg}(\mathcal{P}_{X}(A)/\D_{X_{DR}}(A)) |^{t}
$$

\

\begin{rmk}\label{2.2.8}
It is worth pointing out that while $\mathcal{P}_{X}$ and $D_{X_{DR}}$
are not stacks, all
$\mathbf{\underline{Pol}}^{t}(\mathcal{P}_{X}/\D_{X_{DR}}, n)$,
$\DR(\mathcal{P}_{X}/\D_{X_{DR}})$, and
$\underline{\DR}^{t}(\mathcal{P}_{X}/\D_{X_{DR}})$ are stacks
(\cite[Corollary 2.4.9]{cptvv}).
\end{rmk}

Analogously (see Section \ref{diff}), if we consider the $\s$-category
$\C'_{\textrm{Ind}}$ of prestacks on $\dAff_{k}/X_{DR}$ with values in
$\mathrm{Ind}(\epsilon-\mathbf{dg}^{gr}_{k})$, then (see Section
\ref{diff}) $\D_{X_{DR}}(\s) \in \mathbf{CAlg}(\C'_{\textrm{Ind}})$,
and $\mathcal{P}_X(\s) \in
\D_{X_{DR}}(\s)/\mathbf{CAlg}(\C'_{\textrm{Ind}})$. We let
$\C_{\textrm{Ind}}:=
\D_{X_{DR}}(\s)-\mathbf{Mod}_{\C'_{\textrm{Ind}}}$.  Then
$\mathcal{P}_X (\s)\in \mathbf{CAlg}(\C_{\textrm{Ind}}) \simeq
\D_{X_{DR}}(\s)/\mathbf{CAlg}(\C'_{\textrm{Ind}})$, and, for any $n
\in \mathbb{Z}$, we may
consider $$\mathbf{Pol}^{\C'_{\textrm{Ind}}}(\mathcal{P}_{X}(\s)/\D_{X_{DR}}(\s),
n)= \mathbf{Pol}^{\C_{\textrm{Ind}}}(\mathcal{P}_{X}(\s), n) \in
\mathbb{P}^{gr}_{\C_{\textrm{Ind}}, n+1}-\mathbf{Alg} \, ,$$ $$
\DR^{\C'_{\textrm{Ind}}}(\mathcal{P}_{X}(\s)/\D_{X_{DR}}(\s))=
\DR^{\C_{\textrm{Ind}}}(\mathcal{P}_{X}(\s)) \in
\epsilon-\mathbf{CAlg}_{\C_{\textrm{Ind}}}^{gr}.$$ And we also have
the following prestacks on $\dAff/X_{DR}$ obtained by standard
realizations:
$$
\xymatrix@R-1.7pc{
  \mathbf{\underline{Pol}}(\mathcal{P}_{X}(\s)/\D_{X_{DR}}(\s), n): (\dAff_{k}/X_{DR})^{op} \ar[r] &  \mathbb{P}_{n+1}-\mathbf{cdga}_{k}^{gr}, \\
  (\Spec \, A \to X_{DR}) \ar@{|->}[r] &
  | \mathbf{Pol}^{\mathrm{Ind}(\medg)}(\mathcal{P}_{X}(\s)(A)/\D_{X_{DR}}(\s)(A), n) |
}
$$   
and
$$
\xymatrix@R-1.7pc{
  \underline{\DR}(\mathcal{P}_{X}(\s)/\D_{X_{DR}}(\s)) : (\dAff_{k}/X_{DR})^{op} \ar[r] &
  \epsilon-\mathbf{cdga}_{k}^{gr}\\
  (\Spec \, A \to X_{DR}) \ar@{|->}[r] & | \DR^{\mathrm{Ind}(\medg)}(\mathcal{P}_{X}(\s)(A)/\D_{X_{DR}}(\s)(A)) |
}
$$

\

\begin{rmk} By Proposition \ref{dunno} and Remark \ref{sometimesiso},
  we get equivalences
  $$
  \begin{aligned}
  \DR^{\C'}(\mathcal{P}_{X}/\D_{X_{DR}}) & \simeq
  \underline{\DR}^{t}(\mathcal{P}_{X}/\D_{X_{DR}}) \simeq
  \underline{\DR}(\mathcal{P}_{X}(\s)/\D_{X_{DR}}(\s))\,, \\[-0.1pc]
 \mathbf{\underline{Pol}}^{t}(\mathcal{P}_{X}/\D_{X_{DR}}, n)
  & \simeq \mathbf{\underline{Pol}}(\mathcal{P}_{X}(\s)/\D_{X_{DR}}(\s),
 n)
\end{aligned}
  $$
  but notice that
  $\mathbf{Pol}^{\C'}(\mathcal{P}_{X}/\D_{X_{DR}}, n)$ is not in
  general equivalent to the previous ones.
\end{rmk}

\

\smallskip

\noindent \textbf{The Formal Localization theorem.}  We have already
defined the mixed graded $k$-cdgas $\DR(X/k)$, and $\DR(X/X_{DR})$
(Definition \ref{DRglobal}). It is an easy consequence of the
equivalence $\mathbb{L}_X \simeq \mathbb{L}_{X/X_{DR}}$ (Remark
\ref{easy}), that $\DR(X/k)\simeq \DR(X/X_{DR})$ in $\mecdga_k$. We
can give a similar, general definition of shifted polyvectors on $X$,
\emph{at least as a graded} $k$-\emph{cdga}.

\begin{df}\label{polstacky} Let $F \to G$ be a map between derived stacks, both
  having cotangent complexes (so that $\mathbb{L}_{F/G}$ exists,
  too). We define the graded $k$-dg module of \emph{$n$-shifted
    relative polyvectors} as
  $$
  Pol(F/G,n) \simeq \bigoplus_{p\geq 0}
  (\underline{Hom}_{\mathbf{QCoh}(F)}(\otimes^{p}
  \mathbb{L}_{F/G},\mathcal{O}_F[-pn]))^{h\Sigma_p} 
  \in \cdga^{gr}_k \, .
  $$
\end{df}
 
In the above definition, $\mathbf{QCoh}(F)$ is regarded as a
dg-category over $k$, and $\underline{Hom}_{\mathbf{QCoh}(F)}$ denotes
its $k$-dg-module of morphisms; the $\Sigma_p$-action on $
\underline{Hom}_{\mathbf{QCoh}(F)}(
\otimes^{p}\mathbb{L}_{F/G},\mathcal{O}_X[pn])$ is induced by
$\Sigma_p$ acting in the standard way on
$\otimes^{p}\mathbb{L}_{F/G}$, and via $(-1)^n$ times the sign
representation on $\OO_{F}[-pn] = \OO_{F}[-n]^{\otimes^{p}}$.

\begin{rmk}\label{natural}
  (1) Again by Remark \ref{easy}, we have an equivalence
  $Pol(X/X_{DR}, n) \simeq Pol(X/k, n)$ in $\cdga^{gr}_k$.

  \

  \noindent
  (2) When $\mathbb{L}_{F/G}$ is perfect over $F$ (e.g. for $F=X$
  derived Artin stack lfp over $k$, and $G=X_{DR}$), then we may
  express $Pol(F/G,n)$ using the dual relative tangent complex
  $\mathbb{T}_{F/G}$ as (see Remark \ref{yok0}) $$Pol(F/G, n) \simeq
  \bigoplus_{p\geq 0} \Gamma(F, Sym^{p}(\mathbb{T}_{F/G}[-n]) \in
  \cdga_k^{gr}$$
 \end{rmk}

The problem with Definition \ref{polstacky} is that, in general, it is
\emph{impossible} to directly endow $Pol(F(/G, n)$, as defined, with a
bracket and give it the structure of a graded
$\mathbb{P}_{n+1}$-differential graded algebras over $k$. This is
where the next result comes to rescue.

 \begin{thm} \label{yeah} \emph{(Formal localization for $X \to X_{DR}$)}
 
 Let $X$ be an Artin derived stack locally of finite presentation
  over $k$.
\begin{enumerate}
\item There is a natural equivalence 
of $\s$-categories
$$
\mathbf{Perf}(X) \simeq \mathcal{P}_X-\mathbf{Mod}_{\C}^{\mathsf{perf}}\, ,
$$
where $\C$ was defined right after Remark \ref{Chev}, and $\mathcal{P}_X-\mathbf{Mod}_{\C}^{\mathsf{perf}}$ is
the full sub-$\s$-category of $\mathcal{P}_X-\mathbf{Mod}_{\C}$, consisting of prestacks $E$ of graded mixed
$\mathcal{P}_X$-modules on $\dAff/X_{DR}$ satisfying the following two
conditions:
\begin{itemize}
\item For all $\Spec\, A \longrightarrow X_{DR}$, the graded mixed
  $\mathcal{P}_X(A)$-module 
$E(A)$ is equivalent, just as a graded $\mathcal{P}_X(A)$-module, to 
$\mathcal{P}_X(A) \otimes_{A_{\red}}E_{0}$, for some $E_{0} \in
\mathbf{Perf}(A_{\red})$.
\item $E$ is \emph{quasi-coherent} in the sense that: for all $\Spec\, B
  \longrightarrow \Spec\, A$ in $\dAff_{k}/X_{DR}$, the induced morphism
$
E(A) \otimes_{\mathcal{P}_X(A)}\mathcal{P}_X(B) \longrightarrow E(B)
$
is an equivalence.
\end{itemize}
\item There are natural equivalences of graded mixed cdga's over $k$
$$
\DR(X/X_{DR})\simeq \DR(X/k) \simeq
\Gamma(X_{DR},\underline{\DR}(\mathcal{P}_X/\D_{X_{DR}}))
\simeq \Gamma(X_{DR},\underline{\DR}^{t}(\mathcal{P}_X/\D_{X_{DR}})).
$$
\emph{(}where $\Gamma$ denotes derived global sections , i.e.
$\Gamma(X_{DR}, \mathcal{F})= \lim_{\Spec \, A \to X_{DR} } \mathcal{F}(A)$,
the limit being taken in the $\s$-category where $\mathcal{F}$ is valued\emph{)}.
\item For each $n\in \mathbb{Z}$, there are natural equivalences of
  graded dg-modules over $k$
$$
Pol(X/X_{DR},n)\simeq Pol(X/k,n) \simeq
\Gamma(X_{DR},\underline{\Pol}^t(\mathcal{P}_X/\D_{X_{DR}},n))
\simeq \Gamma(X_{DR},\underline{\Pol}(\mathcal{P}_X(\s)/\D_{X_{DR}}(\s),n)).
$$
\item Let $\C_{\mathrm{Ind}}$ be the $\s$-category defined right after Remark \ref{2.2.8}. The natural $\s$-functor 
$$
\mathcal{P}_X-\mathbf{Mod}_{\C}^{\mathsf{perf}} \longrightarrow
\mathcal{P}_X(\s)-\mathbf{Mod}_{\C_{\textrm{Ind}}}^{\mathsf{perf}}, 
$$ 
induced by the base change $(-)\otimes k(\s)$, 
is an equivalence.
\item For each $n\in \mathbb{Z}$, there are canonical equivalences in $\T$
  $$
  \mathsf{Symp}(X,n) \simeq \mathsf{Symp}(\B_X/\D_{X_{DR}}, n)
  \simeq \mathsf{Symp}(\B_X(\s)/\D_{X_{DR}}(\s),n)
  $$
  where $ \mathsf{Symp}(\B_X/\D_{X_{DR}}, n) $ and
  $\mathsf{Symp}(\B_X(\s)/\D_{X_{DR}}(\s),n)$ are defined  as in
  Definition \ref{symplabstrdef} \emph{(}with $\C$ and $\C_{\textrm{Ind}}$,
  respectively, as defined in the previous paragraph\emph{)}.
\end{enumerate}
\end{thm}
 The proof of Theorem \ref{yeah} can be found in \cite[Corollary 2.4.12, Proposition 2.4.15]{cptvv}. 

 \

 \begin{rmk}\label{whyFL}
   Let us explain why the formal localization Theorem \ref{yeah} is important and useful.

   \

   \noindent
   (1) Points 1 and 4 in Theorem \ref{yeah} are absolutely crucial since they 
   allow us to completely recover perfect complexes on $X$ as
   certain, explicitly identified, mixed graded modules over
   $\mathcal{P}_X$ or $\B_{X}(\s)$. This makes manipulations on
   perfect complexes much easier, and will allow us to go from a quantization of $\B_{X}(\s)$
   to a quantization of $\mathbf{Perf}(X)$ (see Section 3.1).

   \

   \noindent
   (2) Point 2 in the above Theorem lets us completely recover (with
   its full structure of mixed graded algebra) the descent-theoretic
   definition \ref{DRglobal} of the de Rham algebra on $X$ in terms of
   $\mathcal{P}_X/\D_{X_{DR}}$.

   \

   \noindent
   (3) Point 3 is important because it allows us not only to recover
   the 'geometrical'' polyvectors of Definition \ref{polstacky}, but
   also \emph{to establish a full graded}
   $\mathbb{P}_{n+1}$\emph{-algebra structure} on them. This is
   essential in order to be able to \emph{define} shifted Poisson
   structures on $X$ (Definition \ref{defPoiss}).

   \

   \noindent
   (4) The last point of Theorem \ref{yeah} lets us completely recover
   shifted symplectic forms on $X$ in terms of shifted symplectic form
   on $\mathcal{P}_X/\D_{X_{DR}}$ (or
   $\mathcal{P}_X(\s)/\D_{X_{DR}}(\s)$). This have the effect to
   enable a definition of shifted Poisson and shifted symplectic
   structures in terms of the \emph{very same} object
   $\mathcal{P}_X/\D_{X_{DR}}$ (or
   $\mathcal{P}_X(\s)/\D_{X_{DR}}(\s)$), thus opening the way for a
   comparison between them (see Theorem \ref{comparison}).
  
\end{rmk}
 
%


\subsection{Shifted Poisson structures}\label{sps}

We are finally ready to define shifted Poisson structures on a derived
Artin stack $X$ lfp over $k$. In the previous Section (Definition
\ref{crys}), we constructed two prestacks $\D_{X_{DR}}$, and
$\mathcal{P}_{X}:= \D_{X/X_{DR}}$ of graded mixed cdga's on $X_{DR}$ ,
together with a map of prestacks $\D_{X_{DR}} \to \mathcal{P}_X$,
exhibiting $\mathcal{P}_X$ as a prestack of $\D_{X_{DR}}$-linear
graded mixed cdga's on $X_{DR}$. By passing to Tate realization, we
obtain the prestack $\mathbf{\underline{Pol}}^{t}(\mathcal{P}_X /
\D_{X_{DR}} ,n)$ of graded $\mathbb{P}_{n+1}$-cdga' on $\dAff/X_{DR}$.

\begin{df}\label{defPoiss}
  If $X$ is a derived Artin stack lfp over $k$, and $n \in \mathbb{Z}$, we define 
\begin{itemize}
\item the graded $\mathbb{P}_{n+1}$-cdga over $k$,
  $\mathsf{Pol}(X,n):= \Gamma (X_{DR},
  \mathbf{\underline{Pol}}^{t}(\mathcal{P}_X / \D_{X_{DR}} ,n) )$ of
  $n$\emph{-shifted polyvectors on} $X$;
\item the space
  $\mathsf{Poiss}(X,n):=\mathsf{Map}_{\dglie^{gr}_k}(k(2)[-1],
  \mathsf{Pol}(X,n+1)[n+1]) $ of $n$\emph{-shifted Poisson structures
    on} $X$.  An $n$\emph{-shifted Poisson structure on} $X$ is an
  element $\pi \in \pi_0 (\mathsf{Poiss}(X,n))$.
\end{itemize}
\end{df}

\

\noindent
In the second item of the previous definition, $\dglie_{k}^{gr}$ is
the $\s$-category of graded $k$-linear dg-Lie algebras, and $k(2)[-1]$
denotes $k$ sitting in cohomological degree $1$, in weight degree $2$,
endowed with trivial bracket and trivial differential.

Since
$\mathbb{L}_X$ is perfect, by Theorem \ref{yeah} (3) and Remark
\ref{natural}, there is an equivalence $$\mathsf{Pol}(X,n) \simeq
Pol(X/k, n) \simeq \bigoplus_{p\geq 0} \Gamma(X,
Sym^{p}(\mathbb{T}_{X}[-n])$$ of graded mixed dg-modules over $k$.
This justifies the use of the word \emph{polyvectors} for
$\mathsf{Pol}(X,n)$.

The intuition behind our definition of $\mathsf{Poiss}(X,n)$ is that
if $X$ is a smooth scheme, $n=0$, and we replace
$\mathsf{Map}_{\dglie^{gr}_k}$ with the usual, strict
$Hom_{\mathsf{dglie}^{gr}_k}$, one gets that an element in
$Hom_{\mathsf{dglie}^{gr}_k}(k(2)[-1], \mathsf{Pol}(X,n+1)[n+1])$ is
exactly a bivector field $\pi$, such that $[\pi,\pi]=0$, $[-,-]$ being
the usual Schouten-Nijenhuis bracket on algebraic polyvector fields on $X$; in
other words, such a $\pi$ is a usual algebraic Poisson bivector on
$X$. See also Example \ref{expoiss} (1) below.

We can give an alternative description of the space
$\mathsf{Poiss}(X,n)$. Recall from the previous section that the
stabilized versions $\D_{X_{DR}}(\infty)$ and
$\mathcal{P}_{X}(\infty)$ of $\D_{X_{DR}}$ and $\mathcal{P}_{X}$ are
both prestacks of commutative monoid objects in $\mathsf{Ind}(\medg)$
on $X_{DR}$, and that there is an analogous canonical map
$\D_{X_{DR}}(\infty) \to \mathcal{P}_{X}(\infty)$. We let
$\mathbb{P}_{n+1}(\mathcal{P}_{X}(\infty)/\D_{X_{DR}}(\infty))$ be the
space of those $\mathbb{P}_{n+1}$-algebras structures on
$\mathcal{P}_{X}(\infty)$, in the $\infty$-category $\C$ of prestacks,
on $X_{DR}$, of $\D_{X_{DR}}(\infty)$-modules inside
$\mathsf{Ind}(\medg)$, which are compatible the \emph{given}
commutative algebra structure on $\mathcal{P}_{X}(\infty)$ in $\C$. An
elaboration of Corollary $\ref{mel}$ yields the following comparison

\begin{thm}\label{melaniext} \emph{(\cite[Theorem 3.1.2]{cptvv})}
For any derived Artin stack $X$ lfp over $k$, and any $n \in
\mathbb{Z}$, we have a canonical equivalence $\mathsf{Poiss}(X,n)
\simeq \mathbb{P}_{\C ,
  n+1}(\mathcal{P}_{X}(\infty)/\D_{X_{DR}}(\infty))$ in $\T$.
\end{thm}
This theorem should be viewed as a vast generalization to derived
Artin stacks of the equivalence between the notion of Poisson
bivectors on $X$ and Poisson brackets on $\mathcal{O}_X$, well-known
for smooth schemes $X$.

\begin{exs}\label{expoiss}
  \emph{(1) If $X$ is a smooth scheme, then the space
    $\mathsf{Poiss}(X,0)$ is discrete and equivalent to the vector
    space of usual algebraic Poisson brackets on $\OO_X$.}

    \

\noindent    
\emph{(2) If $G$ is a reductive group scheme over $k$, and $\mathfrak{g}$ its Lie algebra,
then one has (\cite[3.1]{cptvv})}
$$
\pi_0 (\mathsf{Poiss}(\mathsf{B}G,n))\simeq \begin{cases} 
\wedge^{3}_{k}(\mathfrak{g})^{G} \, , & n=1\\
Sym_{k}^{2}(\mathfrak{g})^{G} \, , & n=2\\
0 \, , & n\neq 1, 2 
        \end{cases}
$$
\end{exs}
The comparison theorem in the next section, together with the
coisotropic (combined with the results of Section \ref{symplex}), and
the intersection theorem in Section \ref{coisotropic}, will provide
more examples of shifted Poisson structures.

\subsection{Comparison between non degenerate shifted
  Poisson structures and shifted symplectic structures} \label{poiss-sympl}

In this Section, we establish a derived analog of the usual
equivalence between classical non-degenerate Poisson structures and
symplectic structures.

Let $X$ be a derived Artin stack lfp over $k$, $n \in \mathbb{Z}$, and
$\pi \in \pi_0 (\mathsf{Poiss}(X,n))$ be an $n$-shifted Poisson
structure on $X$ (Definition \ref{defPoiss}). By considering the
``forgetful'' map $$\mathsf{Map}_{\dglie^{gr}_k}(k(2)[-1],
\mathsf{Pol}(X,n+1)[n+1]) \longrightarrow
\mathsf{Map}_{\dg^{gr}_k}(k(2)[-1], \mathsf{Pol}(X,n+1)[n+1]),$$
together with the equivalence in $\dg^{gr}_{k}$
$$
\mathsf{Pol}(X,n+1)[n+1] \simeq \bigoplus_{p}
\Gamma(X,Sym^{p}_{\OO_X}(\mathbb{T}_{X}[-n-1])[n+1]
$$
(see Remark \ref{natural}), $\pi$ induces a morphism $ k(2) \to
\oplus_{p}\Gamma(X,Sym^{p}_{\OO_X}(\mathbb{T}_{X}[-n-1])[n+2]$ in
$\dg^{gr}_{k}$, and thus defines an element $\alpha_{\pi} \in
H^{-n}(X,\Phi_n^{(2)}(\mathbb{T}_{X})),$ where
$$ 
\Phi_{n}^{(2)} (\mathbb{T}_{X}) := \left \{
  \begin{array}{c}
    Sym^{2}_{\OO_{X}}\mathbb{T}_{X}, \,\, \mathrm{if \; n\; is \; odd} \\
    \wedge^{2}_{\OO_{X}}\mathbb{T}_{X},\,\,  \mathrm{if \; n\; is \; even.}
  \end{array}
\right. 
$$ We denote by $\pi^{\sharp}$ the map $\mathbb{L}_{X} \to
\mathbb{T}_{X}[-n]$ induced, via adjunction, by $\alpha_{\pi}$.

\begin{df} Let $X$ be a derived Artin stack lfp over $k$, and $n \in \mathbb{Z}$.
  An $n$-shifted Poisson structure $\pi \in \pi_0
  (\mathsf{Poiss}(X,n))$ is \emph{non-degenerate} if the induced map
  $\pi^{\sharp}: \mathbb{L}_{X} \to \mathbb{T}_{X}[-n]$ is an
  equivalence. We denote by $\mathsf{Poiss}^{nd}(X,n)$ the subspace of
  $\mathsf{Poiss}(X,n)$ whose connected components are non-degenerate
  $n$-shifted Poisson structures on $X$.
\end{df}


We are now ready to construct the comparison map between the space of
shifted Poisson structures and the space of shifted symplectic
structures on derived Artin stacks.\\ Recalling from the previous
Section, let $\C_{\textrm{Ind}}'$ is the category of prestacks on
$\dAff/X_{DR}$ with values in $\mathsf{Ind}(\epsilon-\dg^{gr}_{k})$,
and $\C_{\textrm{Ind}}$ the category of $D_{X_{DR}}(\infty)$-modules
in $\C'$. Then $A:=\mathcal{P}_{X}(\infty)$ belongs to
$\mathbf{CAlg}_{\C}$. By Proposition \ref{W'} and Th. \ref{melaniext},
we have an equivalence $\mathbb{P}^{nd}_{\C,
  n+1}(\mathcal{P}_{X}(\infty)/\D_{X_{DR}}(\infty)) \simeq
\mathsf{Poiss}^{nd}(X,n)$, where
$\mathbb{P}^{nd}_{n+1}(\mathcal{P}_{X}(\infty)/\D_{X_{DR}}(\infty))$
is the subspace of non degenerate $\mathbb{P}_{n+1}$-algebra
structures on $\mathcal{P}_{X}(\infty)/\D_{X_{DR}}(\infty)$,
compatible with the underlying commutative
$\D_{X_{DR}}(\infty)$-algebra structure on
$\mathcal{P}_{X}(\infty)$.

Recall from Proposition \ref{W}, that for any $\C$, and any $A\in
\mathbf{CAlg}_{\C}$ with a dualizable cotangent complex, we have a map
$\mathbb{P}^{nd}_{\C , n+1}(A) \to \mathsf{Sympl}(A,n)$, from the
moduli space of those non-degenerate $\mathbb{P}_{n+1}$-algebra
structures in $\C$ on $A$ which are compatible with the given
commutative algebra structure on $A$, to the moduli space of
$n$-shifted symplectic structures on $A$. With our current choice of
$\C$, we thus get a map $$\mathsf{Poiss}^{nd}(X,n) \simeq
\mathbb{P}^{nd}_{\C, n+1}(\mathcal{P}_{X}(\infty)/\D_{X_{DR}}(\infty))
\longrightarrow
\mathsf{Sympl}(\mathcal{P}_{X}(\infty)/\D_{X_{DR}}(\infty),n).$$ But
by Theorem \ref{yeah}, we have
$\mathsf{Sympl}(\mathcal{P}_{X}(\infty)/\D_{X_{DR}}(\infty),n) \simeq
\mathsf{Sympl}(X,n)$, so we obtain a comparison map $$\psi:
\mathsf{Poiss}^{nd}(X,n) \longrightarrow \mathsf{Sympl}(X,n)$$ from
the space of non-degenerate $n$-shifted Poisson structures on $X$, to
the space of $n$-shifted symplectic structures on $X$. One of the main
result of \cite{cptvv}, and the key to quantize all the shifted
symplectic moduli spaces constructed in Section \ref{symplex}, is the
following

\begin{thm}\label{comparison} \emph{\cite[Theorem 3.2.4]{cptvv}}
  Let $X$ be a derived Artin stack lfp over $k$, and $n \in
  \mathbb{Z}$. The canonical map $\psi : \mathsf{Poiss}^{nd}(X,n) \to
  \mathsf{Sympl}(X,n)$ is an equivalence in $\T$.
\end{thm}

\

In spite of being expected, and somehow very natural, Theorem
\ref{comparison} has a rather technical and non-trivial proof. It is
not too difficult showing that $\psi$ induces isomorphisms on all the
homotopy groups $\pi_i$'s for $i \geq 1$. More difficult is proving
that $\psi$ is an isomorphism on $\pi_0$: this is achieved in
\cite{cptvv} by first showing that the functors $A\mapsto
\mathbb{P}^{nd}_{\C,
  n+1}(\mathcal{P}_{X}(\infty)(A)/\D_{X_{DR}}(\infty)(A))$, and $A
\mapsto
\mathsf{Sympl}(\mathcal{P}_{X}(\infty)(A)/\D_{X_{DR}}(\infty)(A),n)$
are both formal derived stacks (Definition \ref{d10-}), and then
showing that it is enough to prove the isomorphism on
$\dAff^{\red}/X_{\red}$. This reduced case is then specifically
handled by using pairings and copairings on $L_{\s}$-algebras, the
$L_{\s}$-algebra being given, for any $u: \Spec \, A \to X$, with $A$
reduced, by an $L_{\s}$-model for $(u^*\mathbb{L}_X)^{\vee}[-1]$. Similar techniques have been used in \cite{cogw}.

\subsection{Coisotropic structures}\label{coisotropic}

In this section we discuss briefly the notion of a coisotrpic
structure on a map to a general $n$-shifted Poisson target. This is
analogous to the notion of a lagrangian structure from
Definition~\ref{lagr}. A new feature of the Poisson context is that
the definition of coisotropic structure itself requires a non-trivial
statement - Rozenblyum's \emph{additivity theorem}.  This additivity
theorem asserts that for every $n\geq 1$ and every symmetric monoidal
$\infty$-category $\C$, satisfying our standard hypotheses, there is a
natural equivalence between the $\infty$-category $\mathbb{P}_{\C,
  n+1}-\mathbf{Alg}$ and the $\infty$-category
$\mathbf{Alg}(\mathbb{P}_{\C, n}-\mathbf{CAlg})$ of unital and
associative algebras in the category $\mathbb{P}_{\C,
  n}-\mathbf{Alg}$ \footnote{The same result holds for the operads $\mathbb{E}_n$ of little $n$-disks, and it is known as Dunn-Lurie additivity  \cite[5.1.2.2]{lualg}.}(note that $\mathbb{P}_{\C, n}$ is a Hopf operad,
hence $\mathbb{P}_{\C, n}-\mathbf{Alg}$ inherits a natural symmetric
monoidal structure).

The additivity equivalence is functorial in $\C$,
with respect to symmetric monoidal $\s$-functors, and commutes with
the forgetful functors to $\C$ (see \cite[Remark 3.4.2]{cptvv} for more details). 
The main utility of the additivity theorem is that it allows us to
make sense of $\mathbb{P}_{n+1}$-structures on morphisms between
commutative algebras in $\C$.  More precisely, if we write
$\mathbb{P}_{\C, (n+1,n)}-\mathbf{Alg}$ for the $\s$-category of
pairs $(A,B)$ consisting of an object $A\in
\mathbf{Alg}(\mathbb{P}_{\C, n}-\mathbf{Alg})$ and an $A$-module $B$
in $\mathbb{P}_{\C, n}-\mathbf{Alg}_{\C}$ , then by the additivity
theorem $\mathbb{P}_{\C, (n+1,n)}-\mathbf{Alg}$ comes equipped with
two forgetful $\infty$-functors $(A,B) \mapsto A$, and $(A,B) \mapsto
B$ to $\mathbb{P}_{\C, n+1}-\mathbf{Alg}$, and $\mathbb{P}_{\C,
  n}-\mathbf{Alg}$, respectively. Furthermore, the forgetful functor
$\mathbb{P}_{\C, n}-\mathbf{Alg} \to \mathbf{Alg}_{\C}$ induces a
natural forgetful functor from $\mathbb{P}_{\C,
  (n+1,n)}-\mathbf{Alg}_{\C}$ to the $\s$-category of pairs $(A,B)$
where $A \in \mathbf{Alg}(\mathbf{CAlg}_{\C})$ and $B$ is an
$A$-module in $\mathbf{CAlg}_{\C}$. The $\s$-category of such pairs is
naturally equivalent to the $\s$-category
$\mathbf{Mor}(\mathbf{CAlg}_{\C})$ of morphisms in
$\mathbf{CAlg}_{\C}$. In particular, given a morphism $\varphi : A \to
B$ between commutative algebras in $\C$, we can define the space of
$\mathbb{P}_{\C, (n+1,n)}$-structures on $\varphi$ as the fiber of the
$\infty$-functor $\mathbb{P}_{\C, (n+1,n)}-\mathbf{Alg} \to
\mathbf{Mor}(\mathbf{Alg}_{\C})$ over $\varphi$. We will write
\[
\mathbb{P}_{\C, (n+1,n)}(\varphi) := \mathbb{P}_{\C, (n+1,n)}-\mathbf{Alg}
\times_{\mathbf{Mor}(\mathbf{CAlg}_{\C})} \{ \varphi \} 
\]
for this space. Note that by construction the space $\mathbb{P}_{\C,
  (n+1,n)}(\varphi)$ projects naturally both to the space
$\mathbb{P}_{\C, n+1}(A)$ of $\mathbb{P}_{n+1}$-structures on the
source $A$, and to the space $\mathbb{P}_{\C, n}(B)$ of
$\mathbb{P}_{n}$-structures on the target $B$.

Let $f : X \to Y$ is a morphism of derived stacks locally of finite
presentation over $k$. We specialize the construction above to the
following case: 
\begin{itemize}
\item $\C:= \C_{X, \textrm{Ind}}$ is the $\s$-category defined in Section
  \ref{poiss-sympl} 
(i.e. if $\C_{X, \textrm{Ind}}'$ is the $\s$-category of prestacks on
$\dAff/X_{DR}$ with values in $\mathsf{Ind}(\epsilon-\dg^{gr}_{k})$,
then $\C_{X, \textrm{Ind}}$ is the $\s$-category of $D_{X_{DR}}(\infty)$-modules
in $\C_{X, \textrm{Ind}}'$).
\item $\varphi$ is the induced map $f_{\B}^* : f_{DR}^*(\B_Y(\s)) \to \B_{X}(\s)$.
\end{itemize} 
Note that the map $f^*_{DR}(\D_{Y_{DR}}(\s)) \to \D_{X_{DR}}(\s)$ is
an equivalence and so $f_{\B}^*$ 
may indeed be  considered as a morphism in $\mathbf{CAlg}_{\C}$. 

Now, if $Y$ is endowed with an $n$-shifted Poisson structure $\pi$,
then $\B_{Y}(\s)$ is canonically an object in $\mathbb{P}_{\C_{Y,
    \mathrm{Ind}}, n+1}-\mathbf{Alg}$ (Theorem \ref{melaniext}), and thus
its pull-back $f_{DR}^*(\B_Y(\s))$ is canonically an object in
$\mathbb{P}_{\C, n+1}-\mathbf{Alg}$. We denote this object by
$\B^{\pi}_{Y, f}$ (in order to distinguish it from
$f_{DR}^*(\B_Y(\s))$ as an object in $\mathbf{CAlg}_{\C}$). To ease
notation, we will write $\mathbb{P}_{(n+1,n)}(f_{\B}^*)$ for
$\mathbb{P}_{\C, (n+1,n)}(f_{\B}^*)$, and
$\mathbb{P}_{n+1}(f_{DR}^*(\B_Y(\s)))$ for $\mathbb{P}_{\C,
  n+1}(f_{DR}^*(\B_Y(\s)))$.  So it makes sense to consider the fiber
of the projection map $\mathbb{P}_{(n+1,n)}(f_{\B}^*) \to
\mathbb{P}_{n+1}(f_{DR}^*(\B_Y(\s))) $ over $\B^{\pi}_{Y, f}$.

\begin{df} \label{df-coisotropic} Let $f : X \to Y$ be a morphism of
  derived stacks locally of finite presentation over $k$ and assume
  that $Y$ is equipped with an $n$-shifted Poisson structure $\pi$.
The 
{\em    space of coisotropic structures on $f$ relative to $\pi$} is the
  fiber
\[
\mathsf{Cois}(f,\pi) :=\mathbb{P}_{(n+1,n)}(f_{\B}^*)
\times_{\mathbb{P}_{n+1}(f_{DR}^*(\B_Y(\s)))}
\{ \B^{\pi}_{Y, f} \}.
\]
A \emph{coisotropic structure on} $f$ 
\emph{relative to} $\pi$ is an element in $\pi_0 \, 
\mathsf{Cois}(f,\pi)$.
\end{df}

\

In other words, a coisotropic structure on $f:X \to Y$ consists of the
datum of a $\D_{X_{DR}}(\s)$-linear $\mathbb{P}_{n}$-algebra structure
on $\B_{X}(\s)$ (the target of $f_{\B}^*$), together with a suitably
compatible structure of module over $\B^{\pi}_{Y, f}$ (the source of
$f_{\B}^*$), inside the $\s$-category of $\D_{X_{DR}}(\s)$-linear
graded mixed $\mathbb{P}_{n}$-algebras on $X_{DR}$.

\begin{rmk}\label{rmk-cois}
  This notion of coisotropic structure has the expected geometric behavior:
\begin{itemize}
\item[(i)] Using the other projection map
  $\mathbb{P}_{(n+1,n)}(f_{\B}^{*}) \to \mathbb{P}_{n}(\B_{X}(\s))$
  (i.e. the map keeping only the target of $f_{\B}^*$ ), and
  Theorem \ref{melaniext}, we get that a choice of a coisotropic
  structure on $f : X \to Y$ relative to an $n$-shifted Poisson
  structure $\pi$ on $Y$, tautologically induces an $(n-1)$-shifted
  Poisson structure on the source $X$.
\item[(ii)] If $\pi_{\omega}$ is a non-degenerate Poisson structure
  corresponding to an $n$-shifted symplectic structure $\omega$ on
  $Y$, then, keeping the notations of Definition \ref{df-coisotropic}, one
  expects a natural equivalence of spaces $\mathsf{Lag}(f,\omega)
  \cong \mathsf{Cois}(f,\pi)^{\text{nd}}$ between the space of
  lagrangian structures on $f$ (see Definition~\ref{lagr}) and the
  space of suitably non-degenerate coisotropic structures on $f$. This
  is being investigated by Melani and Safronov (see
  \cite{melani-cois}).
\item[(iii)] The Lagrangian intersection theorem, Theorem~\ref{Lagr}
  was recently extended to the Poisson context in
  \cite{melani-cois}. Suppose $(Y,\pi)$ be an $n$-shifted Poisson
  Artin stack locally of finite presentation over $k$.  Let $f_{i} :
  X_{i} \to Y$, $i = 1, 2$ be maps of derived Artin stacks each
  endowed with coisotropic structures relative to $\pi$. Then, Melani
  and Safronov prove that the derived fiber product $X_{1}\times_{Y}
  X_{2}$ has a natural, induced $(n-1)$-shifted Poisson structure such
  that the natural map $X_{1}\times_{Y} X_{2} \to X_{1}\times X_{2}$
  is a morphism of $(n-1)$-shifted Poisson stacks, where in the target
  $X_{2}$ is endowed with the $(n-1)$-shifted Poisson structure from
  point (i) above, and $X_{1}$ with the corresponding opposite
  $(n-1)$-shifted Poisson structure (see \cite{melani-cois} for
  details).  A classical, i.e. $0$-shifted, and purely cohomological
  precursor of this result was proved in \cite{gb}.  Aside from its
  conceptual significance, the coisotropic intersection theorem of
  \cite{melani-cois} has many purely utilitarian corollaries. It
  allows us to extend the list of examples at the end of
  Section~\ref{sps}, by providing many more examples of shifted
  Poisson structures on moduli stacks, hence of moduli stacks
  admitting natural deformation quantizations (see Section
  \ref{quantization}). For instance, recently Spaide \cite{spaide}
  applied coisotropic inersections to construct and characterize
  shifted Poisson structures on moduli spaces of framed sheaves in
  arbitrary dimension as well as on the moduli of monopoles.
\end{itemize}
\end{rmk}

\section{Deformation quantization} \label{quantization}

Recall that for an ordinary smooth scheme $X$ over $k$, a classical
(unshifted) Poisson structure $\pi$ can be viewed as the infinitesimal
to the deformation of $\mathcal{O}_{X}$ as a sheaf of associative
algebras on $X$.  According to the algebraic deformation quantization
results of Kontsevich \cite{ko} and Yekutieli \cite{yekutieli} every
ordinary Poisson scheme $(X,\pi)$ admits a quantization as a stack of
algebroids. That is we can always find a stack of algebroids
$\mathcal{X}$ defined over $k[[\hbar]]$ with $(\mathcal{X} \text{ mod
} \hbar) = X$ and with infinitesimal $\pi$. Moreover
\cite{ko,yekutieli} all possible quantizations with a given
infinitesimal depend on a choice of a \emph{formality quasi-isomorphism} 
(Drinfeld associator) and are classified by deformation of $(X,\pi)$
as a Poisson scheme over $k[[\hbar]]$. In particular the
\emph{trivial} Poisson deformation corresponding to the
$k[[\hbar]]$-linear Poisson bivector $\hbar\cdot \pi$ gives rise to a
preferred quantization of $(X,\pi)$. This preferred quantization is
Kontsevich's canonical quantization, or in the case of a
non-degenerate $\pi$ is the algebraic Fedosov canonical quantization of
Bezrukavnikov-Kaledin \cite{bezka}.

In this section we discuss the extension of the deformation
quantization problem to shifted Poisson structures on derived Artin
stacks. We argue that the canonical $n$-shifted quantization always
exists as long as $n \neq 0$ and again depends on the choice of a
Drinfeld associator. Ineterestingly enough the special case when $n=0$
remains the hardest case and the best existing quantization results
are still those of \cite{ko,yekutieli}. The natural question of
extending the \cite{ko,yekutieli} quantization of smooth Poisson
schemes to $0$-shifted Poisson derived Artin stacks requires new ideas
and will not be treated here.

\subsection{Weak and strong quantization} \label{weak_and_strong}
Informally, shifted Poisson structures arise when we study
deformations of $X$ in which we allow only partial non-commutativity
in the deformed product structure.  More precisely, an $n$-shifted
Poisson structure can be viewed as the infinitesimal for deforming the
commutative ($= \mathbb{E}_{\infty}$) algebra structure on
$\mathcal{P}_{X}(\s)$ to an $\mathbb{E}_{n+1}$-algebra structure.

To spell this out, recall that for $n \geq 1$ the operad $\mathbb{E}_{n}$ of chains little
$n$-dimensional disks is a $k$-linear dg-operad which is given by the
chain complexes $C_{\bullet}(FM_n,k)$ of the Fulton-MacPherson's
topological operad $FM_{n}$. By definition, the space of operations of
$FM_{n}$ labeled by a finite set $I$ is the Fulton-MacPherson
compactified configuration space $FM_{n}(I)$ of $I$-labeled
configurations of points in $\mathbb{R}^{n}$. For $n \geq 2$ the
Postnikov tower of the spaces $FM_{n}(I)$ induces a filtration of
$\mathbb{E}_n$ whose associated graded is the graded $k$-linear
homology operad $H_{\bullet}(FM_n,k)$ of $FM_{n}$ which is known to be
the operad $\mathbb{P}_n$ controlling $(n-1)$-shifted Poisson
algebras. We can now apply the Rees construction to the filtration on
$\mathbb{E}_n$ to obtain a dg-operad $\mathbb{BD}_{n}$ (for Beilinson-Drinfeld) which is linear
over $k[h]$ and deforms the filtered operad $\mathbb{E}_n$ to its
associated graded $\mathbb{P}_n$. This deformation makes sense for
$n=1$ as well. In this case, $\mathbb{E}_{1}$ is the operad controlling
associative algebra structures. For every finite $I$, the $k$-module
of operations $\mathbb{E}_{1}(I)$ labeled by $I$ is the
non-commutative polynomial algebra $k\langle x_{i} | i \in I \rangle$
which is naturally filtered by monomial degree. The associated graded
to this filtration is the commutative polynomial algebra $k[x_{i} | i
  \in I ]$ equipped with the induced Lie bracket, i.e. we have
$\text{gr} \mathbb{E}_{1}(I) = P_{1}(I)$. Again applying the Rees
construction to the filtration gives a $k[h]$-linear operad
$\mathbb{BD}_{1}$ which interpolates between $\mathbb{E}_{1}$ and
$\mathbb{P}_{1}$.  The difference between this case and the case $n
\geq 2$ is that $\mathbb{P}_{1}$ is \emph{not} the homology of
$\mathbb{E}_{1}$. In fact $\mathbb{E}_{1}$ is already formal and
isomorphic to $H_{\bullet}(FM_{1},k)$.

Nevertheless, for any $n\geq1$ we constructed a $k[\hbar]$-linear
dg-operad operad $\mathbb{BD}_n$ such that
$\mathbb{BD}_n\otimes_{k[\hbar]}k\cong\mathbb{P}_n$ and
$\mathbb{BD}_n\otimes_{k[\hbar]}k[\hbar,\hbar^{-1}]\cong
\mathbb{E}_n[\hbar,\hbar^{-1}]$. With this in place, we are now ready
to formulate the quantization problem in the shifted setting.  Suppose
$X$ is a derived Artin stack, locally of finite presentation over $k$. We use again the notations from Section 2.2.2:
$\C'_{\textrm{Ind}}$ is the $\s$-category of prestacks on $\dAff_{k}/X_{DR}$ with values in
$\mathrm{Ind}(\epsilon-\mathbf{dg}^{gr}_{k})$, so that $\D_{X_{DR}}(\s) \in \mathbf{CAlg}(\C'_{\textrm{Ind}})$, 
$\mathcal{P}_X(\s) \in
\D_{X_{DR}}(\s)/\mathbf{CAlg}(\C'_{\textrm{Ind}})$, and  we define 
$\C_{\textrm{Ind}}:=
\D_{X_{DR}}(\s)-\mathbf{Mod}_{\C'_{\textrm{Ind}}}$.

By Theorem~\ref{melaniext} specifying an $n$-shifted Poisson structure
$\pi$ on $X$ is equivalent to specifying a
$\mathbb{P}_{\C_{\textrm{Ind}}, n+1}$-algebra structure on the Tate stack of principal parts
$\B_{X}(\infty)$, compatible with its  given commutative $\D_{X_{DR}}(\s)$-algebra structure.
Suppose $n \geq 0$, then we have two
flavors of the quantization problem:
\begin{description}
\item[(strong quantization)] Construct a $\D_{X_{DR}}(\infty)$-linear
  $\mathbb{BD}_{n+1}$-algebra structure on
  $\B_{X}(\infty)\otimes k[\hbar]$, such that after tensoring
  with $\otimes_{k[\hbar]} k$ we get the $\mathbb{P}_{n+1}$-structure
  given by $\pi$.
\item[(weak quantization)] Construct a 
$\mathbb{BD}_{n}$-monoidal structure on the $\s$-category  (Theorem \ref{yeah}, (1) and (4))
  \[
  \mathbf{Perf}(X)\otimes_{k} k[\hbar]  \cong
  \mathcal{P}_X(\s)-\mathbf{Mod}_{\C_{\textrm{Ind}}}^{\mathsf{perf}}\otimes_{k} k[\hbar]
  \]
which after $\otimes_{k[\hbar]} k$ recovers the
$\mathbb{P}_{n}$-monoidal structure on $\mathcal{P}_X(\s)-\mathbf{Mod}_{\C_{\textrm{Ind}}}^{\mathsf{perf}}$ corresponding
to $\pi$ via Rozenblyum's additivity theorem.
\end{description}

\begin{rmk} \label{rm-sgivesw}
It is natural to expect that a solution to the strong quantization
problem yields a solution to the weak quantization problem by passing
to the category of perfect complexes over the
$\mathbb{BD}_{n+1}$-algebra provided by the strong quantization. For
this to make sense we need a $\mathbb{BD}$-version of the additivity
theorem. In other words we need to know that, for any $k$-linear presentable stable symmetric monoidal $\s$-category $\mathcal{N}$, there exists a natural
equivalence of $\infty$-categories:
\begin{equation} \label{eq-BDadditivity}
  \mathbb{BD}_{n+1}-\mathbf{Alg}_{\mathcal{N}} \cong
  \mathbf{Alg}\left(\mathbb{BD}_{n}-\mathbf{Alg}_{\mathcal{N}}\right)
\end{equation}
which specializes to Rozenblyum's additivity at $\hbar = 0$ and to
Dunn-Lurie's additivity from \cite[5.1.2.2]{lualg} at $\hbar =
1$. Rozenblyum recently proved that the additivity equivalence
\eqref{eq-BDadditivity} exists and so to any the strong shifted quantization we can indeed associate a weak shifted quantization.
\end{rmk}

Our main result in this setting is the following unobstructedness theorem

\begin{thm}\emph{(\cite[Theorem~3.5.4]{cptvv})} \label{thm-unobs}
  Let $X$ be a derived Artin stack locally of finite presentation over
  $k$, equipped with an $n$-shifted Poisson structure $\pi$.  If $n>0$
  then there is a canonical strong 
  quantization.
\end{thm}

\

\noindent
The above theorem is analogous to the existence of canonical quantization for
unshifted smooth schemes.  In fact, at this stage, the proof of the
theorem is almost a tautology.  Since the operad $\mathbb{E}_{n+1}$ is
formal and for $n>0$ its homology is $\mathbb{P}_{n+1}$, we can choose
a \emph{formality} equivalence of $k$-dg-operads $\alpha_{n+1} :
\mathbb{E}_{n+1} \simeq \mathbb{P}_{n+1}$.  The map $\alpha_{n+1}$
induces an equivalence
$\mathbb{BD}_{n+1}\simeq\mathbb{P}_{n+1}\otimes_kk[\hbar]$ which is
the identity mod $\hbar$.
Therefore one can consider $\B_{X}(\s)\otimes_kk[\hbar]$ as a stack of
$\D_{X_{DR}}(\s)$-linear graded mixed $\mathbb{BD}_{n+1}$-algebras on
$X_{DR}$, and by construction this stack is a strong deformation quantization
of $\B_{X}(\s)$.\\


By specializing the $\mathbb{BD}_{n+1}$-algebra
structure at $\hbar = 1$ we can then view $\B_{X}(\s)$ as a $\mathbb{E}_{n+1}$-algebra in $\C_{\textrm{Ind}}$. By  \cite[5.1.2.2 and 5.1.2.7]{lualg},  the $\s$-category $\mathcal{P}_X(\s)-\mathbf{Mod}_{\C_{\textrm{Ind}}}$ has an induced $\mathbb{E}_{n}$-monoidal structure, and one checks that its full sub-category $\mathcal{P}_X(\s)-\mathbf{Mod}_{\C_{\textrm{Ind}}}^{\mathsf{perf}}$ inherits an $\mathbb{E}_{n}$-monoidal structure that  we will denote by $(\mathcal{P}_X(\s)-\mathbf{Mod}_{\C_{\textrm{Ind}}}^{\mathsf{perf}})_{\mathbb{E}_{n},\, \pi}$. The
subscript $\pi$ indicates that the
$\mathbb{E}_{n+1}$-algebra structure on $\B_{X}(\s)$, hence the induced $\mathbb{E}_{n}$-monoidal structure on $\mathcal{P}_X(\s)-\mathbf{Mod}_{\C_{\textrm{Ind}}}^{\mathsf{perf}}$, depends on $\pi$, while the
subscript $\mathbb{E}_{n}$ records the $\mathbb{E}_{n}$-monoidal structure. This is exactly the
deformation of $\mathbf{Perf}(X) \cong \mathcal{P}_X(\s)-\mathbf{Mod}_{\C_{\textrm{Ind}}}^{\mathsf{perf}}$ that we were looking for. We
record this fact in the following

\begin{df}\label{dquant}
  With the notation above, and $n>0$, \emph{the weak quantization of
    $X$ with infinitesimal $\pi$} is the $\mathbb{E}_{n}$-monoidal
  $\s$-category
  $$
  \mathbf{Perf}(X,\pi):=(\mathcal{P}_X(\s)-\mathbf{Mod}_{\C_{\textrm{Ind}}}^{\mathsf{perf}})_{\mathbb{E}_{n},\, \pi}  $$
\end{df}




\

\noindent {\bfseries Quantization for $n<0$.} \ 
The quantization problem  for $n$-shifted Poisson structures with $n <
0$ can be understood by looking at parameter spaces which are themselves
dg schemes. 
%
%
Concretely, let $n <
0$ and let $\pi$ be an $n$-shifted Poisson structure on some derived
stack $X$ lfp over $k$.  Let $\hbar_{2n}$ a formal
variable of cohomological degree $2n$, and consider the stack
$\B_{X}(\s)[\hbar_{2n}]$ of Ind-objects in graded
$k(\s)[\hbar_{2n}]$-linear mixed cdgas over $X_{DR}$. Because of the
homological shift it is equipped with a $k(\s)[\hbar_{2n}]$-linear
$\mathbb{P}_{1-n}$-structure, induced by $\hbar_{2n}\cdot \pi$ (Theorem~\ref{melaniext}). Since
$n<0$, this brings us back to the situation of positively shifted
Poisson structures.

Proceeding as before, we choose a formality equivalence of
$k$-dg-operads $\alpha_{1-n} : \mathbb{E}_{1-n} \simeq
\mathbb{P}_{1-n}$, and thus view $\B_{X}(\s)[\hbar_{2n}]$ as a
$k(\s)[\hbar_{2n}]$-linear $\mathbb{E}_{1-n}$-algebra. Again by using
Dunn-Lurie's additivity \cite[5.1.2.2 and 5.1.2.7]{lualg}, the
$\s$-category $\mathcal{P}_X(\s)-\mathbf{Mod}_{\C_{\textrm{Ind}}}^{\mathsf{perf}}$
comes equipped with an induced 
$\mathbb{E}_{-n}$-monoidal structure.  We will write
$(\mathcal{P}_X(\s)-\mathbf{Mod}_{\C_{\textrm{Ind}}}^{\mathsf{perf}})_{\mathbb{E}_{-n},\,\pi}$
for this $\mathbb{E}_{-n}$-monoidal category. Thus for $n < 0$ we can define 
 the \emph{weak quantization of $X$}
    with infinitesimal $\pi$ as the $\mathbb{E}_{-n}$-monoidal
  $\s$-category
$$
\mathbf{Perf}(X,\pi):=(\mathcal{P}_X(\s)-\mathbf{Mod}_{\C_{\textrm{Ind}}}^{\mathsf{perf}})_{\mathbb{E}_{-n},\,\pi}.
$$

\

\noindent
As before, the underlying $\s$-category of $\mathbf{Perf}(X,\pi)$ is
$\mathbf{Perf}(X)\otimes_{k}k[\hbar_{2n}]=:\mathbf{Perf}(X)[\hbar_{2n}]$.
Hence for $n < 0$ our weak quantization of $X$ consists then of the
\emph{datum of a $\mathbb{E}_{-n}$-monoidal structure on}
$\mathbf{Perf}(X)[\hbar_{2n}]$, and by the strong version of
Rozenblyum's additivity, such a quantization can be considered as a
deformation of the standard symmetric (i.e. $\mathbb{E}_{\s}$-) monoidal structure on
$\mathbf{Perf}(X)[\hbar_{2n}]$. Note that this standard symmetric
monoidal structure on $\mathbf{Perf}(X)[\hbar_{2n}]$ recovers the
standard symmetric monoidal structure on $\mathbf{Perf}(X)$ after base
change along the canonical map $k[\hbar_{2n}] \rightarrow k$.

\begin{rmk}\label{kapustin} This quantization answers a conjecture of Kapustin \cite[3.2]{kap} which concerns the $n=-1$ case. Note that Kapustin considers $\mathbb{Z}/2$-graded derived categories, and therefore the fact that we work over $k[\hbar_{2n}]$ is immaterial: we really obtain a quantization of the $\mathbb{Z}/2$ perfect derived category of $X$, since $\mathbf{Perf}(X)[\hbar_{2n}]$ and $\mathbf{Perf}(X)$ coincide after $\mathbb{Z}/2$-periodization.
\end{rmk}

\subsection{Examples and vistas}

\

\noindent
{\bfseries (a)} \ {\bfseries Quantization of moduli.} 
The equivalence of shifted symplectic and
non-degenerate shifted Poisson structures from Theorem~\ref{comparison} combined
with the $n > 0$ (or $n < 0$) quantization scheme described in the
previous section provides a canonical $\mathbb{E}_{n}$-monoidal (or
$\mathbb{E}_{-n}$-monoidal) deformation of the $\s$-category of
perfect complexes on the various shifted symplectic  moduli stacks listed
at the end of section~\ref{sps}. For example:
\begin{itemize}
\item For a derived Artin stack $X$ locally of finite presentation we
  obtain a canonical quantization of the shifted cotangent stack
  $T^{*}X[n]$ for $n \neq 0$. The shifted cotangent stack $T^{*}X[n]$
  has a natural $n$-shifted symplectic form \cite{cal2}. If
  we denote the corresponding non-degenerate $n$-shifted Poisson
  structure by $\pi_{n}$, then the modules over the
  $\mathbb{E}_{n}$-monoidal (or $\mathbb{E}_{-n}$-monoidal) category
  $\mathbf{Perf}(T^{*}X[n],\pi_{n})$ will be the modules over the
  $n$-shifted differential operators on $X$.  
 \item  For a complex reductive
  group $G$ we obtain canonical quantizations of:
\begin{itemize}
\item the derived stack $\mathbb{R}\mathsf{Loc}_{G}(M)$ of
$G$-local systems on a compact oriented topological manifold $M$ of
  dimension $\neq 2$;
\item the derived stack $\mathbb{R}\mathsf{Loc}_{G}^{DR}(X)$ of algebraic $G$-local
  systems on a smooth complex projective variety $X$ of dimension $> 1$;
\item the derived stack $\mathbb{R}\mathsf{Higgs}_{G}(X)$ of algebraic
  $G$-Higgs bundles on a smooth complex projective variety $X$ of dimension $>
  1$;
\item the derived stack $\mathbb{R}\mathsf{Bun}_{G}(X)$ of algebraic
  $G$-torsors on a smooth compact Calabi-Yau variety $X$ of dimension
  $\neq 2$.
\end{itemize}

Similarly we get quantizations of the stack of perfect complexes on a
compact oriented topological manifold $M$ of dimension $\neq 2$, of
the derived stack of perfect complexes over $X_{DR}$ for a smooth
complex projective variety $X$ of dimension $> 1$, of the derived 
stack of perfect complexes of Higgs bundles on a smooth complex
projective variety of dimension $> 1$, and on the derived stack of
perfect complexes on a smooth compact Calabi-Yau variety of dimension
$\neq 2$.
\item For a smooth compact Calabi-Yau dg category $T$ of dimension
  $\neq 2$ we get a canonical weak quantization of the derived moduli
  stack $\mathcal{M}_{T}$ of compact objects in $T$. For instance we
  can take $T$ to be the category of graded matrix factorizations of a
  cubic polynomial $f$ in $3n$ variables with $n \neq 2$. Applying the
  general quantization procedure to this setting we get an
  $(n-2)$-shifted quantization of the moduli stack of graded matrix
  factorizations of $f$.
\end{itemize}

\

\smallskip

\noindent
    {\bfseries (b)} \ {\bfseries Quantization formally at a point.} \
    Let $(X,\pi)$ be an $n$-shifted Poisson derived
Artin stack locally of finite presentation and let $x : * = \mathbf{Spec}\, 
k \to X$ be a closed point. It can be checked
\cite[Lemma~3.6.1]{cptvv} that any $n$-shifted Poisson structure on
  $X$ induces an $n$-shifted Poisson structure on the formal
  completion $\widehat{X}_x$ at $x$.

As a (non-mixed) graded cdga over
$k$, $\B_{\widehat{X}_x}$ is equivalent to
$$
Sym(\mathbb{L}_{*/\widehat{X}_x}[-1])\cong
Sym(x^*\mathbb{L}_{\widehat{X}_x})\cong Sym(x^*\mathbb{L}_X)\,.
$$ We therefore get a graded mixed $\mathbb P_{n+1}$-algebra structure
on $Sym(x^*\mathbb{L}_X)$, whose underlying graded mixed cdgas is
$\B_{\widehat{X}_x}$.  After a choice of formality $\alpha_{n+1}$, we
get a graded mixed $\mathbb{E}_{n+1}$-structure on
$Sym(x^*\mathbb{L}_X)$ whenever $n>0$. When $\pi$ is non-degenerate at
$x$ and the induced $n$-shifted Poisson structure on
$\B_{\widehat{X}_x}$ is strict and constant, then the graded mixed
$\mathbb{E}_{n+1}$-structure on $Sym(x^*\mathbb{L}_X)$ can be
described explicitly in terms of Kontsevich's graph complex
\cite[3.6.1]{cptvv}. When the underlying dg-Lie algebra is formal, the explicit formula then identifies the
$\mathbb{E}_{n+1}$-structure on $Sym(x^*\mathbb{L}_X)$ with the Weyl
$n$-algebra recently introduced by Markarian \cite{markarian}.

\

\smallskip

\noindent
{\bfseries (c)} \ {\bfseries Quantization of $BG$.} \ Suppose $G$ is an
affine group scheme, and let $X = BG$ be the classifying stack of $G$.
Note that $X_{DR}=B(G_{DR})$. Let $x:*\to BG$ be the classifying map
of the unit $e:*\to G$. We have a fiber sequence of groups
$\widehat{G}_e\longrightarrow G\longrightarrow G_{DR}$, 
and hence $\widehat{BG}_x\simeq B(\widehat{G}_e)$. 

As we noted in {\bfseries (b)} the pull-back of
$\B_X$ along $x_{DR}:*\to BG_{DR}$ is $\B_{\widehat{X}_x}$.
Thus the symmetric monoidal 
$\infty$-category
$$
\mathbf{Perf}(BG)\simeq \B_X-\mathbf{Mod}_{\epsilon-\mathbf{dg}^{gr}}^{\mathbf{Perf}}
$$
is equivalent to the symmetric monoidal $\infty$-category of
$G_{DR}$-equivariant objects in
$$
\mathbf{Perf}(B\widehat{G}_e)\simeq
\B_{\widehat{X}_x}-\mathbf{Mod}_{\epsilon-\mathbf{dg}^{gr}}^{\mathbf{Perf}}\,. 
$$
In view of this the quantization of an $n$-shifted Poisson
structure on $BG$ will be determined completely by the
$G_{DR}$-equivariant graded mixed $\mathbb{E}_{n+1}$-algebra structure
on $\B_{\widehat{X}_x}$ obtained from the equivalence
$\alpha_{n+1}:\mathbb{P}_{n+1}\simeq\mathbb{E}_{n+1}$.

This algebra structure can be analyzed in concrete terms. Before we
look more closely at the $1$ and $2$-shifted cases it is useful to
observe that as a graded mixed cdga over $k$ the algebra
$\B_{\widehat{X}_x}\simeq\mathbb{D}(B\widehat{G}_e)$ admits an
explicit description. Indeed, in \cite[3.6.2]{cptvv} it is proven that
$\mathbb{D}(B\widehat{G}_e)$ is actually equivalent to the
Chevalley-Eilenberg graded mixed cdga $\mathsf{CE}(\mathfrak{g})$ of
the Lie algebra $\mathfrak{g}= \mathsf{Lie}(G)$.\\

\noindent \textbf{The case $n=1$ for a reductive $G$.}  For a reductive group $G$
  the $1$-shifted Poisson structures on $BG$ are simply elements in
  $\wedge^3(\mathfrak g)^G$.  If $\pi$ is such an element, then the
  induced $1$-shifted Poisson structure on the graded mixed cdga
  $\mathsf{CE}(\mathfrak g)$ is given explicitly as a semi-strict
  $\mathbb{P}_{2}$-structure (see \cite{mel}): all structure
  $2$-shifted polyvectors are trivial except for the $3$-ary one which
  is constant and given by $\pi$.

The weak $1$-shifted deformation quantization in particular gives rise
to a deformation of the category $\mathsf{Rep}^{fd}(\mathfrak{g})$ of finite dimensional representation of $\mathfrak{g}$ as a monoidal
category. For specific choices of $\pi$ we recover familiar monoidal
deformations:

\begin{ex}\label{ex:Z}
Given a non-degenerate invariant pairing $\langle\,,\,\rangle$ on $\mathfrak g$,
we can choose $\pi$ as the dual of the $G$-invariant linear form
$$
\wedge^3\mathfrak g\longrightarrow k,
\quad\quad(x,y,z)\longmapsto \langle x,[y,z]\rangle\,.
$$
Alternatively, any invariant symmetric $2$-tensor $t\in
Sym^2(\mathfrak g)^G$ leads to such an element
$\pi=[t^{1,2},t^{2,3}]\in\wedge^3(\mathfrak g)^G$.  In these cases the
deformation of $\mathsf{Rep}^{fd}(\mathfrak{g})$ as a monoidal
category can be obtained by means of a deformation of the
associativity constraint only (see \cite{Dr1}), which then looks like
$$
\Phi=1^{\otimes3}+\hbar^2\pi+
o(\hbar^2)\in U(\mathfrak{g})^{\otimes3}[[\hbar]]\,.
$$
\end{ex}
\noindent \textbf{The case $n=2$ for a reductive $G$.}  For a reductive group $G$
  the equivalences classes of $2$-shifted Poisson structures on $BG$
  are in bijection with elements $t\in Sym^2(\mathfrak g)^G$.  The
  induced $2$-shifted Poisson structure on the graded mixed cdga
  $\mathsf{CE}(\mathfrak g)$ is strict and constant.  The graded mixed
  $\mathbb{E}_3$-structure on $\mathsf{CE}(\mathfrak g)$ given by our
  deformation quantization then takes the form of a Weyl $3$-algebra,
  as described in \cite{markarian}.

Note that  this graded mixed
$\mathbb{E}_3$-structure is $G_{DR}$-equivariant by construction, so
it leads to an $\mathbb{E}_2$-monoidal deformation of
$\mathsf{Perf}(BG)$.  This in particular leads to a braided monoidal
deformation of $\mathsf{Rep}^{fd}(\mathfrak{g})$.

Such deformation quantizations of $BG$ have already been constructed: 
\begin{itemize}
\item when $\mathfrak g$ is reductive and $t$ is non-degenerate, by
  means of purely algebraic methods: the quantum group
  $U_\hbar(\mathfrak g)$ is an explicit deformation of the enveloping
  algebra $U(\mathfrak g)$ as a quasi-triangular Hopf algebra.
\item without any assumption, by Drinfeld \cite{Dr2}, using
  transcendental methods similar to the ones that are crucial in the
  proof of the formality of $\mathbb{E}_2$.
\end{itemize} 
It is known that Drinfeld's quantization is equivalent to the quantum 
group one in the semi-simple case (see e.g. \cite{kas2} and references
therein).

\begin{rmk}
  It is interesting to note that our quantization, in contrast to Drinfeld's,  relies on the
  formality of $\mathbb{E}_3$ rather than on the formality of
  $\mathbb{E}_2$. 
\end{rmk}

\

\smallskip


\noindent {\bfseries (d)} \ {\bfseries Relative and absolute
  quantization.} \ An important question that is not addressed in this
paper or in \cite{cptvv} is the question of quantizing lagrangian
structures on maps with shifted symplectic targets or quantizing
coisotropic structures on maps with shifted Poisson targets. The
quantization problem in this relative stting can be formulated in a
manner similar to the absolute quantization from
section~\ref{weak_and_strong}. 

Suppose $(Y,\pi)$ is an $n$-shifted derived Artin stack locally of
finite presentation, and $f : X\to Y$ is a morphism of derived stacks
furnished with a coisotropic structure $\kappa$ relative to $\pi$.
The $\s$ categories $\mathbf{Perf}(X)$ and $\mathbf{Perf}(Y)$ are
symmetric monoidal categories and via the pullback functor $f^{*} :
\mathbf{Perf}(Y) \to \mathbf{Perf}(X)$ the category $
\mathbf{Perf}(X)$ becomes a module over $\mathbf{Perf}(Y)$ so that its
$\mathbb{E}_{\infty}$-monoidal structure becomes linear over the
$\mathbb{E}_{\infty}$-monoidal structure on $\mathbf{Perf}(Y)$. In
other words $f^{*}$ makes $\mathbf{Perf}(X)$ is an
$\mathbb{E}_{\infty}$-algebra over $\mathbf{Perf}(Y)$.  Assume for
simplicity $n > 0$. Then the weak quantization of $(Y,\pi)$ gives a
deformation of $\mathbf{Perf}(Y)$ to an $\mathbb{E}_{n}$-monoidal
category $\mathbf{Perf}(Y,\pi)$. The {\bfseries weak quantization
  problem for $f$} is to find a concurrent deformation of
$\mathbf{Perf}(X)$ as an algebra over $\mathbf{Perf}(Y)$. In other
words, we need to deform $\mathbf{Perf}(X)$ to an
$\mathbb{E}_{n-1}$-monoidal category $\mathbf{Perf}(X,\kappa)$, so
that the functor $f^{*}$ deforms to a functor $\mathbf{q}f^{*} :
\mathbf{Perf}(Y,\pi) \to \mathbf{Perf}(X,\kappa)$ exhibiting
$\mathbf{Perf}(X,\kappa)$ as a module, inside
$\mathbb{E}_{n-1}$-monoidal $\s$-categories, over
$\mathbf{Perf}(Y,\pi)$, viewed as an $\mathbb{E}_{1}$-algebra in
$\mathbb{E}_{n-1}$-monoidal categories. Here again we use Dunn-Lurie
additivity \cite[\S 5.1.2]{lualg} asserting the equivalence between
$\mathbb{E}_{1}$-algberas in $\mathbb{E}_{n-1}$-algebras and
$\mathbb{E}_{n}$-algebras, in any base symmetric monoidal
$\s$-category.

Note that the extension from the absolute to the relative case is not
tautological as the \emph{swiss-cheese operad}, which governs the deformations
of pairs of an $\mathbb{E}_{n}$ algebra and an $\mathbb{E}_{n-1}$
module over it, is not formal.  Nevertheless we expect that such
relative quantizations are again unobstructed for $n > 0$ and we are
currently investigating the problem. 

Another interesting problem in this regard is the question of
compatibility of quantizations with our standard constructions. A
simple instance of this goes as follows. Suppose $Y$ is an $n$-shifted
symplectic derived stack (with $n > 1$) and suppose $f _{1} : X_{1}
\to Y$ and $f_{2} : X_{2} \to Y$ be two morphisms equipped with
lagrangian structures. The derived intersection $Z =
X_{1}\times^{h}_{Y} X_{2}$ carries a natural $(n-1)$-shifted
symplectic form. We expect that the derived intersection persists in
quantizations, i.e. that absolute quantization of $Z$ is the homotopy
fiber product of the relative quantizations of $f_{1}$ and
$f_{2}$. More precisely, write $\pi$ for the non-degenerate shifted
Poisson structure corresponding to the symplectic structure on $Y$ and
$\eta$ for the induced $(n-1)$-shifted non-degenerate Poisson
structure on $Z$. Let $\kappa_{2} \in \pi_{0}\mathsf{Cois}(f_{2},\pi)$
be the non-degenerate coisotropic structure relative to $\pi$
corresponding to the lagrangian structure on $f_{2}$, and let
$\kappa_{1} \in \pi_{0}\mathsf{Cois}(f_{1},-\pi)$ be the
non-degenerate coisotropic structure relative to $-\pi$ corresponding
to the lagrangian structure on $f_{1}$. Then
$\mathbf{Perf}(X_{2},\kappa_{2})$ is an $\mathbb{E}_{n-1}$-algebra
over the $\mathbb{E}_{n}$-algebra $\mathbf{Perf}(Y,\pi)$. Similarly
$\mathbf{Perf}(X_{1},\kappa_{1})$ is an $\mathbb{E}_{n-1}$-algebra
over the $\mathbb{E}_{n}$-algebra $\mathbf{Perf}(Y,-\pi)$ or
equivalently $\mathbf{Perf}(X_{1},\kappa_{1})$ is an
$\mathbb{E}_{n-1}$-algebra over the opposite of the
$\mathbb{E}_{n}$-algebra $\mathbf{Perf}(Y,\pi)$. Conjecturally the
quantized $\mathbb{E}_{n-1}$-monoidal category $\mathbf{Perf}(Z,\eta)$
is reconstructed from the $\mathbb{E}_{n-1}$-monoidal category
$\mathbf{Perf}(X_{1},\kappa_{1})\otimes_{\mathbf{Perf}(Y,\pi)}
\mathbf{Perf}(X_{2},\kappa_{2})$.

 \

\smallskip

\noindent {\bfseries (e)} \ {\bfseries Vistas.
} \ We conclude our paper by  short list of few other directions of
investigation, just to stimulate the reader's interest.

First of all we would like to mention the proposal of
\cite[4.2.2]{cal} for a category of \emph{lagrangian correspondences}
$\mathbf{Lagr}_n$, based on Thm. \ref{Lagr}. Roughly speaking its
objects are $n$-shifted symplectic stacks, and morphisms from $X$ to
$Y$ are maps $L \to X\times Y$ equipped with lagrangian structures
(where $X$ is considered with the opposite of the given symplectic
structure). Some truncated versions of this category were already
considered in \cite{amorim}.  The details of a complete construction
of $\mathbf{Lagr}_n$ as an $\s$-category together with its natural
symmetric monoidal structure remain to be written down, but we have no
doubts that such a construction exists. This will be an important
step toward the study of extended TQFT's (as in \cite{luetft}) with
values in lagrangian correspondences. Building on Remark
\ref{rmk-cois}, one can imagine an analogous construction by replacing
shifted symplectic stacks with shifted Poisson stacks, and lagrangian
correspondences with \emph{coisotropic correspondences}. Such a
construction is currently being considered by Melani and Safronov, and
it might shed some light on Weinstein's original proposal \cite{we}.

Another promising research program related to the topics treated in
this review is the one being pursued, since a few years, by D. Joyce
and various collaborators. As part of their project, they use shifted
symplectic structures to study Donaldson-Thomas moduli spaces of
Calabi-Yau fourfolds, vanishing cycles, and various categorifications
of Donaldson-Thomas invariants (see e.g. \cite{bbdjs, jsaf}).

The geometry of coisotropic structures for shifted Poisson structures is in a very early stage of development, and
a lot of new phenomena need to be properly explored. Just to give one example, the identity map is always endowed with a canonical
coisotropic structure, and this produces a map from $n$-shifted Poisson structures to $(n-1)$-shifted Poisson structures, which is worth investigating.
If not trivial, this could e.g. connect the $n=2$ to the $n=1$ examples of quantizations of $BG$ \textbf{(c)} above. More generally, having a definition of coisotropic structures that is equivalent to Definition \ref{df-coisotropic} but avoids reference to the additivity theorem, would be very useful, especially in applications. Some important steps in this direction have been done by P. Safronov \cite{saf}, and more recently by V. Melani and P. Safronov \cite{melani-cois}.

Finally, it would be interesting to have a version of Theorem
\ref{MAP} with target $\mathsf{Perf}$ and a \emph{stratified}
topological space as a source. A possible way to include the
stratification in our theory is through the use of the
MacPherson-Lurie \emph{exit paths} $\s$-category \cite[A.6]{lualg}. If this
can be accomplished, then a corresponding relative version, as in
Remark~\ref{damienboundary}, could be relevant for some ideas and
conjectures about moduli spaces of constructible sheaves with singular
support in a legendrian knot (see \cite{stz}). A promising preliminary
step in this direction is the theory of left and right
Calabi-Yau structures on functors developed recently by To\"en, and Brav-Dyckerhoff (\cite[p. 228]{toenems}, \cite{brav}).

\

\smallskip 

\noindent
Tony Pantev, {\sc University of Pennsylvania}, tpantev@math.upenn.edu

\smallskip

\noindent
Gabriele Vezzosi, {\sc Universit\`a di Firenze},
gabriele.vezzosi@unifi.it

\end{document}